\documentclass[a4paper, 10pt ]{article}

\usepackage[a4paper,left=2cm,right=2cm,top=2.5cm,bottom=2.5cm]{geometry}
\usepackage[backend=biber,style=numeric]{biblatex} 
\usepackage{hyperref}
\hypersetup{colorlinks,linkcolor={blue},citecolor={blue},urlcolor={blue}}
\usepackage{subfig}
\usepackage{pgfplots}
\usepackage{pgfplotstable}
\usepackage{mathtools}
\usepackage{multicol}
\usepackage{comment}
\usepackage{booktabs}
\pgfplotsset{compat=1.5}
\usepackage{amssymb}
\usepackage{url}
\usepackage{bm}
\usepackage{enumitem}
\usepackage{siunitx}
\usepackage{pdfsync}
\usepackage{float}
\usepackage{tabularx}
\usepackage{array}
\usepackage{xspace}
\usepackage{tikz}
\usepackage{tikz-cd}
\usepackage{tikzsymbols}
\usepackage[labelfont=bf]{caption}
\usepackage[normalem]{ulem}
\usetikzlibrary{calc,trees,positioning,arrows,chains,shapes.geometric,%
    decorations.pathreplacing,decorations.pathmorphing,shapes,%
    matrix,shapes.symbols, decorations.markings, patterns,fit}

\usepackage{siunitx}

\usepackage{authblk}

\usepackage[draft,inline,marginclue]{fixme}
\usepackage{mathrsfs}
\usepackage{color, colortbl}
\usepackage{multirow}

\usepackage{amsthm}
\usepackage{amsmath}
\usepackage{stmaryrd}
\usepackage[ruled,vlined,linesnumbered]{algorithm2e}
\usepackage{graphicx}
\usepackage{booktabs}
\usepackage{cleveref}
\usepackage{lineno}

\FXRegisterAuthor{fr}{afr}{\color{blue}FR}

\FXRegisterAuthor{fg}{afg}{\color{orange}FG}

\FXRegisterAuthor{ac}{aac}{\color{purple}AC}

\newcommand\deleted{\bgroup\markoverwith{\textcolor{red}{\rule[0.5ex]{2pt}{0.4pt}}}\ULon}

\newtheorem{theorem}{Theorem}

\newtheorem{corol}{Corollary}

\newtheorem{problem}{Problem}
\theoremstyle{remark}
\newtheorem{rmk}{Remark}

\newcolumntype{C}[1]{>{\centering\arraybackslash}m{#1}}

\newcommand{\dofu}{d_{\mathbf u}}
\newcommand{\dofp}{d_p}
\newcommand{\nvert}{n_{p,\mathcal T}}
\newcommand{\nvertref}{n_{p,\mathcal S}}
\newcommand{\shape}[1]{\mathcal T_{#1}}
\definecolor{Gray}{gray}{0.9}

\newcommand{\x}{\ensuremath{\mathbf{x}}}

\newcommand{\ub}{\ensuremath{\mathbf{u}}}
\newcommand{\wb}{\ensuremath{\mathbf{w}}}

\newcommand{\phib}{\ensuremath{\boldsymbol{\phi}}}

\DeclareMathOperator*{\argmin}{\arg\!\min}

\newcommand{\XS}{X_{\mathcal S}}
\newcommand{\FS}{F_{\mathcal S}}
\newcommand{\lS}{l_{\mathcal S}}
\newcommand{\NS}{{\mathbf N}_{\mathcal S}}
\newcommand{\XT}{X_{\mathcal T}}
\newcommand{\FT}{F_{\mathcal T}}
\newcommand{\lT}{l_{\mathcal T}}

\newcommand\rsvd{\mathrel{\stackrel{\makebox[0pt]{\mbox{\normalfont\tiny rSVD}}}{=}}}
\newcommand{\bu}{\mathbf{u}}
\begin{document}

\title{
Data assimilation performed with robust shape registration and graph neural networks: application to aortic coarctation}

\author[1]{Francesco~Romor\footnote{francesco.romor@wias-berlin.de}}
\author[2]{Felipe Galarce}
\author[3]{Jan Brüning}
\author[3]{Leonid Goubergrits}
\author[1]{Alfonso Caiazzo\footnote{alfonso.caiazzo@wias-berlin.de}}

\affil[1]{Weierstraß Institute, Mohrenstr. 39 10117, Berlin, Germany}
\affil[2]{School of Civil Engineering, Pontificia Universidad Católica de Valparaíso, Valparaíso, Chile.}
\affil[3]{Institute of Computer-assisted Cardiovascular Medicine, Deutsches Herzzentrum der Charité, Augustenburger Platz 1, 13353 Berlin, Germany}

\maketitle

\begin{abstract}
  Image-based, patient-specific modelling of hemodynamics can improve diagnostic capabilities and provide complementary insights to better understand the hemodynamic treatment outcomes. However, computational fluid dynamics simulations remain relatively costly in a clinical context. Moreover, projection-based reduced-order models and purely data-driven surrogate models struggle due to the high variability of anatomical shapes in a population. A possible solution is shape registration: a reference template geometry is designed from a cohort of available geometries, which can then be diffeomorphically mapped onto it. This provides a natural encoding that can be exploited by machine learning architectures and, at the same time, a reference computational domain in which efficient dimension-reduction strategies can be performed. We compare state-of-the-art graph neural network models with recent data assimilation strategies for the prediction of physical quantities and clinically relevant biomarkers in the context of aortic coarctation.
\end{abstract}

\tableofcontents

\section{Introduction}
\label{sec:intro}
Computational fluid dynamics (CFD) of the cardiovascular system can provide valuable tools for analyzing complex biological phenomena in the absence or scarcity of medical data,  supporting the image-based diagnosis and diseases staging, as well as the treatment selection and postoperative monitoring (see, e.g.~\cite{quarteroni24_cardio, lesage2023mapping,Morris18}). 
While computational models are becoming increasingly relevant in clinical practice, computational costs remain a significant bottleneck. In specific cases, such as the onset of turbulent flows in large vessels or the simulation of complex flows in realistic geometries, CFD simulations demand substantial computational resources.

Physics-based reduced-order models (ROMs) (see, e.g.~\cite{benner_model_2017, hesthaven2016certified, rozza2022advanced}) and data-driven models obtained from machine learning have proven to be able to reduce the computational costs in different applications. However, especially in the context of hemodynamics, these models are often highly patient-specific and cannot be efficiently extrapolated to computational domains of new patients without requiring a substantial amount of new training simulation data, thereby hindering the development of effective predictive models: they may be employed to monitor post-surgery biomarkers, but generally they could not be extended to inter-patient studies. Often the shape variability is introduced locally via small parametrized deformations, and on a fixed computational domain~\cite{Ballarin2016}. Projection-based ROMs are applied for CFD simulations with parametric inflows or optimal control problems~\cite{Girfoglio2022, Zainib2019} on fixed geometries, and they cannot be extrapolated to different patients without recomputing all the training snapshots. Purely data-driven surrogates are also often designed on single geometries~\cite{Fresca2020} and are affected by the same limitations. All these methodologies could benefit from the registration techniques we are going to introduce.

In this work, we present a non-parametric shape registration approach to handle the geometric variability in the case of realistic aortic coarctation and
use it to design a data assimilation pipeline to efficiently predict the three-dimensional patient-specific velocity and pressure fields from velocity measurements or from a patient-specific geometrical encoding
and additional boundary conditions.

\begin{figure}[!htp]
  \centering
  \includegraphics[width=1\textwidth]{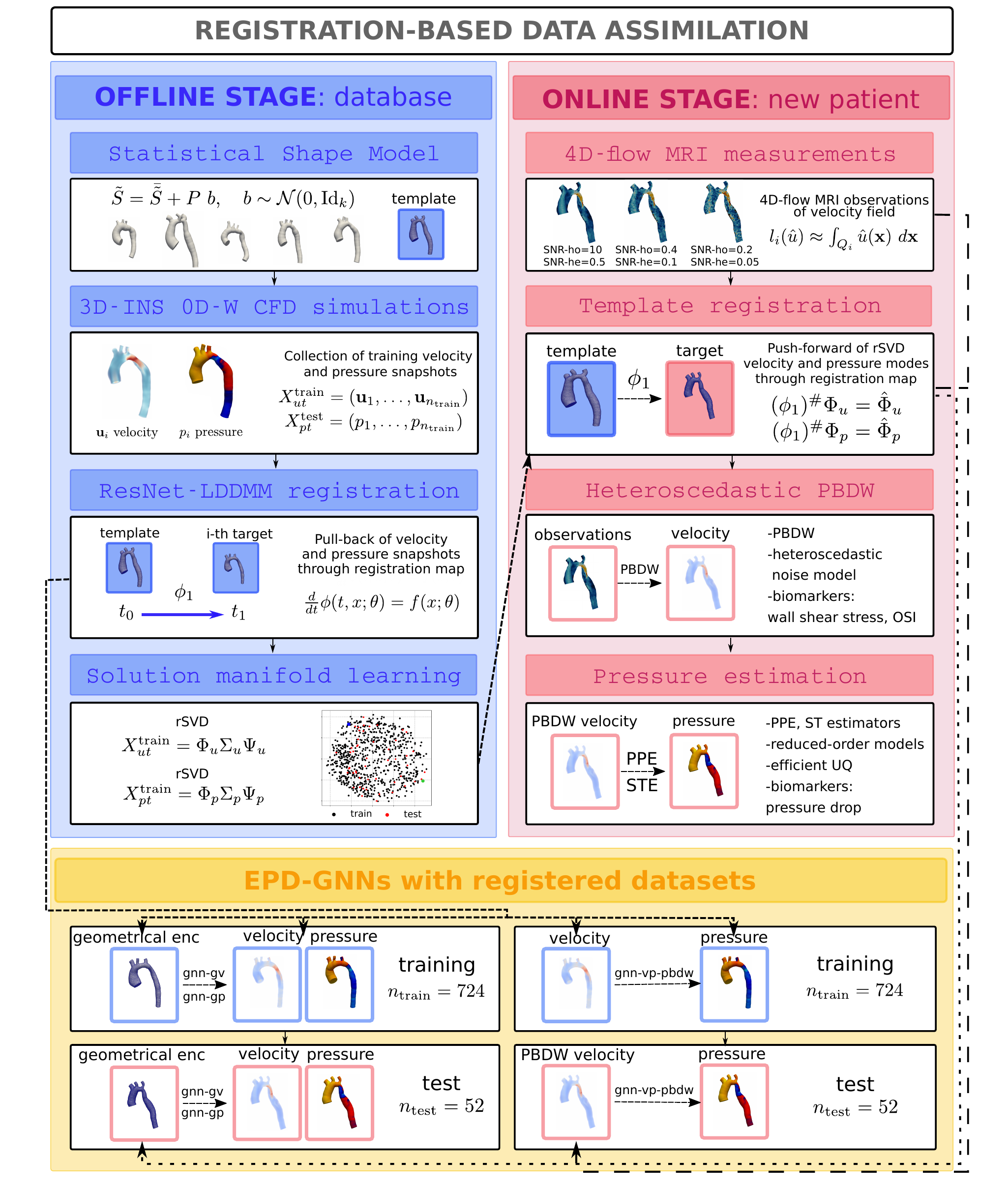}
  \caption{Overview of the registration-based data assimilation process.}
  \label{fig:scheme}
\end{figure}

The method is composed of an \textit{offline stage}, in which a database of training solutions is prepared and preprocessed on a template geometry, and of an \textit{online stage} in which 
velocity and pressure fields and related quantities of interest (e.g. time-averaged wall shear stress, oscillatory shear index, pressure drops) are inferred on a new patient. 
An overview of the data assimilation processes presented in this paper is shown in Figure~\ref{fig:scheme}. 
%

\paragraph*{Offline stage} In the offline stage, the first step is the design of a Statistical Shape Model (SSM) based on a centerline encoding of aorta geometries directly acquired from healthy subjects and aortic coarctation patients, which is used to generate additional synthetic geometries representing realistic anatomical features. An initial database containing 776 
geometries -- including both healthy and stenotic aortas -- is employed for generating velocity and pressure fields solving the three-dimensional
incompressible Navier--Stokes equations with 3-elements (lumped-parameter) Windkessel models on the main outlet branches, and with a variational multiscale turbulence model \cite{bazilevs2007variational}. The individual simulations are set up using patient-specific boundary conditions (b.c.),
based on measured inlet flow rates and tuning the parameters of the Windkessel models to control available or estimated flow rates.
Next, we use a registration algorithm to map the point clouds of the available patient geometries onto a reference shape, chosen within the cohort.
The registration algorithm is based on a large deformation diffeomorphic metric mapping (LDDMM)~\cite{trouve1998diffeomorphisms, trouve2005local}, which defines
the map between source and target points as the flow of a differential equation, whose vector field is parametrized by a residual neural network (ResNet)~\cite{amor2022resnet}.
We propose a cost functional tailored to the case of aortic meshes, taking into account additional information on the vessel centerline, the common topology of the shapes (inlet, outlets, and vessel wall),
and the normals to the boundary faces. 
As reported in \cite{pajaziti2023shape}, handling large meshes, as those required in the simulation of blood flows, leads to prohibitive computational costs for the training of the neural networks.
To overcome this problem we employ a multigrid optimization strategy, gradually
increasing the mesh size over the epochs. 
Crucial for the computational efficiency of our methodology is also the employment of GPUs with \texttt{pytorch3d}~\cite{ravi2020pytorch3d}. 
The multigrid ResNet-LDDMM registration is used to pullback the database of blood flow solutions onto the same reference shape. These data can be used for the solution manifold learning, to enable efficient linear dimension reduction methods such as randomized singular value decomposition (rSVD), and to define a geometrical encoding across shapes.
This implicit geometry encoding is used to accelerate the training of different Encode-Process-Decode GNNs~\cite{pfaff2020learning} that infer velocity and pressure solely from the geometry, or the pressure from the velocity data.

\paragraph*{Online stage} In the online stage, we assume to have available a new patient shape, possibly also with related velocity measurements. 
The geometry is first mapped on the reference shape via the ResNet-LDDMM registration, using the computed map to transport the global rSVD basis on the new domain: this enables the implementation of efficient data assimilation techniques to infer velocity, pressure, and related quantities of interest.
First, we focus on the Parametrized-Background Data-Weak (PBDW) method to reconstruct the high-resolution velocity field from partial observations minimizing the distance from the physics-informed space defined by the global rSVD basis. The PBDW, originally proposed in~\cite{MPPY2015}, has been recently used in different contexts to tackle data assimilation problems in hemodynamics \cite{galarce2022state} as well as to handle different noise models~\cite{gong2019pbdw}. We consider a generalized PBDW formulation extending the approach of~\cite{gong2019pbdw} to the case of heteroscedastic noise, to account for measurement data whose quality degrades close to the vessels' boundaries. 
Next, we investigate the estimation of pressure fields and pressure drops from velocity observations or solely from the geometry and b.c., and compare different approaches based on registration: EPD-GNNs against pressure-Poisson equation and Stokes pressure estimators (see, e.g.~\cite{bertoglio2018relative}). Using the global rSVD basis, the uncertainty quantification of these methodologies is also carried out.\newline

A classical LDDMM registration method to handle shape variability in model-order reduction for hemodynamics was first proposed in~\cite{guibert2014group} in the context of pulmonary blood flow, where a dataset of $17$ patients has been registered to a common template and then used to design a time-dependent ROM with a proper orthogonal decomposition.
The approach, however, was restricted to computational meshes with the same topology and validated only with limited variability of boundary conditions.
In~\cite{pajaziti2023shape} a dataset of $2800$ truncated healthy aortic arches has been generated with SSM and mapped with parametric non-rigid deformations, which, 
unlike LDDMM, do not guarantee the bijectivity of the registration maps. CFD simulations, limited to the stationary case, have been considered, and a shallow neural network was trained as a surrogate model in the space of proper orthogonal decomposition coordinates, in order to predict the velocity and pressure fields from the encoded description of the geometry in a shape vector. 
A proof-of-concept pipeline to handle shape variability was recently proposed in~\cite{galarce2022state}, employing PBDW and registration. The approach involves parametric deformations of cylinders and a registration based on LDDMM, while a nearest neighbor criteria is used to select at the online stage the proximal ROMs with a Hausdorff-like metric for the data assimilation.
A computational framework to address shape variability in hemodynamics was recently proposed in~\cite{Tenderini2024}, where registration was used to train
deep neural operators from a geometrical encoding obtained with conditioned neural ODEs, in a LDDMM fashion. The dataset was generated with the Radial Basis Function interpolation (RBF) of small deformations from a dataset of $20$ healthy patient-specific aortas. 
The procedure, tested on $358$ training and $39$ testing healthy geometries, was limited to the prediction of velocity field at systolic peak and considered limited variability of boundary conditions.
The use of machine learning methods for estimating the pressure field from reduced geometric representations of the hemodynamics has been discussed, among others, in~\cite{pegolotti2024learning,iacovelli2023novel, yevtushenko2021deep, versnjak2024deep},
focusing on the inference of pressure from one-dimensional vessel centerline information, and in~\cite{10.1115/1.4055285}, considering two-dimensional representations.
A deep learning frameworks for estimation of pressure for three-dimensional hemodynamics has been recently 
proposed in~\cite{nannini2025learninghemodynamicscalarfields}, without employing registration. An alternative approach, employing domain decomposable projection-based ROMs, has been proposed in~\cite{PEGOLOTTI2021113762}.

We focus on a realistic three-dimensional patient-specific representation of the computational domain as a general framework of interest in other research fields,
since some clinically relevant biomarkers are inherently two- or three-dimensional, for example the wall shear stress or the full description of the velocity field in proximity of a stenosis. Moreover, the evaluation of biomarkers from 3d fields is more interpretable than the use of surrogate models based on reduced 1d representations: non-physiological results can be assessed with more confidence looking at the 3d fields.

The main contributions of this work can be summarized as follows. 
Firstly, we propose a registration algorithm tailored to the case of aortic coarctation meshes including the main branches, that employs a multigrid optimization strategy, and that can handle realistic computational meshes.
Secondly, we investigate in detail the approximation properties and the potential of reduced rSVD bases created across different patient shapes for realistic meshes, time-dependent, three-dimensional hemodynamics, and patient-specific boundary conditions. 
Thirdly, we propose an EPD-GNN where the registration is used as a pre-processing step to speed up the training, in order to infer the velocity and pressure fields directly from the geometrical encoding, as well as the pressure field from velocity data.
Finally, we use the PBDW method built on the global rSVD basis and extended to the case of heteroscedastic noise, to account for 
the higher uncertainty of 4DMRI images close to the vessel boundaries.
We then validate and compare the different algorithms for data assimilation to reconstruct velocity and pressure solutions, as well as related biomarkers, from coarse velocity observations 
mimicking 4DMRI data. In particular, we investigate a pressure estimator that combines PBDW velocity reconstruction with EPD-GNN inference to reconstruct pressure quantities of interest
from coarse-grained and noisy velocity data, validating it against state--of--the art pressure estimators.
Finally, we assess the benefit of employing registration as a pre-processing step to speed up the training of GNNs and compare their prediction accuracy with PBDW for the inference of the velocity field and with PPE and STE for the inference of the pressure field.

The results show that (i) the physical variability of the problem is more complex than what has been shown in the literature, (ii) the intrinsic dimensionality of the solution manifold is high, more than $500-1000$ rSVD modes are used to estimate, with a relative accuracy under $10\%$, the pressure and velocity fields, (iii) while the employment of linear ROMs is possibly not feasible, EPD-GNNs have good potential but require more data and a higher computational budget with respect to that available for our studies, (iv) direct inference from the geometry is a harder problem than exploiting velocity observations, (v) standard pressure estimators are outperformed even in our limited data regime.

The rest of the paper is structured as follows. Section~\ref{sec:setting} introduces the key elements of the computational framework for aortic shape modelling and blood flow simulations. The registration algorithm is described in section~\ref{sec:registration}, while the detailed validation and employment of registration for solution manifold learning is presented in section~\ref{sec:sml}. In section~\ref{ssec:pres-gnn}, we present the hyperparameter studies for the EPD-GNNs architecture and how we employ the registered datasets in this context. Sections~\ref{sec:da} and~\ref{sec:prec} focus on the estimation of velocity and pressure fields from medical imaging data. Finally, in section~\ref{sec:discussions} we discuss the results and limitations of the approach and in section~\ref{sec:conclusions} we draw the conclusions and present possible future directions of research.

\section{Forward computational hemodynamics}
\label{sec:setting}
\subsection{Statistical shape modelling of patients with aortic coarctation}
\label{subsec:ssm}

The data used in this study were obtained from a cohort of patients with coarctation of the aorta (CoA), augmented synthetically using statistical shape models (SSM). The procedure is briefly outlined below. For the detailed methodology, we refer the reader, e.g. to 
~\cite{goubergrits2022ct, thamsen2021synthetic,thamsen2020unsupervised,versnjak2024deep, yevtushenko2021deep}.

The initial database contained $228$ surfaces acquired from 3D steady-state free-precession (SSFP) magnetic resonance imaging (MRI)
(acquired resolution $\SI{2}{mm}\times \SI{2}{mm}\times \SI{4}{mm}$, reconstructed resolution used for surface reconstruction $\SI{1}{mm}\times \SI{1}{mm}\times \SI{2}{mm}$) and segmented with \texttt{ZIB Amira}~\cite{stalling2005amira}. In total, 106 CoA patients (32 female) and 85 healthy subjects were acquired (25 female). 
For 37 (8 female) of the 106 CoA patients also post-treatment image data were available, thus increasing the database. The median age was 21 years with interquartile range (IQR) of 32 years.
The considered region of interest comprises the vessel surface of aortic arch up to the thoracic aorta (TA), including three main branches 
and the corresponding boundary surfaces (brachiocephalic artery, BCA, left common carotid artery, LCCA, left subclavian artery, LSA). Few available cases with two or four branches of the aortic arch were not included into the database.

\begin{figure}[!htp]
  \centering
  \includegraphics[width=0.9\textwidth]{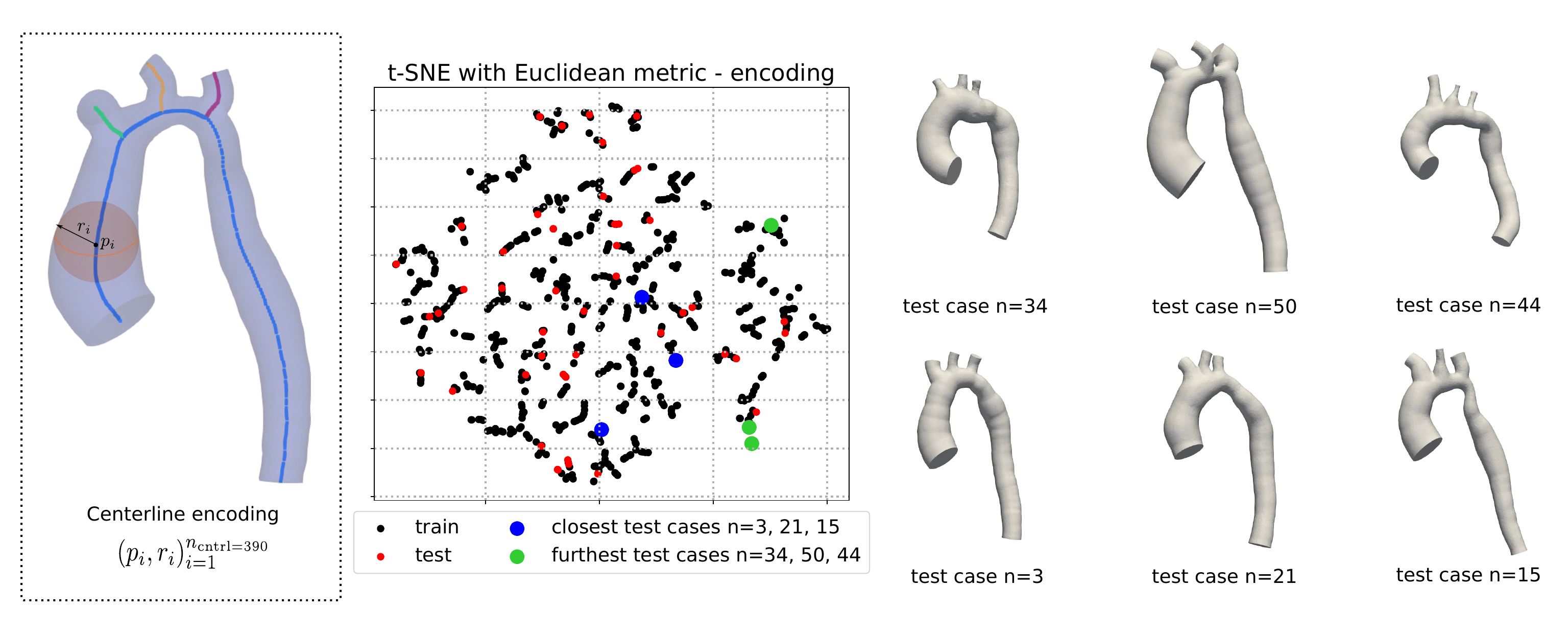}
  \caption{
  \textbf{Left:} Sketch of the centerline encoding (points $p_i$ and radius $r_i$ of the associated inscribed sphere, for $i\in\{1,\dots,390\}$). 
  \textbf{Center:} Clustering of the considered training ($n=724$) and test ($n=52$) shapes using t-SNE
  with the Euclidean distance on a geometrical encoding, based on the distance of each point from the centerline after shape registration (see section~\ref{subsec:sml_correlations} for details).
  \textbf{Right:} Visualization of the furthest shapes (top, test cases $34$, $50$, and $44$) and the closest ones (bottom, test cases $3$, $21$, $15$) according to the metric in the center plot.}
  \label{fig:clustergeometries}
\end{figure}

Additionally, pointwise linear centerlines for the aorta and the three branching vessels have been obtained along with the radii of the inscribed spheres using the vascular modelling toolkit \texttt{VMTK}~\cite{antiga2008image} (see the sketch in figure~\ref{fig:clustergeometries} (left) for an example).

The procedure resulted in $300$ centerline points for the aorta and $30$ centerline points for each branching vessel, for a total of $n_{\text{cntrl}}=390$ points and corresponding radii of inscribed spheres for each considered shape.
These data allow to encode the morphology of each shape into a matrix $S_{\text{SSM}}\in\mathbb{R}^{n_{\text{cntrl}}\times (3+1)}$ containing the spatial coordinates of the $n_{\text{cntrl}}$ centerline points and the associated radii. Closed triangulated surfaces are then generated from this skeletal representation~\cite{yevtushenko2021deep}. 
Each geometry is rigidly moved towards the mean shape $S_{\rm mean}$, 
minimizing the least-squares distance between points with the closest point algorithm using \textsf{mcAlignPoints} package of the \texttt{ZIB Amira} software.
No scaling is performed. New shapes 
are generated through SSM with Principal Component Analysis (PCA):
\begin{equation*}
  \Tilde{S}_{\text{SSM}} = S_{\rm mean}  + P_{\text{SSM}}\ b_{\text{SSM}}, 
\end{equation*}
where $P_{\text{SSM}}\in\mathbb{R}^{(n_{\text{cntrl}}\times (3+1))\times k}$ contains $k>0$ truncated modes of the correlation matrix of the training shapes used for SSM, and 
$b_{\text{SSM}}\in\mathbb{R}^k$ is the vector of coefficients. For the SSM development only pre-treatment CoA shapes ($93$ cases) and healthy aorta ($65$ cases) were used.

A database of more than $10000$ shapes is generated sampling $b_{\text{SSM}}$ from a normal distribution. Unrealistic shapes, e.g., containing self-intersection, small vessel radius (below $1.0$ mm), or 
excessive degree of stenosis (less than $20\%$ or greater than $80\%$) have been removed. 
As aortic length and aortic inlet diameter are correlated with age, when age ranges deduced from these two morphological parameters did not overlap, the corresponding shapes were discarded. Further shapes have been removed performing a preliminary CFD analysis of peak systole flow using \texttt{STAR-CCM+}~\cite{yevtushenko2021deep}, discarding 
those resulting in unphysical quantities of interest.

This procedure resulted in a cleaned database of $1312$ (real and synthetic) shapes, represented by triangulated surface meshes, centerline points, and centerline radii.
Along the curse of this study, $437$ additional geometries have been removed based on the results of time dependent simulations (see section \ref{ssec:blood_flow}) and further $99$ due to 
inaccurate registrations (see section \ref{sec:registration}). The remaining $776$ cases have been split in $724$ training and $52$ testing shapes.
Figure \ref{fig:clustergeometries}, center and right plots, show qualitatively the extent of the training and test datasets, based on a 
T-distributed Stochastic Neighbor Embedding~\cite{van2008visualizing} (t-SNE) with the Euclidean distance on a geometrical encoding of the shapes
that relies on the shape registration map (further details will be given in section~\ref{subsec:sml_correlations}).

\subsection{Blood flow modelling}\label{ssec:blood_flow}
\begin{figure}[!htp]
  \centering
  \includegraphics[width=0.4\textwidth, trim={0 0 0 20}, clip]{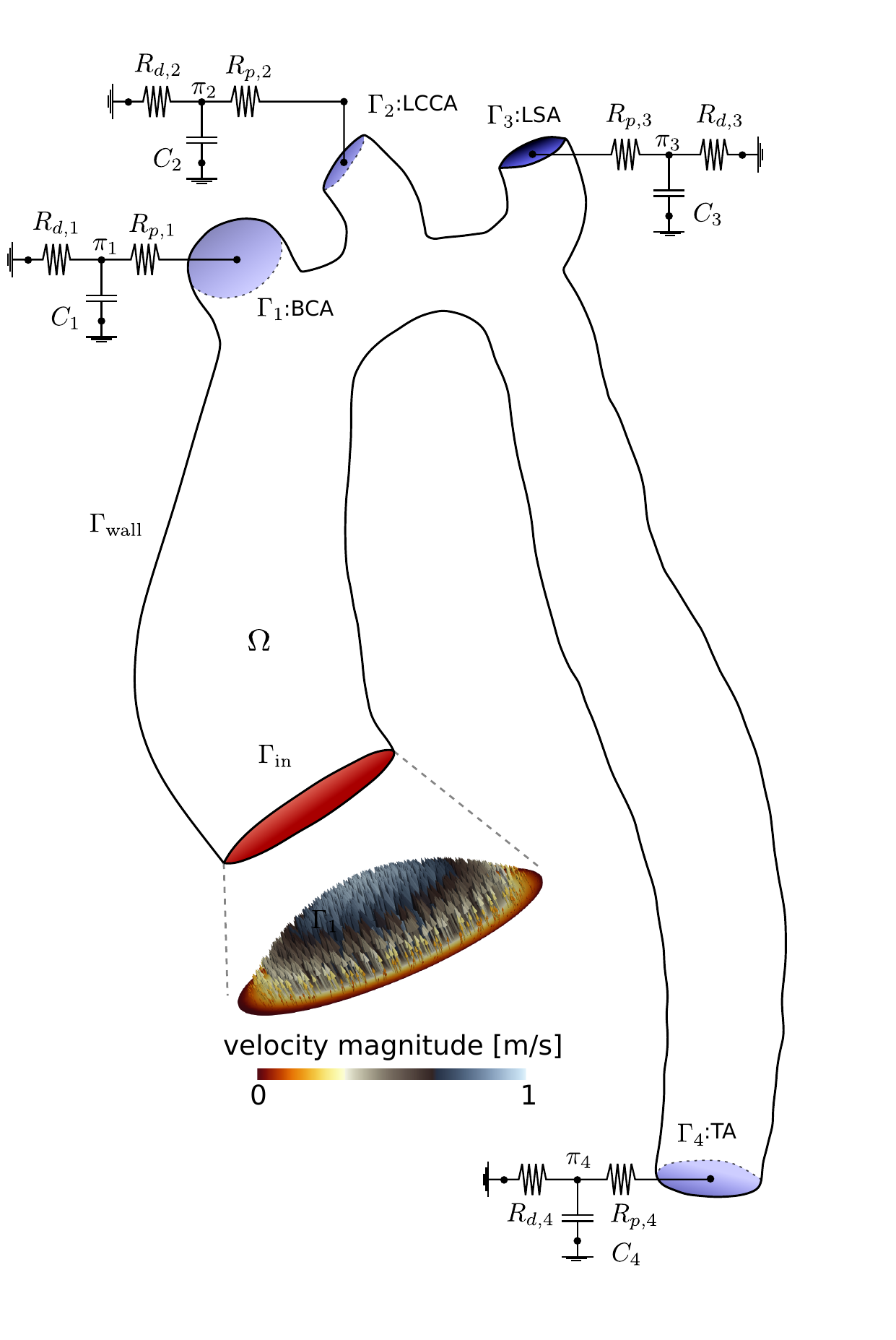}
  \caption{Example of a computational domain. A parabolic profile is imposed on the inlet boundary (ascending aorta), based on a given peak flow rate. 
  Windkessel models are used at the outlets (thoracic aorta, brachiocephalic artery, left common carotid artery, and left subclavian artery).}
  \label{fig:domain
}
  \label{fig:domain}
\end{figure}

Let us denote with $\Omega \subset \mathbb R^3$ the computational domain representing a generic shape from the considered dataset, whose boundaries can be decomposed as
\begin{equation}\label{eq:omega_bnd}
\partial \Omega = \Gamma_{\text{wall}} \cup \Gamma_{\rm in} \cup \left( \bigcup_{i=1}^4 \Gamma_i \right),
\end{equation}
distinguishing between the vessel \text{wall} $\Gamma_{wall}$, the inlet boundary $\Gamma_{\rm in}$, and the four outlet boundaries (BCA, LCCA, LSA, TA), as depicted in figure \ref{fig:domain}. We assume that the blood flow in the considered vessels behaves as an incompressible Newtonian fluid and thus describes the hemodynamics via the incompressible Navier--Stokes equations for the velocity $\bu:\Omega \to \mathbb R^3$ and the pressure fields $p:\Omega \to \mathbb R$:
\begin{equation}\label{eq:3dnse}
\left\{
\begin{aligned}
\rho \partial_t \mathbf{u}+\rho \mathbf{u}\cdot\nabla\mathbf{u}+\mu\Delta\mathbf{u}-\nabla p=\mathbf{0},\qquad&\text{in}\ \Omega,\\
    \nabla\cdot\mathbf{u}=0,\qquad&\text{in}\ \Omega,
\end{aligned}
\right.
\end{equation}
where $\rho=\SI{1.06e3}{\kilogram\per\meter^3}$ stands for the blood density, and $\mu=\SI{3.5e-3}{\second\cdot\pascal}$ is the dynamic viscosity.

Equations~\eqref{eq:3dnse} are complemented by homogeneous Dirichlet boundary conditions for the velocity on $\Gamma_{\rm \text{wall}}$, i.e., neglecting  fluid-structure interactions between the blood flow and the vessel wall, by a Dirichlet boundary condition on $\Gamma_{\rm in}$, imposing a parabolic flow profile at the inlet~\cite{katz2023impact}, and by lumped parameter models on the four outlet boundaries, i.e., 
\begin{equation}
  \label{eq:3dnse-bc}
\left\{
\begin{aligned}
\mathbf{u}&=\mathbf{u}_{\rm in},\qquad &&\text{on}\ \Gamma_{\rm in},\\
\mathbf{u}&=\mathbf{0},\qquad &&\text{on}\ \Gamma_{\text{wall}},\\
    -pI+\mu\left(\nabla\mathbf{u}+\nabla\mathbf{u}^T\right)&=-P_i\mathbf{n},\qquad &&\text{on}\ \Gamma_i,\ i\in\{1,2,3,4\}.
\end{aligned}
\right.
\end{equation}
In the last equation, $\mathbf{n}$ denotes the outward normal vector to the fluid boundary and $P_i$ stands for an approximation of the outlet pressure imposed on the boundary $\Gamma_i$, which is evaluated as a function of the boundary flow rates $Q_i:=\int_{\Gamma_i} \mathbf{u}\cdot\mathbf{n}$, $i=\in\{1,2,3,4\}$, 
via a 3-elements (RCR) Windkessel model~\cite{westerhof2009arterial}
\begin{equation}\label{eq:wk-rcr-i}
\left\{
\begin{aligned}
C_{d,i}\frac{d\pi_i}{dt}+\frac{\pi_i}{R_{d,i}}=Q_i,\qquad\qquad\qquad\qquad\ \;\,\quad&\text{on}\ \Gamma_i,\ i\in\{1,2,3,4\},\\
   P_i=R_{p,i}Q_i+\pi_i,\qquad\qquad\qquad\qquad\ \;\,\quad&\text{on}\ \Gamma_i,\ i\in\{1,2,3,4\},\\
\end{aligned}
\right.
\end{equation}
depending on an auxiliary \textit{distal} pressure $\pi_i$, a \textit{proximal} resistance $R_p$ (modeling the resistance to the flow of the arteries close to the open boundary), a \textit{distal} resistance $R_d$ (modeling the downstream resistance of the rest of the cardiovascular system), and a \textit{capacitance} $C_d$ (modeling the compliance of the cardiovascular system).
The tuning of the Windkessel parameters will be discussed in detail in section \ref{ssec:bc-calibration}.
Peak inlet flow rates for each shape were provided as part of the patient cohort data for each considered shape. These values  were adjusted to match a parabolic shape on the inlet boundary, and multiplied by a time dependent function to obtain the inlet boundary condition over time.

\subsection{Calibration of boundary conditions across the cohort of patients}\label{ssec:bc-calibration}
The Windkessel parameters might have a considerable impact on the solution and the calibration typically depends on the flow regime of interest, on the
particular anatomical details, and on available data. For the purpose of this study, we opted for an approach driven by the flow split across
the different branches.

We introduce the total resistances
$R_i := R_{p,i} + R_{d,i}$ and the equivalent
\textit{systemic} resistance $R_S^{-1}  := \sum_{i=1}^N R_i^{-1}$. Neglecting the contribution of the 3D domain to the total resistance, the quantities
$R_i/ R_S$, for $i=1,2,3,4$, can be used to control the \textit{flow split}, i.e., the ratio of the inlet flow $Q_{\rm in}$ that flows, on average, through each outlet.

We have tuned the template geometry's parameters considering a flow split of $50\%$ for the BCA, and of $25\%$ for LCCA and LSA, as in~\cite{katz2023impact}. For the systematic calibration of the Windkessel parameters on other geometries, we consider a model for the average flow split based on the following steps:
\begin{enumerate}[itemsep=2pt, left=0pt, labelsep=5pt]
  \item A rescaling of the systemic resistance, based on the patient specific inlet,
  \begin{equation*}
  R_S :=  \frac{\hat Q_{\rm in}}{Q_{\rm in}} \hat R_S,
  \end{equation*}
  \item A shape-specific flow split based on the reference mean velocities and the outlet areas of the new patient,
  \begin{equation*}
      \sigma_i := \hat Q_i \frac{A_i}{\hat A_i} \frac{1}{\sum_{j=1}^4 \hat Q_j \frac{A_j}{\hat A_j}} = \frac{\hat u^{\mathrm{mean}}_i A_i }{ \sum_{j=1}^4 \hat u^{\mathrm{mean}}_j A_j},
  \end{equation*}
      \item The approximation of the patient-specific total resistance 
      \begin{equation*}
        R_i = \sigma_i R_S,
      \end{equation*}
      \item A splitting between proximal and distal resistance (the same for all  patients),
  \begin{equation*}
  R_{p,i}  = 0.1 R_i,\ R_{d,i}=0.9 R_i\,.
  \end{equation*}    
\end{enumerate}
 Finally, capacitances are defined proportionally to the area of the outlet boundaries, i.e., 
\begin{equation*}
  C_i = \frac{A_i}{A_{\rm tot}} C_{\rm tot},
\end{equation*}
as a fraction of the total capacitance $C_{\rm tot} = 10^{-8}$ (the same for all shapes).

\begin{figure}[!htp]
  \centering
  \includegraphics[width=0.49\textwidth]{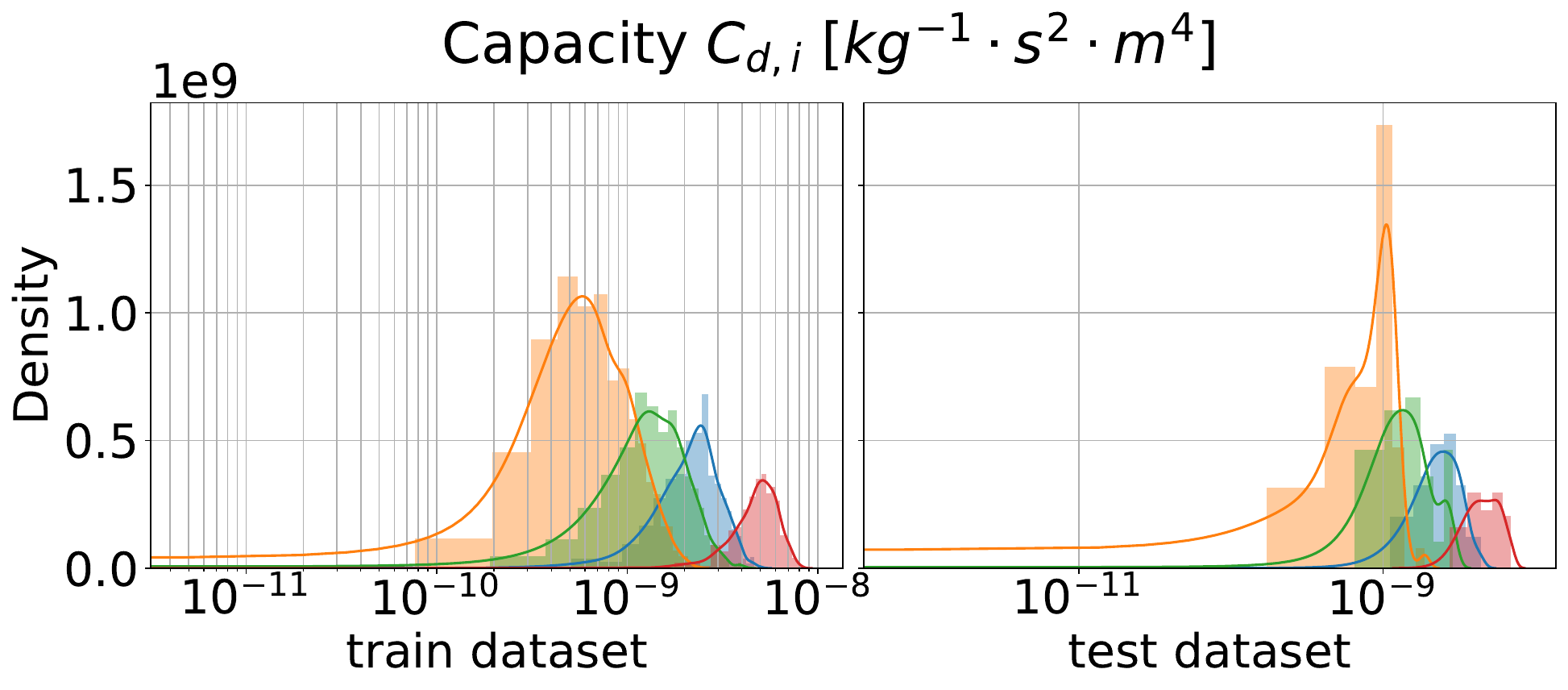}
  \includegraphics[width=0.49\textwidth]{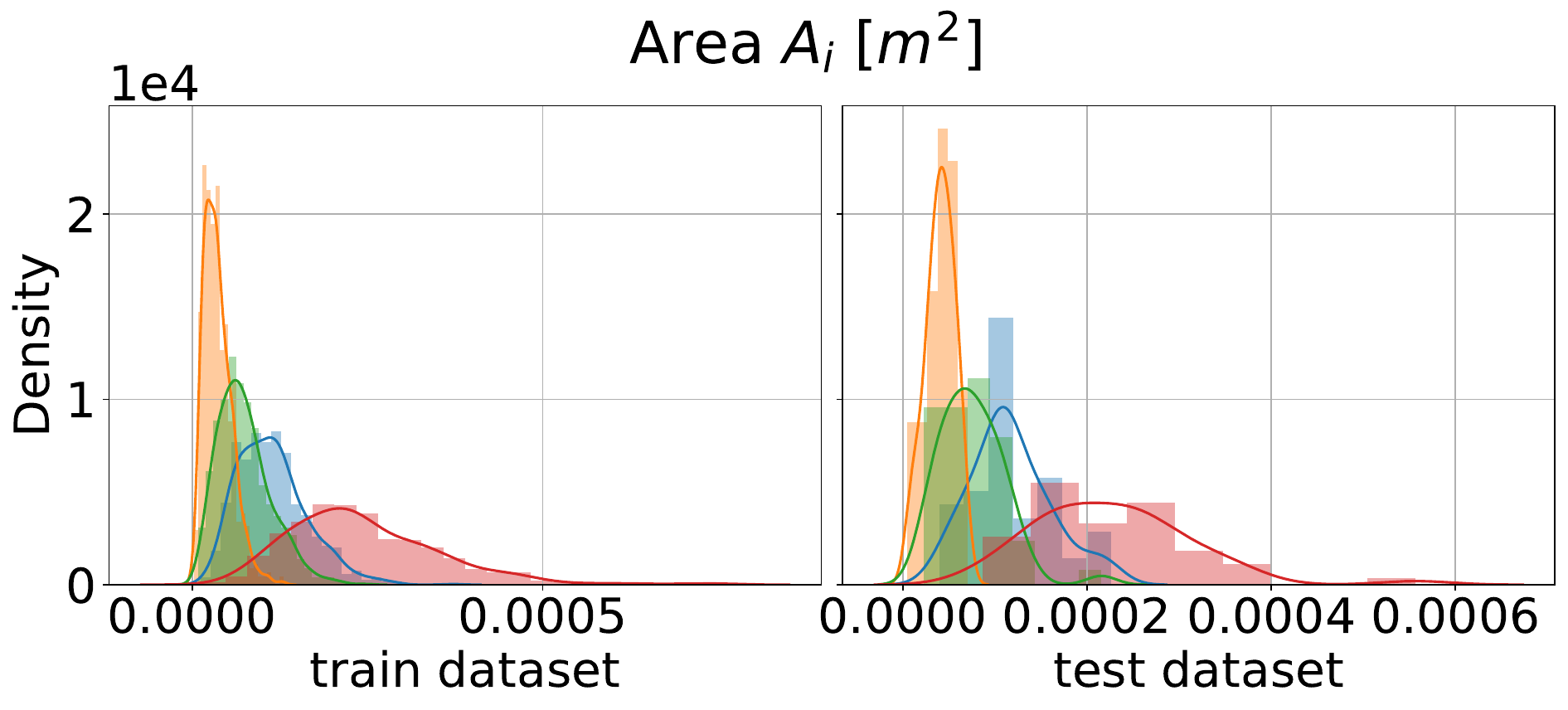}\\
  \includegraphics[width=0.49\textwidth]{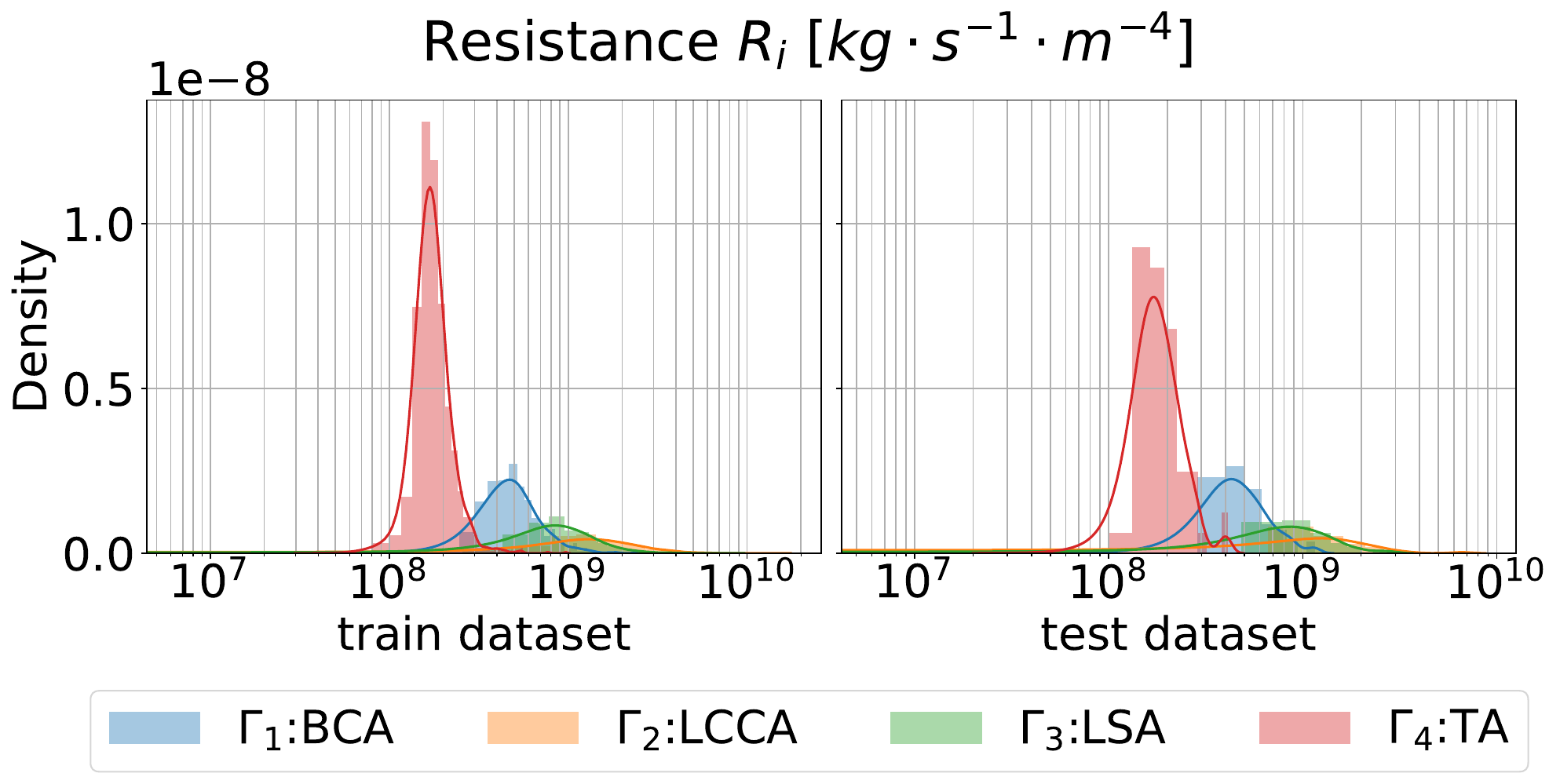}
  \caption{Distribution of Windkessel parameters $C_{d,i}$ (distal capacity), $R_i=R_{d,i}+R_{p,i}$ (Total resistance), and $A_i$ (boundary area) for the different outlets, across the training and test datasets.}
  \label{fig:param_distr}
\end{figure}

The distributions of the total resistances, distal capacities, and outlets areas, for the complete shape database, are shown in figure~\ref{fig:param_distr}.

\begin{rmk}[Calibration based on flow split]
The main motivation behind this approach is the fact that the flow split can be experimentally measured non-invasively on the different sections, or inferred according 
to existing literature data and patient anatomy. Moreover, using the flow split as parameter allows modelling different physiological (rest/exercise)
or pathological (e.g., obstruction of vessels downstream) conditions, and can be used to enrich the solution dataset depending on the context of interest.
\end{rmk}

\subsection{Synthetic dataset of aortic shapes and numerical simulations}
\label{ref:numsim}
We solve numerically system \eqref{eq:3dnse}-\eqref{eq:3dnse-bc} for each of the $776$ considered geometries, discretizing the corresponding volume with a tetrahedral mesh reated from the original surface shape
and imposing shape specific boundary conditions as described in section \eqref{ssec:bc-calibration}. 
We use stabilized equal-order linear finite elements for velocity and pressure and a BDF2 time marching scheme, 
with a semi-implicit treatment of the non-linear convective term and of the VMS turbulence model~\cite{bazilevs2007variational, forti2015semi}.
Further details on the discretization and on the numerical method are provided in  appendix \ref{appendix:weak}.
The ODEs \eqref{eq:wk-rcr-i} are solved using an implicit Euler scheme, and the coupling at the boundary is implemented explicitly, i.e., using the boundary pressures at the previous time iteration to impose Neumann boundary conditions on each outlet. Some velocity snapshots are shown in figure~\ref{fig:snapsu}.
The solver is implemented in the computational framework \texttt{lifex-cfd}~\cite{AFRICA2024109039}, based on the open-source library \texttt{deal.II}~\cite{arndt2022deal}.
Simulations have been run for five hearth beats, only the last period is considered, in order to ensure a quasi-periodic state.

\begin{figure}[!htp]
  \centering
  \includegraphics[width=0.9\textwidth]{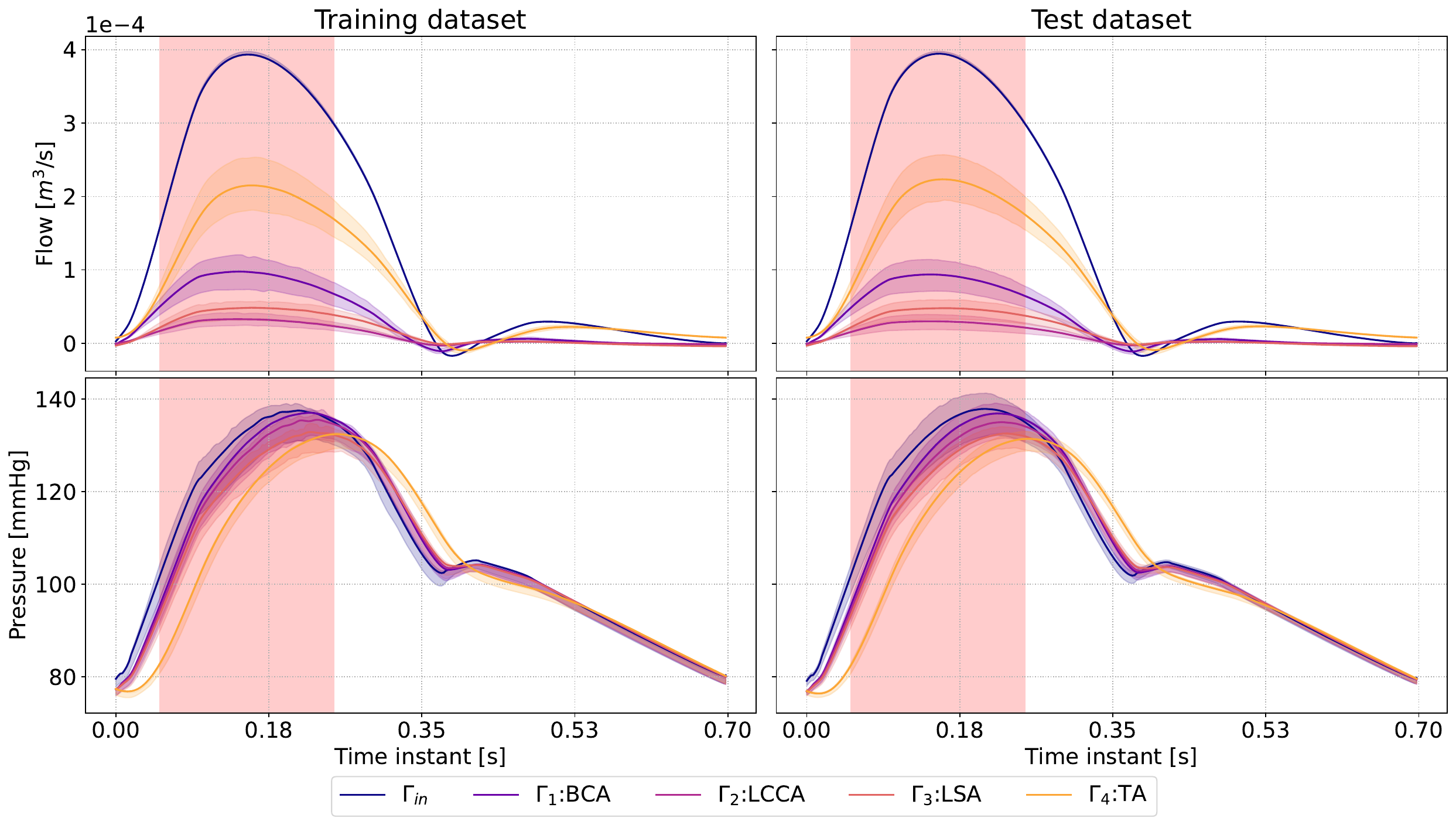}
  \caption{\textbf{Top: } Numerical results for te flow at inlet (AAo, with opposite sign) and outlets (BCA, LCCA, LSA, TA) boundaries. 
  \textbf{Bottom: }  Numerical results for the pressure at inlet (AAo) and outlets (BCA, LCCA, LSA, TA) boundaries. The $25$-th and $75$-th percentile across all the $724$ and $52$ training and test data are shown. The red vertical bands correspond to the time window $t\in[0.05s, 0.25s]$ which is the focus of our data assimilation studies, see remark~\ref{rmk:timewindow}.}
\label{fig:flow_and_pressure}
\end{figure}

The results (figure~\ref{fig:flow_and_pressure}) show a rather uniform distribution of flow and pressure values at boundaries across the dataset.

\section{Registration with ResNet-LDDMM}
\label{sec:registration}
\subsection{Large deformation diffeomorphic metric mapping}\label{ssec:resnet-lddmm-intro}
The registration, or image matching problem, consists in smoothly mapping a \textit{source} (or template) image into a \textit{target} image. 
Our approach is based on the so-called Large Deformation Diffeomorphic Metric Mapping (LDDMM), in which
the map between the source and the target is sought as a diffeomorphic flow of an ODE~\cite{bruveris2017completeness,dupuis1998variational}.
In this section, we present the formulation of the method from a continuous perspective and the main theoretical background. 
The specific implementation to the case of three-dimensional
meshes of aortic shapes will be described in detail in section \ref{subsec:resnetlddmm}.

Formally, let us consider the source and target images defined by the characteristic functions $\chi_S:\mathbb{R}^3\rightarrow\mathbb{R}$ and $\chi_T:\mathbb{R}^3\rightarrow\mathbb{R}$, respectively.
We assume that both images are contained in an open bounded set, i.e., 
\begin{equation*}
\text{supp}(\chi_S)\cup\text{supp}(\chi_T)\subset G\subset\mathbb{R}^3.
\end{equation*}

The goal of LDDMM is to find a diffeomorphic map between the source and the target as a one-parameter group of diffeomorphisms  $\{\phi(t, \cdot)\}_{t\in I}$, $I :=[0,1]$, defined
as the flow of an ODE depending on a vector field $f:I\times\mathbb{R}^3\rightarrow\mathbb{R}^3$, i.e., such that
\begin{equation}\label{eq:phi_t_lddmm}
  d_t \phi(t, \x) = f(t, \phi(t, \x)),\qquad \phi(0, \x) = \x,\qquad\forall\x\in G.
\end{equation}
In particular, $\chi_S\left(\phi(0,\cdot)\right) = \chi_S$ coincides with the source, and  $\chi_S\left(\phi(1,\cdot)\right)$ is the mapped image.
The problem is addressed in an optimal control framework, where the control is the vector field $f$, minimizing the matching error between the mapped source and the target.
The following result ensures the existence of an optimal solution in the considered setting for arbitrary dimension $d$.

\begin{theorem}[LDDMM registration, theorem 3.1~\cite{dupuis1998variational}]
  \label{def:regpb}
  Assume that $S:G\rightarrow\mathbb{R}$ and $T:G\rightarrow\mathbb{R}$ are two bounded measurable functions, and that $S$ 
 \textit{(e.g., the source image)} is continuous almost everywhere. 
 Let us suppose that $f\in L^2(I, H^s)$, with $s>d/2 + 1$ (the Sobolev embedding implies $f\in \mathcal{C}^{1, \alpha}(I, \mathcal{C}^{1, \alpha})$ with $\alpha=s - \tfrac{d}{2} - 1>0$), then the following registration problem
  \begin{equation}\label{eq:reg_miminization}
    \min_{f\in L^2(I, H^s)} \int_{0}^{1}\lVert f(t, \cdot)\rVert^2_{H^s}\ dt + \int_G |S(\phi(1, \x))-T(\x)|^2\ d\x
  \end{equation}
  has a minimizer.
\end{theorem}

In particular, theorem \ref{def:regpb} holds for signed distance functions, as well as for the particular case of characteristic functions. 
The regularity condition $f\in L^2(I, H^s)$, with $s>d/2 + 1$ from the Sobolev embedding theorem, is sufficient to obtain that $\{\phi(t, \cdot)\}_{t\in I}\subset\text{Diff}^1(G)$ is a one-parameter group of diffeomorphisms. 
In general,  for a generic Hilbert space $\mathcal{H}$ such that $\mathcal{H}\hookrightarrow C^1_b(\mathbb{R}^d)$, the group $\mathcal{G}_\mathcal{H}(\mathbb{R}^d)$ of diffeomorphisms generated by vector fields in $L^2(I, \mathcal{H})$ is only strictly contained in $\text{Diff}^1(\mathbb{R}^d)$.
However, if $s>d/2+1$, it holds (see ~\cite{bruveris2017completeness}) 
\[
\mathcal{G}_{H^s}=\mathcal{D}^s(\mathbb{R}^d) =\{\phi\in\text{Diff}^1(\mathbb{R}^d) \mid \phi\in\text{Id}+H^s(\mathbb{R}^d,\mathbb{R}^d)\}.
\]

Using the universal approximation theorems of neural networks, one can infer the existence of neural networks that approximate the vector field $f$ arbitrarily well, 
with convergence estimates depending on its regularity.
\begin{theorem}[Existence of ResNet-LDDMM vector field]
  \label{theo:existreg}
Let  $f:I\times\mathbb{R}^d\rightarrow\mathbb{R}^d$ be a minimizer for problem \eqref{eq:reg_miminization}, let
$\{\phi(t, \cdot)\}_{t\in I}$ the corresponding group of diffeomorphisms \eqref{eq:phi_t_lddmm}, and let $\epsilon>0$. 
\begin{itemize}
\item[(i)]There exists a neural network (NN) $f_\epsilon: I\times \mathbb{R}^d\rightarrow  I\times \mathbb{R}^d$, $f_\epsilon\in\mathcal{C}^{0, 1}(I, \mathcal{C}^{0, 1})$,  with ReLU activations and one hidden layer, such that:
  \begin{equation*}
    \lVert \phi_\epsilon - \phi\rVert_{L^2(I, L^2)}\rVert \leq C_1 \lVert f_\epsilon - f\rVert_{L^2(I, L^2)} \leq \epsilon,
  \end{equation*}
  where  $\{\phi_\epsilon(t, \cdot)\}_{t\in I}$ satisfy
 \begin{equation*}
    d_t \phi_\epsilon(t, \x) = f_\epsilon(t, \phi_\epsilon(t, \x)),\qquad \phi_\epsilon(0, \x) = \x,\qquad\forall\x\in G,
  \end{equation*}
and $C_1>0$ is independent of $\epsilon$.
  
\item[(ii)] Let $f_N\in\mathcal{C}^{2, 1}(I, \mathcal{C}^{2, 1})$ denote a deep NN, with a  fixed number of layers, ReCU activations, $N$ non-zero weights, and
let $\{\phi_N(t, \cdot)\}_{t\in I}$ be the flow of the corresponding ODE, i.e., 
  \begin{equation*}
    d_t \phi_N(t, \x) = f_N(t, \phi_N(t, \x)),\qquad \phi_N(0, \x) = \x,\qquad\forall\x\in G.
  \end{equation*}
For any $n\in\{0, 1\}$, and $\forall m\in\mathbb{N}$, $m\geq n+1$, the following estimate holds
  \begin{equation*}
    \lVert \phi_N - \phi\rVert_{L^2(I, H^n)}\rVert \leq C_2 \lVert f_N - f\rVert_{L^2(I, H^n)} \leq C_3\cdot N^{-\frac{m-n}{2(d+1)}},
  \end{equation*}
  with positive constants $C_2, C_3>0$ independent of $N$.
  \end{itemize}
\end{theorem}

\begin{proof}
See appendix~\ref{appendix:convergence}.
\end{proof}

In the practical implementation (see section~\ref{subsec:resnetlddmm}) we will 
consider only autonomous vector fields $f:\mathbb{R}^3\rightarrow \mathbb{R}^3$, i.e., that do not depend on $t\in I$,  and
NNs with ReLU activations, $7$ hidden layers, and where the first layer's input is augmented with Random Fourier Features (RFF).
Moreover, the matching error between mapped source and target images will be measured using a metric based on the Chamfer distance, which is a natural choice
when considering discrete point clouds, rather than the $L^2$-norm between characteristic functions used in \eqref{eq:reg_miminization}.
Up to our best knowledge, general existence results using the Chamfer distance are not available. However, it can be
shown that the solutions to a discrete version of~\eqref{eq:reg_miminization} with the discrete $L^2$-norm as discrepancy metrics convergence to the solution of the continuous registration problem (see appendix~\ref{appendix:convergence}). This result motivates also the usage of a multigrid optimization.

\subsection{Multigrid ResNet-LDDMM for aortic shape meshes}
\label{subsec:resnetlddmm}

\paragraph{Preliminaries} 
In this section, we address the shape registration between 3d computational domains representing different aortic surfaces, discretized by triangular meshes.  
We assume that, for all considered shapes, it is possible to subdivide each surface mesh into an inlet boundary $\Gamma_{\rm in}$, four outlet boundaries $\Gamma_1, \hdots,\Gamma_4$, 
and the wall boundary $\Gamma_{wall}$ (see equation \eqref{eq:omega_bnd}).
All shapes have are also characterized by a piecewise centerline $l\in\mathbb{R}^{n_{\text{cntrl}}\times 3}$ with a main branch (ascending and thoracic aorta) and three minor branches.
By construction (see section \ref{subsec:ssm}), we also assume that all centerlines have the same number of points $n_{\text{cntrl}}=390$.
These assumptions are motivated by the fact that all computational domains shall provide suitable discretizations of the physical model of interest 
(Equations \eqref{eq:3dnse}, with boundary conditions \eqref{eq:3dnse-bc}).

However, \textit{no assumptions} are made on the number of vertices, edges, or faces in each mesh nor on the connectivity of the elements.
In what follows, a shape $\mathcal S$ will be generally defined by the corresponding surface mesh, i.e., a pair $(X_{\mathcal S}, F_{\mathcal S})$, where
\begin{itemize}
  \item $\XS=\{\mathbf x_i\}_{i=1}^{n_{p,\mathcal S}}\subset\mathbb{R}^3$ is a point cloud with cardinality $n_{p,\mathcal S}>0$.
  \item $\FS = \{(a_i, b_i, c_i), a_i \neq b_i, a_i \neq c_i, b_i \neq c_i \}_{i=1}^{n_{f,\mathcal S}}\subset\mathbb{N}^3$, is a set of triangular faces with cardinality $n_{f,\mathcal S}>0$, where the element $(a_i, b_i, c_i)$ corresponds to the face 
defined by three (different) vertices $(\mathbf x_{a_i}, \mathbf x_{b_i}, \mathbf x_{c_i})$,
\end{itemize}
and its centerline $l_{\mathcal S}\in\mathbb{R}^{n_{\text{cntrl}}\times 3}$ ($n_{\text{cntrl}} = 90$, equal amount for all shapes).
We also introduce the set 
\begin{equation*}
\NS := \{\mathbf n_i^{\mathcal S}, i=1,\hdots,n_{p,\mathcal S}\} 
\end{equation*}
of normal vectors to the mesh vertices defined, for each $\mathbf x_i\in \XS$, as the linear combination of 
the normals of all faces adjacent to the vertex $\mathbf x_i$, weighted by the arccosine of the angles corresponding to $\mathbf x_i$ in the respective adjacent triangular faces, and renormalized such that $\lVert \mathbf n_i^{\mathcal S}\rVert=1$.

Let us also introduce, for two arbitrary point clouds $X_{\mathcal{S}}$ and $X_{\mathcal{T}}$, the \textit{Chamfer} distance
\begin{equation}
  \label{eq:classical_Chamfer}
  \mathcal{D}_{\text{Chamfer}}(X_{\mathcal{S}}, X_{\mathcal{T}}) = \frac{1}{n_{p,\mathcal{S}}}\left(\sum_{\mathbf x\in X_{\mathcal{S}}} \min_{ \mathbf y\in X_{\mathcal{T}}} \lVert \mathbf{\mathbf x} - \mathbf{\mathbf y} \rVert_2\right) + \frac{1}{n_{p, \mathcal{T}}}\left(\sum_{\mathbf y \in X_{\mathcal{T}}} \min_{\mathbf x \in X_{\mathcal{S}}} \lVert \mathbf{x} - \mathbf{y} \rVert_2\right).
\end{equation}
We consider a metric inspired by \eqref{eq:classical_Chamfer} but tailored to the case of closed meshes, accounting
for the comon anatomical features of all shapes and for the common subdivision of the boundary.  

Specifically, for a pair of shapes $(\XS,\FS;l_{\mathcal S})$ and $(X_{\mathcal T},F_{\mathcal T}; l_{\mathcal T})$, 
we define a modified Chamfer distance computed separately on each subdomain, and with an additional term related to the orientation of the faces on the vessel wall:
\begin{equation}\label{eq:mesh-chamfer}
\begin{aligned}
  \mathcal{D}_{\text{Chamfer}}^*&\left((\XS,\FS),(X_{\mathcal T},F_{\mathcal T})\right) := \\
& \mathcal{D}_{\text{Chamfer}}(\XS(\Gamma_{\text{wall}}), \XT(\Gamma_{\text{wall}})) \\
& +   \lambda_n \left( 
  \frac{1}{n_{p,\mathcal{S}}}\sum_{\mathbf x\in \XS(\Gamma_{\text{wall}})}   \left(1-\lvert \mathbf{n}_{\mathbf x}\cdot\mathbf{n}^*_{(\XT(\Gamma_{\text{wall}}),{\mathbf x})} \rvert \right)
   +    \frac{1}{n_{p,\mathcal{T}}}\sum_{\mathbf x\in \XT(\Gamma_{\text{wall}})}   \left(1-\lvert \mathbf{n}_{\mathbf x}\cdot\mathbf{n}^*_{(\XS(\Gamma_{\text{wall}}),{\mathbf x})} \rvert \right)
   \right)\\
&  + \sum_{i=1}^{4} \mathcal{D}_{\text{Chamfer}}(\XS(\Gamma_i), \XT(\Gamma_i)) + \mathcal{D}_{\text{Chamfer}}(\XS(\Gamma_{\text{in}}), \XT(\Gamma_{\text{in}})) \\
& +  \sum_{i=1}^{4} \mathcal{D}_{\text{Chamfer}}(\XS(\Gamma_i\cap\Gamma_{\text{wall}}), \XT(\Gamma_i\cap\Gamma_{\text{wall}}))+\mathcal{D}_{\text{Chamfer}}(\XS(\Gamma_{\text{in}}\cap\Gamma_{\text{wall}}), \XT(\Gamma_{\text{in}}\cap\Gamma_{\text{wall}})),
\end{aligned}
\end{equation}
where $n_{p, \mathcal S}$ and $n_{p,\mathcal T}$ denote the cardinalities of the two point clouds,
$X(\Gamma)$ stands for the subset of point clouds whose vertices belong to the boundary subset $\Gamma$, 
$\mathbf{n}^*_{(X,{\mathbf x})}$ is the normal at the point of $X$ closest to $\mathbf x$, and the constant $\lambda_n>0$ is an additional hyperparameter.

\paragraph*{Transformation map} 
Following the approach introduced in section \ref{ssec:resnet-lddmm-intro}, we seek a map between source and target shapes as a diffeomorphism $\phi:[0,1]\times\mathbb{R}^3\rightarrow\mathbb{R}^3$ defined as the flow of an autonomous ordinary differential equation:
\begin{equation}\label{eq:resnet-fnn-ode}
  \dot{\mathbf x} = f(\mathbf x;\theta),\qquad \frac{d}{dt}\phi(t, \mathbf x;\theta)=f(\mathbf  x;\theta),
\end{equation}
where the vector field $f:\mathbb{R}^3\rightarrow\mathbb{R}^3$ is approximated by a feed-forward neural network (FNN) with ReLU activations as in~\cite{amor2022resnet},
and where $\theta$ represents the dependency on a generic set of hyperparameters of the FNN. We will use the abbreviation $\phi_1=\phi(1,\cdot):\mathbb{R}^3\rightarrow\mathbb{R}^3$.

In practice, the ODE \eqref{eq:resnet-fnn-ode} is discretized with a forward Euler method with $10$ time steps in the interval $I=[0,1]$ ($\Delta t = \frac{1}{10}$), resulting in a ResNet architecture~\cite{amor2022resnet} taking as input only the points of the source surface mesh $X|_{t=0}=\XS$:
\begin{equation}
  \label{eq:resnet}
  X|_{t=\Delta t\cdot (i+1)} = X|_{t=\Delta t\cdot i} + \Delta t\cdot f_{FNN}(\psi(X|_{t=\Delta t\cdot i});\theta),\quad \forall i\in\{0, \dots, N-1\}.
\end{equation} 
The map  $\psi:\mathbb{R}^3\rightarrow\mathbb{R}^{3+6\cdot n_{rff}}$, 
\begin{equation*}
  \psi(X|_{t=\Delta t\cdot i}) = (X|_{t=\Delta t\cdot i}, \{\cos(2^i\cdot X|_{t=\Delta t\cdot i}), \sin(2^i\cdot X|_{t=\Delta t\cdot i})\}_{i=0}^{7}),
\end{equation*}
is used to augment the inpus,  at each time iteration, with random fourier features~\cite{tancik2020fourier}. We used $n_{rff}=8$ in our implementation.

We employ vectorization, so the feed-forward neural network $f_{FNN}:\mathbb{R}^{3+6\cdot n_{rff}}\rightarrow\mathbb{R}^3$ that defines the vector field $f:\mathbb{R}^3\rightarrow\mathbb{R}^3$, $f(x)\mapsto f_{FNN}(\psi(x); \theta)$ is evaluated on the rows of $\psi(X|_{t=\Delta t\cdot i})\in\mathbb{R}^{n_S\times (3+6\cdot n_{rff})}$ in equation~\eqref{eq:resnet}, and $\psi:\mathbb{R}^3\rightarrow\mathbb{R}^{3+6\cdot n_{rff}}$ acts on the rows of $X|_{t=\Delta t\cdot i}\in\mathbb{R}^{n_{p, S}\times 3}$, for every $i$. 

The architecture of the FNN that we use is fixed, but its weights change for every pair of source-target aorta geometries: it has six hidden layer of dimension $500$ with ReLU activations, an input dimension of $51$ and an output dimension of $3$.

\begin{rmk}[Bijectivity of the transformation map]
In general, ResNets as defined in equation~\eqref{eq:resnet} are not invertible. 
  Using the Banach fixed point theorem, a sufficient condition to have bijectivity is 
  $\{f_{FNN}(\psi(X|_{t=\Delta t\cdot i});\theta)\}_{i\in\{0, \dots, N-1\}}$ to be 1-Lipshitz. 
  In practice, this condition can be verified as a post-processing step once the registration map has been computed, without the need of 
  additional computational costs associated with techniques, architectures, and optimization methods that enforce the invertibility exactly. The bijectivity is necessary from the theoretical point of view to define the \textit{pllback} and \textit{pushforward} of the velocity and pressure fields, and from the practical point of view to implement more robust registration algorithms, since the bijectivity property acts as an additional regularization with respect to non-rigid deformations~\cite{scarpolini2023enabling}.
\end{rmk}

\paragraph{Objective function} 
Let us denote with $\mathcal A:=\{(X, F, l)\in\mathbb{R}^{\nvert \times 3}\times\mathbb{N}^{n_f\times 3}\times\mathbb{R}^{n_{\text{cntrl}}\times 3}\}$ the set of all possible combinations of aortic shapes (vertices, faces, centerlines) in $\mathbb R^3$.

The objective function $\mathcal{L}:\mathcal{A}\times\mathcal{A}\rightarrow\mathbb{R}$ has the form
\begin{equation}\label{eq:L_AA}
\begin{aligned}
  \mathcal{L}\left((\XS,\FS, \lS),(\XT,\FT, \lT)\right)  & = 
  \mathcal{D}_{\text{Chamfer}}^*(\phi(1, \XS; \theta), \XT) + \lambda_C\cdot\sum_{i=1}^{n_{\text{cntrl}}}\lVert \phi(1, l_{\mathcal S, i}; \theta)-l_{\mathcal T, i} \rVert^2_2 \\
  & \quad + 
\lambda_{\text{edges}}  \mathcal{R}_{\text{edges}}(\phi(1, \XS; \theta),\FS) +  \lambda_{\text{en}} \mathcal{R}_{\text{energy}}(X_{\mathcal S}; \theta)
\end{aligned}
\end{equation}
and it is composed by the modified Chamfer distance \eqref{eq:mesh-chamfer}
between the mapped source $(\phi(1, \XS),\FS)$ and the target $(\XT,\FT)$, 
the distance in $L^2$-norm between the mapped source centerline $l_{\mathcal{S}}=\{l_{\mathcal{S}, i}\}_{i=1}^{390}$ and the target centerline $l_{\mathcal{T}}=\{l_{\mathcal{T}, i}\}_{i=1}^{390}$, and 
two regularizers.
The first term 
%
\begin{equation*}
  \mathcal{R}_{\text{edges}}\left(X,F\right) := 
  \sum_{(a_i,b_i,c_i) \in F} \left(
  \sum_{(e,f) \in \{(a_i, b_i), (b_i,c_i),(c_i,a_i)\} } \lVert \mathbf x_e - \mathbf x_f \rVert^2_2
   \right),
\end{equation*}
penalizes the presence of stretched edges in each face of the mesh, whilst the second term imposes the minimization of the kinetic energy along the discrete trajectory 
$\{X|_{t=\Delta t\cdot i}\}_{i=0}^{N-1}=\{\phi(\Delta t\cdot i, X_S)\}_{i=0}^{N-1}$:
\begin{equation*}
  \mathcal{R}_{\text{energy}}(\XS; \theta) := \sum_{i=0}^{N-1} \lVert f_{FNN}(\psi(X|_{t=\Delta t\cdot i});\theta)\rVert^2_2.
\end{equation*}

The constants $\lambda_C>0$, $\lambda_{\text{en}}$, $\lambda_{\text{edges}}$ are positive parameters. Notice that in the above definitions we have used the fact that, when applying the diffeomorphism $\phi$, only the point clouds (coordinates of the mesh vertices) are mapped, whilst the set of faces $\FS$ remains the same.

\begin{problem}[Shape registration with ResNet-LDDMM]
With the above definitions, we formulate the following discrete surface registration problem. Given a source and a target meshes 
$\mathcal S=(\XS,\FS, \lS)$ and $\mathcal T = (\XT,\FT, \lT)$, find 
  \label{def:resnetlddmm}
  \begin{align*}
    &\argmin_{\theta}\ 
    \mathcal{D}_{\text{Chamfer}}^*(\phi(1, \XS; \theta), \XT) + \lambda_C\cdot\sum_{i=1}^{n_{\text{cntrl}}}\lVert \phi(1, l_{\mathcal S, i}; \theta)-l_{\mathcal T, i} \rVert^2_2 \\
  & \quad + 
 \lambda_{\text{edges}} \mathcal{R}_{\text{edges}}(\phi(1, \XS; \theta),\FS) +  \lambda_{\text{en}} \mathcal{R}_{\text{energy}}(\{\phi(\Delta t\cdot (i+1), A_S;\theta)\}_{i=0}^{N-1})\\
  \end{align*}
such that   $\phi(t,\cdot;\theta)$ is the flow of the discretized ODE defined by the corresponding ResNet vector field:
\begin{align*}
    &\phi(\Delta t\cdot (i+1), \XS;\theta) = X|_{t=\Delta t\cdot i} + \Delta t\cdot f_{FNN}(\psi(X|_{t=\Delta t\cdot i});\theta),\quad \forall i\in\{0, \dots, N-1\},\\
    &\phi(0, X_S;\theta) = \XS.
\end{align*}
\end{problem}

\paragraph{Optimization} To solve Problem \eqref{def:resnetlddmm}, we consider the ADAM optimizer~\cite{kingma2014adam} combined with a \textit{multigrid} strategy, i.e., considering
three level of refinement for the source mesh, while the target mesh is kept fixed. 
This approach is crucial to guarantee 
the convergence of the discrete registration problem: on one hand it speeds up the procedure and on the other hand guarantees an arbitrary small registration error. 
In practice, for a source mesh $(\XS,\FS)$, we will denote as
\begin{equation*}
(\XS^{j}, \FS^{j})\in\mathbb{R}^{\nvert^{j}\times 3}\times\mathbb{N}^{n_f^{j}\times 3},\qquad j\in\{0, 1 , 2\},
\end{equation*}
with an increasing number of vertices $\nvert^{j=0} < \nvert^{j=1} < \nvert^{j=2}$ the three considered refinements, assuming that the set of faces is consistently refined.
The different refinements are obtained imposing a different upper bound for the lengths of face edges and an upper bound for the radii of the surface Delaunay balls.

An analogous \textit{multigrid} approach has been proposed in~\cite{scarpolini2023enabling}. However, since both the transformation map 
$\phi(1, \cdot):\mathbb{R}^3\rightarrow\mathbb{R}^3$ and the vector field $f:\mathbb{R}^3\rightarrow\mathbb{R}^3$ act on the ambient space $\mathbb{R}^3$, we do not need to interpolate from one source mesh refinement to the other. 

\paragraph{Registration of the computational domain: \textit{pullback} and \textit{pushforward} operators}
For each couple of source and target meshes, we store the registration map as the image of the source point cloud, i.e.,  $\phi_1(X_{\mathcal{S}})\subset\mathbb{R}^3$. 
We denote with $\Omega_{\mathcal{S}}\subset\mathbb{R}^3$ and $\Omega_{\mathcal{T}}\subset\mathbb{R}^3$ the computational domains used for the CFD simulations for the source/template and target geometries, respectively. The surface registration maps are computed and interpolated on the source domain $\Omega_{\mathcal{S}}$ 
through RBF interpolation $\phi_{RBF}:\mathbb{R}^3\rightarrow\mathbb{R}^3$ with thin splines as kernels and $\phi_1(X_{\mathcal{S}})$ as centers.
Let $g_{\mathcal S}:\Omega_{\mathcal{S}}\rightarrow\mathbb{R}$ be a function defined on the source domain. The \textit{\textit{pushforward}} of $g_{\mathcal S}$ through the registration map
is defined as
\begin{equation}\label{eq:pushforward}
\phi_{RBF}^{\#}(g_{\mathcal S}):\Omega_{\mathcal{T}}\rightarrow\mathbb{R},\quad  \phi_{RBF}^{\#}(g_{\mathcal S}):=g_{\mathcal S} \circ\phi_{RBF}^{-1},
\end{equation}
Conversely, for a function  $g_{\mathcal T}:\Omega_{\mathcal{T}}\rightarrow\mathbb{R}$, the \textit{pullback} $(\phi_{RBF})_{\#} \left(g_{\mathcal T}\right):\Omega_{\mathcal{S}}\rightarrow\mathbb{R}$ is defined as:
\begin{equation}\label{eq:pullback}
  (\phi_{RBF})_{\#} \left(g_{\mathcal T}\right) :=\left(g_{\mathcal T} \right) \circ\phi_{RBF}.
\end{equation}

\subsection{Shape registration results}
To train the ResNet, the source (or template) mesh is chosen in the training set of shapes and kept fixed throughout the offline stage, 
whilst the target mesh varies among the remaining $723$ training shapes. 
In the online stage, the source is unchanged whilst the target mesh varies among the $52$ test geometries. 

\begin{rmk}[Choice of template geometry]
 The source mesh has been chosen within the training set without any particular criteria. In general, the selection can also be optimized with respect to the reconstruction error of the velocity and pressure fields, see section~\ref{sec:sml}.  
\end{rmk}

The hyperparameters for the regularization in \eqref{eq:L_AA} are chosen as $\lambda_n=5\cdot 10^{-5}$, $\lambda_C = 10^{-5}$, 
and   $\lambda_{\text{en}}=\lambda_{\text{edges}}=1$, for the terms related to face orientation, centerline, energy of the trajectory, and edges, respectively.

We consider a total number of epochs $n_{\text{epochs}}=5000$. The change of source mesh is performed at the epochs $e_{0,1} = 3000$, 
from $\left(\XS^{0}, \FS^{0}\right)$ to $\left(\XS^{1}, \FS^{1}\right)$, and  $e_{1,2} = 4000$, from $\left(\XS^{1}, \FS^{1}\right)$ 
to $\left(\XS^{2}, \FS^{2}\right)$. The source meshes have a number of vertices $n^{j=0}_v=4127$, $n^{j=1}_v=11402$ and  $n^{j=2}_v=110676$.
These values are kept fixed. However, one could also use an adaptive strategy to select the refinement epochs, e.g., based on the ratio between Chamfer distances at consecutive steps.

\begin{figure}[!htp]
  \centering
  \includegraphics[width=0.95\textwidth]{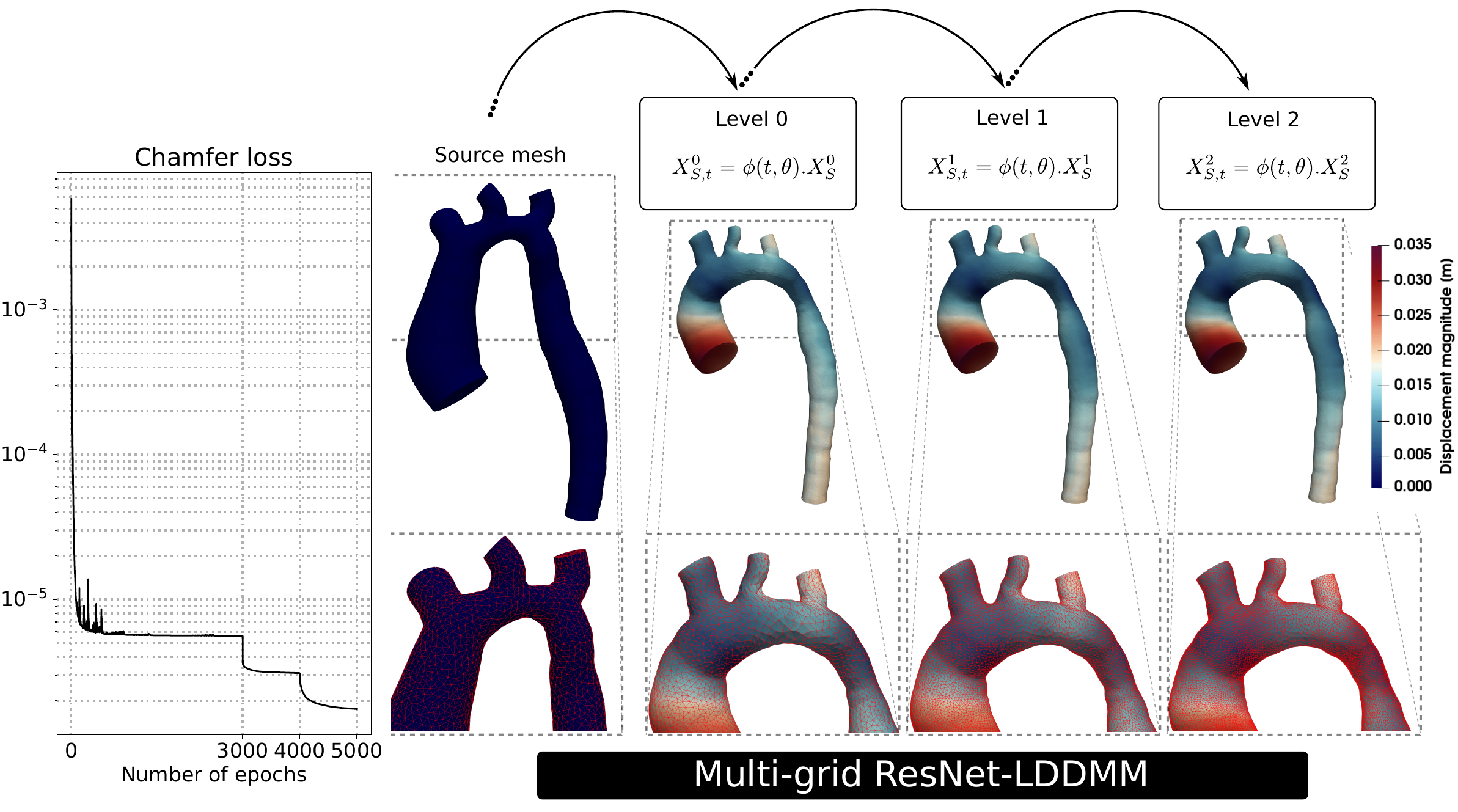}
  \caption{Application of the \textit{multigrid} ResNet-LDDMM. The source surface mesh is refined progressively during the training to guarantee the convergence of the discrete registration problem. 
\textbf{Left}: Decay of the Chamfer loss over the epochs. At epochs $3000$ and $4000$ the mesh is refined. 
\textbf{Right}: Displacement field $\phi(1, X_\mathcal{S}^i;\theta)\text{-}X_\mathcal{S}$ and surface mesh for refinement levels $j\in\{0, 1, 2\}$.}
  \label{fig:multigrid}
\end{figure}

\begin{figure}[!htp]
  \centering
  \includegraphics[width=0.83\textwidth]{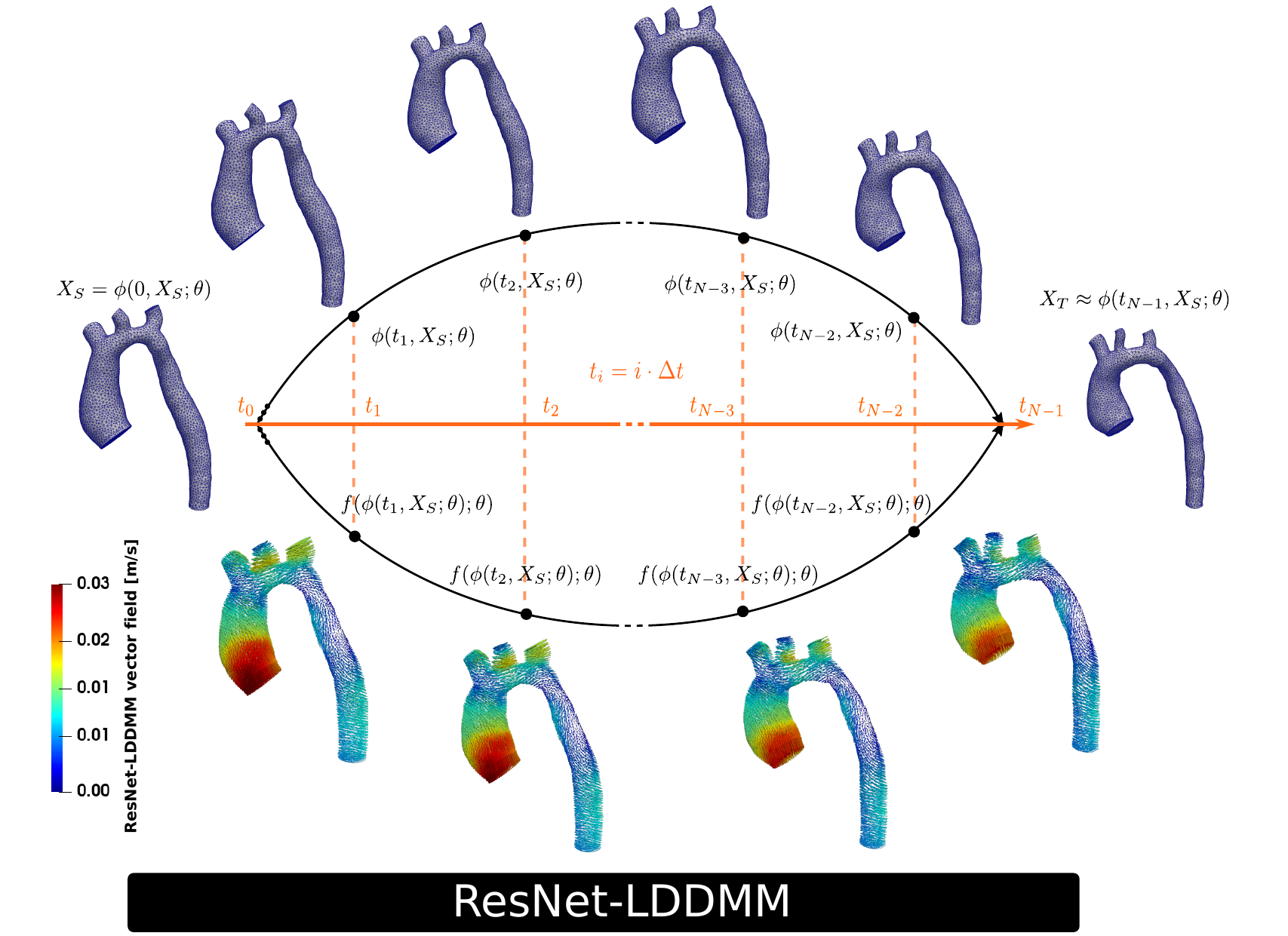}
  \caption{Registration of a source surface mesh $X_\mathcal{S}$ onto the target surface mesh $X_\mathcal{T}$ through the ResNet-LDDMM vector field $f$. 
  \textbf{Top}: registration field $\phi(t_i, \mathbf x;\theta)$ at different intermediate steps $t_i$. \textbf{Bottom}: Vector field $f(\mathbf x;\theta)$ in the configuration $\phi(t_i, X_\mathcal{S};\theta)$.}
  \label{fig:flow}
\end{figure}

\begin{figure}[!htp]
  \centering
  \includegraphics[width=0.95\textwidth, trim={0 0 0 20}, clip]{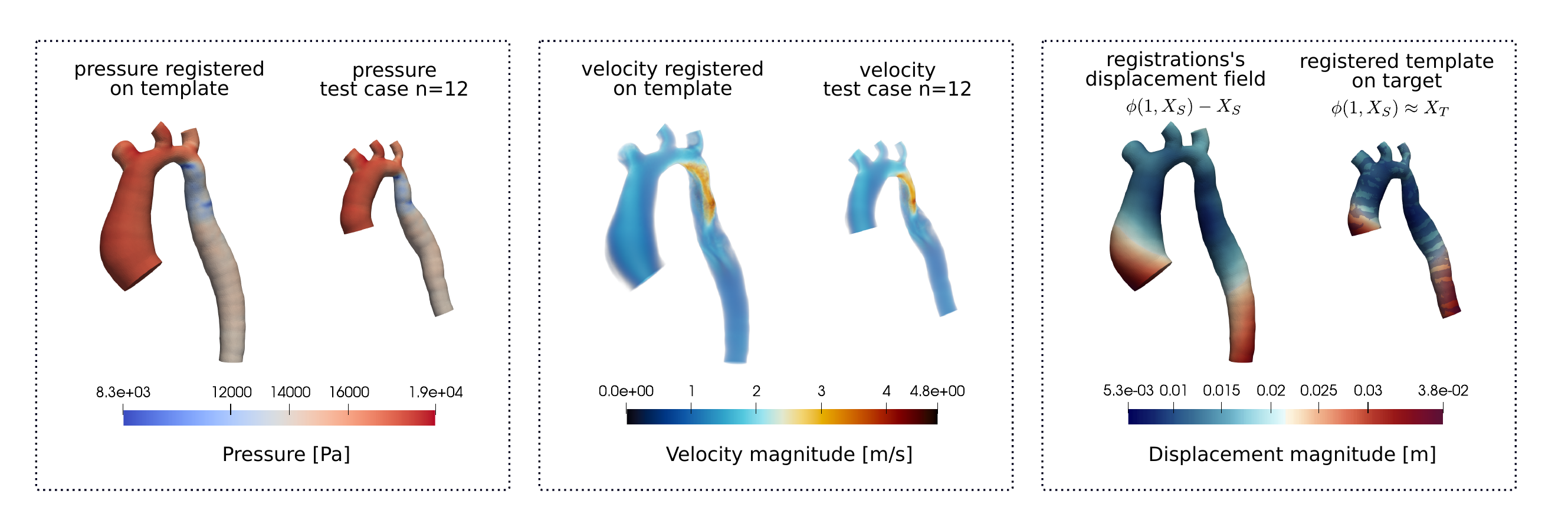}\\
  \includegraphics[width=0.95\textwidth, trim={0 0 0 20}, clip]{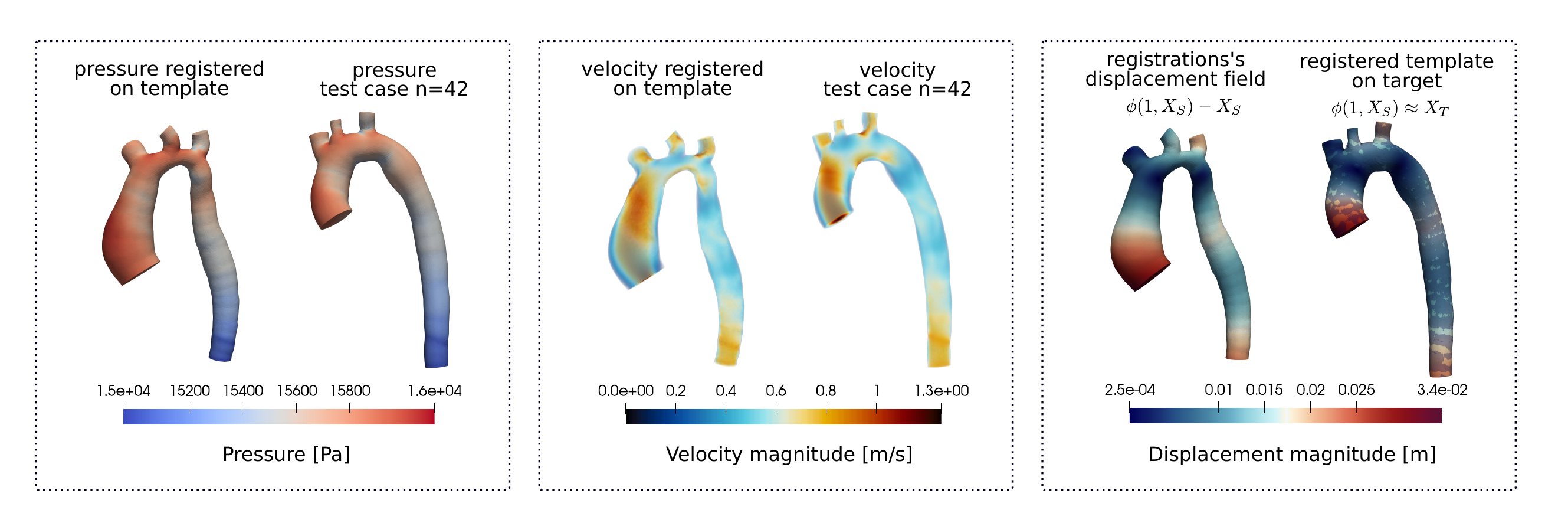}
  \caption{Registration of test case $n=12$ (top row, best) and $n=42$ (bottom row, worst) velocity and pressure fields onto the template. The displacement fields $\phi(1,X_\mathcal{S};\theta)\text{-}X_\mathcal{S}$ are shown on the third column: the registered geometry $\phi(1,X_\mathcal{S})$ is shaded and compared with original template $X_\mathcal{T}$, under the title \textit{registered template on target}.}
  \label{fig:12_42}
\end{figure}

A sketch of the \textit{multigrid} optimization is shown in figure~\ref{fig:multigrid}, displaying both the mesh refinement on the aortic arch and the 
displacement field $\phi(1, X_S^i;\theta)\text{-}X_S^i$ for the different refinement levels. Figure~\ref{fig:multigrid} (left) shows also 
the loss decay on a sample geometry, highlighting the influence of the \textit{multigrid} strategy for convergence.
Figure~\ref{fig:flow} depicts different steps of the registration process between two surface meshes showing
both the registration field $\phi(t_i, X_S;\theta)$ and the vector field $f(x;\theta)$ at different intermediate deformed configurations. See figure~\ref{fig:12_42} for an example of velocity and pressure fields' registration on the template geometry at systolic peak, corresponding to the best $n=12$ and worst $n=42$ test cases from figure~\ref{fig:cluster_v}. 

For validating the registration algorithm, we evaluate the classical Chamfer distance~\eqref{eq:classical_Chamfer} between the point clouds of the 
registered and the target geometries, normalized by the diameter of the target geometry, for each considered shape. The results, shown in Table \ref{tab:reg-chamfer}, confirm the 
robustness and the accuracy of the registration across the whole database.
The computational cost for registering the $723$ training geometries was of circa $\SI{114}{\hour}\approx 723\times\SI{9.45}{\minute}$ (embarrassingly parallel tasks), while the online registration of the source on a single
target required, on average, $9.45$ minutes.
\begin{table}[H]
    \centering
    \begin{tabular}{lcc}
         & \textbf{Training set} (n=723) & \textbf{Test set} (n=52) \\
         \hline
        Average Chamfer Distance & $0.00367$ & $0.00347$ \\
        \hline
        Maximum Chamfer Distance & $0.00605$ & $0.00470$ \\
        \hline
    \end{tabular}
    \caption{Chamfer distance~\eqref{eq:classical_Chamfer} between the registered source and target shapes, normalized by the diameter of the target mesh.}
    \label{tab:reg-chamfer}
\end{table}

%

\section{Application of shape registration for solution manifold learning}
\label{sec:sml}

We formally refer to \textsl{solution manifold} as the set of time-dependent velocity and pressure fields which are solutions of the Navier--Stokes equations~\eqref{eq:3dnse}
for different computational domains and boundary conditions~\eqref{eq:3dnse-bc}. 
The goal of solution manifold learning is to accurately and efficiently approximate the solution manifold using the available training data. 
Widely used approaches consider linear global reduced bases of the discrete finite element spaces 
(see e.g. ~\cite{benner_model_2017,hesthaven2016certified,rozza2022advanced}), time- or geometry-dependent linear bases, as well as 
non-linear approximants, including autoencoders~\cite{Fresca2020, romor2023nonlinear, ROMOR2025113729} or other non-linear dimension reduction techniques from machine learning. 

This section presents two applications of the shape registration method from section \ref{sec:registration} to solution manifold learning. 
First, in subsection~\ref{subsec:sml_correlations}, we propose different metrics to study the correlation between geometries and solutions (pressure/velocity), as well 
as between velocity and pressure, based on mapping the dataset of patient-specific flow data on the same reference shape.
The results are used in the context of data assimilation problems,  where correlations are particularly relevant when considering the estimation of pressure-related quantities from velocity measurements (see sections \ref{sec:da} and \ref{sec:prec}).  
Next, in section~\ref{subsec:sml_rec}, we investigate the accuracy of the global and local rSVD bases constructed registering the snapshots of pressure and velocity solutions from different shapes onto the same reference.

\subsection{Analysis of physics-based and geometric-based correlations}
\label{subsec:sml_correlations}
We evaluate the correlation between dissimilarity matrices computed from the database of shapes and corresponding numerical solutions, exploiting
the fact that all relevant fields can be encoded in the same discrete space, i.e., the finite element space on the computational domain of the reference shape.
In what follows, let us denote with $\mathcal{S}$ the fixed shape, and with $n_{p, \mathcal{S}}$, $\dofu$, and $\dofp$ the number of vertices of the
corresponding finite element mesh and the degrees of freedom of the velocity and pressure solutions, respectively ($n_{p, \mathcal{S}}=110676$, $\dofu = 332028$, $\dofp = 110676$ for the particular selection of the reference mesh considered in this study).
The remaining \textit{target} geometries, that can be registered to $\mathcal{S}$ via the registration maps
$\phi_1^i$, will be denoted
as $\shape{i}$, $i=1,\hdots,n_{\rm geo}$, with $n_{\rm geo} = 724 + 52 = 776$ (both the test and training sets used in section  \ref{sec:registration}).
Each target shape $\shape{i}$ can be encoded using the three spatial coordinates of the reference mesh vertices mapped to the target domain i.e. through the displacement fields
\begin{equation*}
  Y_{\shape{i}} = \phi^i_1(\XS)-\XS, \qquad Y_{\shape{i}}\in\mathbb{R}^{\nvertref \times 3}, 
\end{equation*}
and a further coordinate based on the distance from the centerline $l_{\shape{i}}$ computed on the mapped reference vertices
\begin{equation}
  \label{eq:enc_geo}
  Z_{\shape{i}} = d\left(\phi^i_1\left(\XS\right),l_{\shape{i}} \right),\qquad Z_{\shape{i}} \in \mathbb{R}^{\nvertref}\,,
\end{equation}
where $d$ stands for the Euclidean distance.

For each $\shape{i}$, let us now define the matrices $X^{\mathbf u}_{\shape{i}} \in\mathbb{R}^{\dofu \times n_T}$  and 
$X^{p}_{\shape{i}}\in\mathbb{R}^{\dofp \times n_T}$ containing $n_T=80$ equally spaced snapshots of the velocity and pressure fields in the time interval $[0.05s, 0.25 s]$ (see remark~\ref{rmk:timewindow}) registered on $\mathcal{S}$ and evaluated on the reference finite element mesh.

We then introduce two metrics. Let $A$ and $B$ be two matrices, $A, B \in O(D, r)$  with $r>0$ orthonormal columns of dimension $D>0$. 
Following ~\cite{galarce2022state}, we define the Hausdorff distance as
\begin{equation}\label{eq:d_hausdorff}
 d_H^2(A, B) := \max\left(\max_{a\in \text{col}(A)} \frac{\lVert a-P_{B}a\rVert^2}{\lVert a\rVert^2};\max_{b\in \text{col}(B)} \frac{\lVert b-P_{A}b\rVert^2}{\lVert b\rVert^2}\right),
\end{equation}
where $\text{col}(A), \text{col}(B)$ are the set of columns of the matrices $A$ and $B$, respectively, 
and $P_C$, ($C=A,B$), is the finite-dimensional orthogonal projector onto the space spanned by the columns of $C$.
We also define the Grassmann distance (see e.g.~\cite{daniel2020model}) as
\begin{equation}\label{eq:d_grassman}
  d_{\text{Gr}}^2(A, B) = \sum^{r}_{i=1}\arccos^2(\sigma_i),\qquad A^TB = 
  U\begin{bmatrix}
    \sigma_{1} & & \\
    & \ddots & \\
    & & \sigma_{r}
  \end{bmatrix}V^T,
\end{equation}
where $U\Sigma V^T$ stands for the singular value decomposition of $A^TB$.

We consider then the following \textit{dissimilarity} matrices of dimension $\mathbb{R}^{n_{\text{geo}}\times n_{\text{geo}}}$:
\begin{equation}\label{eq:geo_matrices}
\begin{aligned}
(K^{\text{enc}})_{ij} & = d\left(Z_{\shape{i}}, Z_{\shape{j}}\right), \\
(K^{\phi})_{ij} & = d\left(\phi(1, X_{\mathcal S_i})-X_{\mathcal S_i}, \phi(1, X_{\shape{j}})-X_{\shape{j}}\right) = d\left(Y_{\shape{i}}, Y_{\shape{j}}\right), \\
\end{aligned}
\end{equation}
based on the Euclidean distance between geometric encodings, and 
\begin{equation}\label{eq:up_matrices}
\begin{aligned}
(K_{H}^{\mathbf u})_{ij} & = d_H(X_{\shape{i}}^{\mathbf u}, X_{\shape{j}}^{\mathbf u}),
\quad(K_{{\rm Gr}}^{\mathbf u})_{ij} = d_{\rm Gr}(X_{\shape{i}}^{\mathbf u}, X_{\shape{j}}^{\mathbf u}),\\
(K_{H}^{p})_{ij} & = d_H(X^p_{\shape{i}}, X^p_{\shape{j}}),
\quad(K_{{\rm Gr}}^{p})_{ij} = d_{\rm Gr}(X^p_{\shape{i}}, X^p_{\shape{j}}),
\end{aligned}
\end{equation}
based on the distances between the solutions introduced in \eqref{eq:d_hausdorff} and \eqref{eq:d_grassman}.
These matrices are used to evaluate the correlation between the geometry and the solution, as well as the velocity and the pressure fields using a Mantel test with 
Pearson’s product-moment correlation coefficient $r_m\in[-1, 1]$ and $999$ permutations.

\begin{figure}[!htp]
  \centering
  \includegraphics[width=0.59\textwidth]{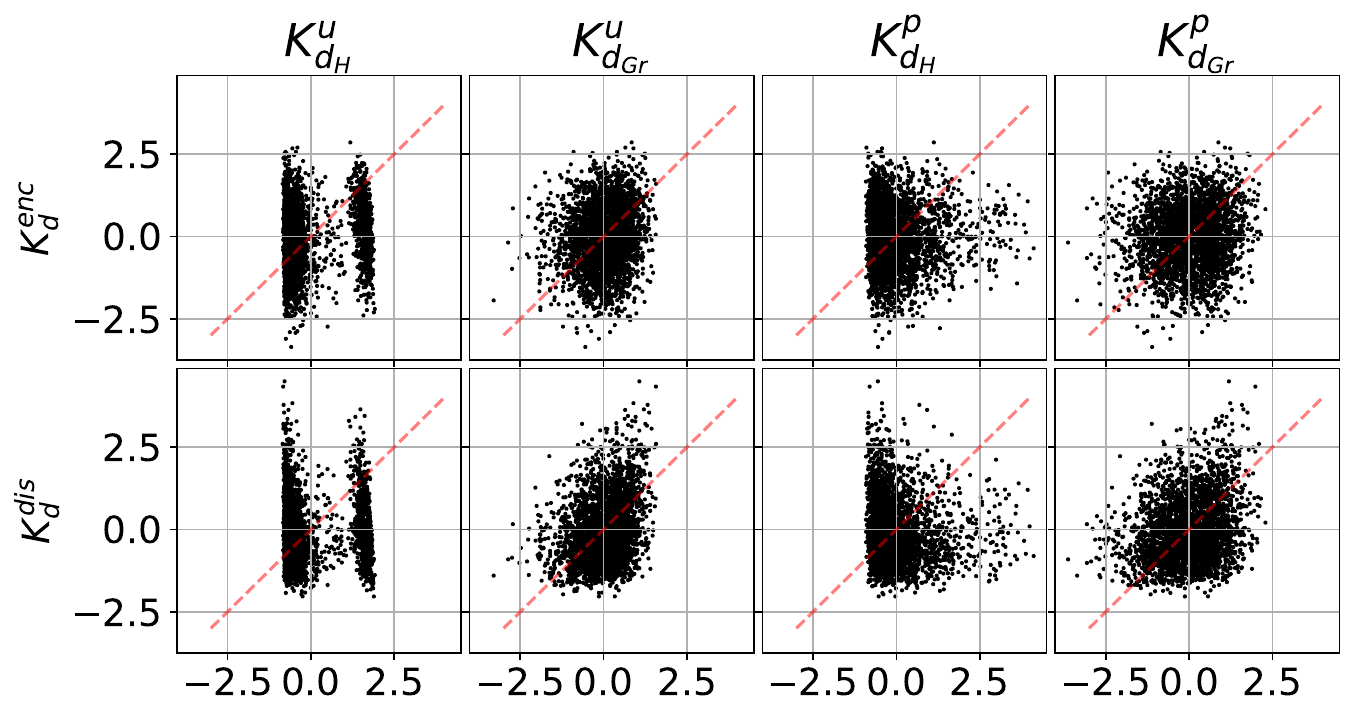}
  \includegraphics[width=.335\textwidth]{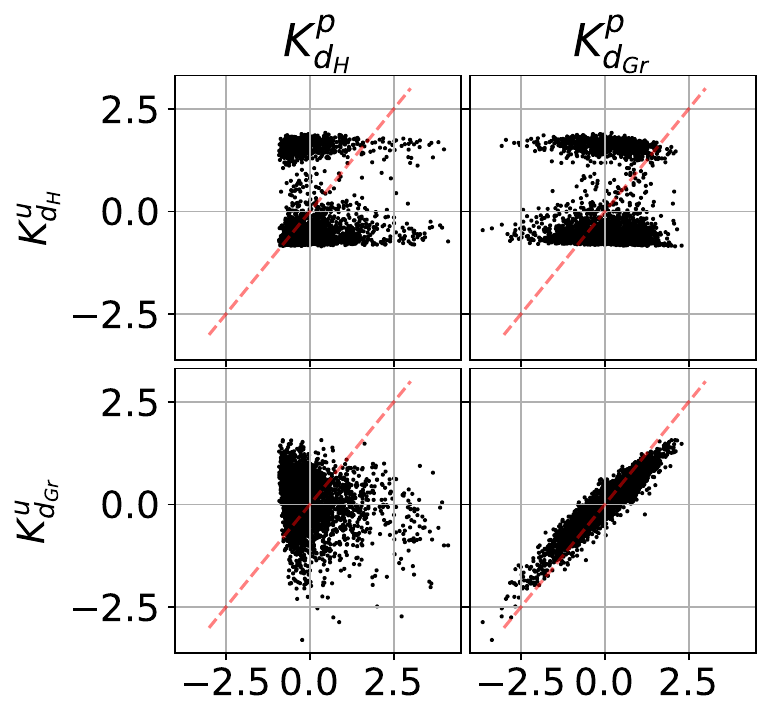}
  \caption{Correlation among dissimilarity matrices: each dot corresponds to an entry of the matrices indicated in the $x$- and $y$-axis. 
Diagonal entries have been omitted.  Only one every $100$ entries among the $n_{\text{geo}}^2\text{-}n_{\text{geo}}$ off-diagonal entries are shown. 
\textbf{Left}: correlation between geometry encoding and velocity/pressure fields. \textbf{Right}: Correlation between velocity and pressure solutions.}
  \label{fig:mantel}
\end{figure}

Table~\ref{tab:mantel1} shows the results for the correlation between geometry and velocity/pressure fields, whilst Table~\ref{tab:mantel2} shows the results for the correlation 
between velocity and pressure.
Since different metrics are used, the dissimilarity matrices are centered, before the correlation coefficient is computed. 
For a qualitative comparison, figure~\ref{fig:mantel} shows the dissimilarity matrices entries omitting the diagonal ones. 

\begin{table}[htp!]
  \centering
  \footnotesize
  \caption{Mantel test results for the correlation between geometry dissimilarity matrices \eqref{eq:geo_matrices} and velocity/pressure dissimilarity matrices \eqref{eq:up_matrices}. Correlation coefficients with less statistical significance ($p$-value$>0.05$) have been omitted. }
  \begin{tabular}{
      l 
      |>{\centering\arraybackslash}p{3.cm} 
      |>{\centering\arraybackslash}p{3.cm} 
      |>{\centering\arraybackslash}p{3.2cm} 
      |>{\centering\arraybackslash}p{3.cm} }
      \textbf{} &$K_{H}^{\mathbf u}$ & $K^{\mathbf u}_{\text{Gr}}$ & $K^p_{d_H}$ &$K^p_{d_{\text{Gr}}}$ \\[3pt]
      \hline
      \hline 
     $K^{\text{enc}}$ & $r_m=0.172, p=0.001$ & $r_m=0.156, p=0.001$ & - & $r_m=0.110, p=0.001$ \\
      \hline
    $K^{\phi}$ & - & $r_m=0.267, p=0.001$ & $r_m=-0.094, p=0.001$ & $r_m=0.217, p=0.001$ \\
      \hline
  \end{tabular}
  \label{tab:mantel1}
\end{table}
\begin{table}[htp!]
  \caption{Mantel test results for the correlation between the velocity and pressure dissimilarity matrices \eqref{eq:up_matrices}. Correlation coefficients with less statistical significance ($p$-value$>0.05$) have been omitted. }
  \centering
    \footnotesize
  \begin{tabular}{
      l 
      |>{\centering\arraybackslash}p{3.5cm} 
      |>{\centering\arraybackslash}p{3.5cm} }
     & $K^p_{H}$ & $K^p_{\text{Gr}}$\\[3pt]
    \hline
    \hline
    $K_{H}^{\mathbf u}$ & $r_m=0.433, p=0.033$ & -\\[2pt]
    \hline
    $K^{\mathbf u}_{\text{Gr}}$ & $r_m=-0.169, p=0.001$ & $r_m=0.943, p=0.001$ \\[2pt]
    \hline
\end{tabular}
  \label{tab:mantel2}
\end{table}
The results suggests that the Grassmann metric is more granular than the Hausdorff metrics proposed in~\cite{galarce2022state}. 
In particular, the correlation plots $K^{\text{enc}}$ \textit{vs.} $K^{\mathbf{u}}_{H}$, $K^{\phi}$ \textit{vs.} $K^{\mathbf{u}}_{H}$, 
$K^{\mathbf{u}}_{H}$ \textit{vs.} $K^{p}_{H}$, and $K^{\mathbf{u}}_{H}$ \textit{vs.} $K^{p}_{\text{Gr}}$ show two clusters of more and less correlated pairs. 
We can also observe that there exists a choice of metrics ($K^{\mathbf{u}}_{\text{Gr}}$ \textit{vs.} $K^{p}_{\text{Gr}}$) for which velocity and pressure fields 
are highly correlated, suggesting that in this setting the data assimilation of pressure from velocity measurements may be more feasible. 
At the same time, based on the correlation study, inferring the velocity or the pressure solution solely from the geometry seems to be a more challenging task.

\begin{figure}[!htp]
 \centering
 \includegraphics[width=0.9\textwidth]{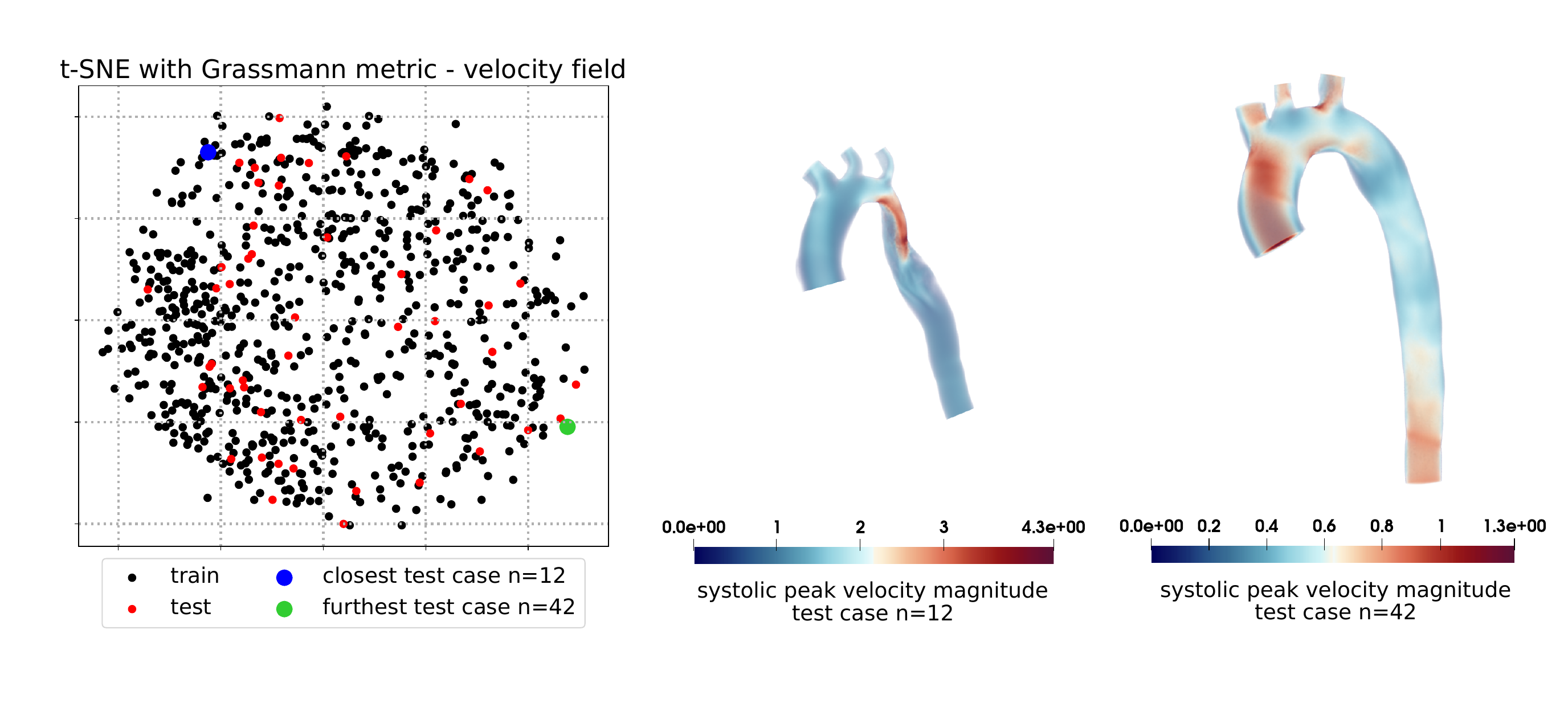}
 \caption{\textbf{Left: } Clustering with MDS of the $724$ training and $52$ test geomeries. \textbf{Right: }  Test case $12$ and $42$ represent the closest and furthest geometries to the training set with respect to the Grassmann distance on the velocity field.}
 \label{fig:cluster_v}
\end{figure}

In figure~\ref{fig:cluster_v}, we show the clustering of the available training and test geometries with MDS and dissimilarity matrix $K^{u}_{d_{G_r}}$ from equation~\eqref{eq:up_matrices}. In particular, we can detect the test geometries with the best ($n=12$) and worst ($n=42$) approximable velocity field.

\subsection{Approximation properties of  global and local rSVD bases}
\label{subsec:sml_rec}
The shape dataset has been split into a  training ($n_{\text{train}}=724$) and a test  ($n_{\text{test}}=52$) set. To obtain a global linear reduced basis, we apply randomized SVD~\cite{halko2011finding} with a given rank $r>0$ to the matrices of velocity and pressure fields on the template geometry ordered column-wise $X^{\mathbf u}_{\text{train}}\in\mathbb{R}^{\dofu\times(n_{\text{train}}n_T)}$ and $X^{p}_{\text{train}}\in\mathbb{R}^{\dofp\times(n_{\text{train}}n_T)}$, respectively: 
\begin{equation*}
\begin{aligned}
X^{\mathbf u}_{\text{train}} \ &\rsvd\  \Phi_{\mathbf u}\Sigma^{r}_{\mathbf u}\Psi_{\mathbf u}
,\quad \Phi_{\mathbf u} \in\mathbb{R}^{\dofu\times r},\ \Sigma^{r}_{\mathbf u}\in\mathbb{R}^{r\times r},\ \Psi_{\mathbf u}\in\mathbb{R}^{r\times (n_{\text{train}}n_T)},\\
X^{p}_{\text{train}} \ &\rsvd\  \Phi_p\Sigma^{r}_p\Psi_p,\quad \Phi_p\in\mathbb{R}^{\dofp\times r},\ \Sigma^{r}_p\in\mathbb{R}^{r\times r},\ \Psi_p\in\mathbb{R}^{r\times (n_{\text{train}}n_T)}.
\end{aligned}
\end{equation*}
The columns of the matrices $\Phi_{\mathbf u}$ and $\Phi_p$ define the orthonormal global rSVD basis. For a $r$-dimensional (reduced) representation of a velocity field
$z_{\mathbf u}^{r} \in \mathbb R^r$ (resp. of a pressure field $z_p^{r}\in\mathbb{R}^r$ ), the corresponding approximation in the full finite element space will be defined by $\Phi_{\mathbf u} z_{\mathbf u}^{r}$ (resp. $\Phi_p z_p^{r}$).

\begin{rmk}[Partitioned \textit{vs.}  monolithic rSVD]
 We considered a partitioned rSVD global basis, i.e., computing the rSVD modes for velocity and pressure from two snapshot matrices. An alternative monolithic approach consists in 
computing the rSVD on a single snapshot matrix of dimension $(\dofu+\dofp) \times n_{\rm train} n_T$ where velocity and pressure solutions are stacked row-wise.
In our case, the choice was dictated by the better performance in terms of accuracy of the partitioned rSVD basis.
However, especially in the context of data assimilation for inferring pressure fields from velocity observations, a monolithic rSVD might have the advantage of handling the coupled latent representation in a single $r$-dimensional variable, which allows to automatically obtain the pressure field from the same reduced variable~\cite{galarce2023displacement}. 
\end{rmk}

Depending on the reduced dimension $r$, we consider the relative $L^2$-reconstruction errors to evaluate the accuracy of the reduced approximation 
\begin{equation}
  \label{eq:rec}
  \epsilon^{r}_{u} (\mathbf u_i(t)) := \frac{\lVert \mathbf u_i(t) -\Phi_{\mathbf u}\Phi_{\mathbf u}^T \mathbf u_i(t)\rVert_2}{\lVert \mathbf u_i(t) \rVert_2},\quad 
  \epsilon^{r}_p (p_i(t)) := \frac{\lVert p_i(t)-\Phi_p\Phi_p^T p_i(t)\rVert_2}{\lVert p_i(t)-\bar{p}_i \rVert_2},
\end{equation}
varying $u_i$ and $p_i$ among the numerical solutions of the training and test geometries, for $i=1,\hdots,724 + 52$, and time instances $t\in \{0.05s+n\cdot\Delta t\ | n\in\{0,\dots, n_T=80\}\}$.
In \eqref{eq:rec}, $\bar{p}_i\in\mathbb{R}$ stands for the average of the considered pressure solution.

\begin{figure}[!htp]
  \centering
  \includegraphics[width=0.85\textwidth]{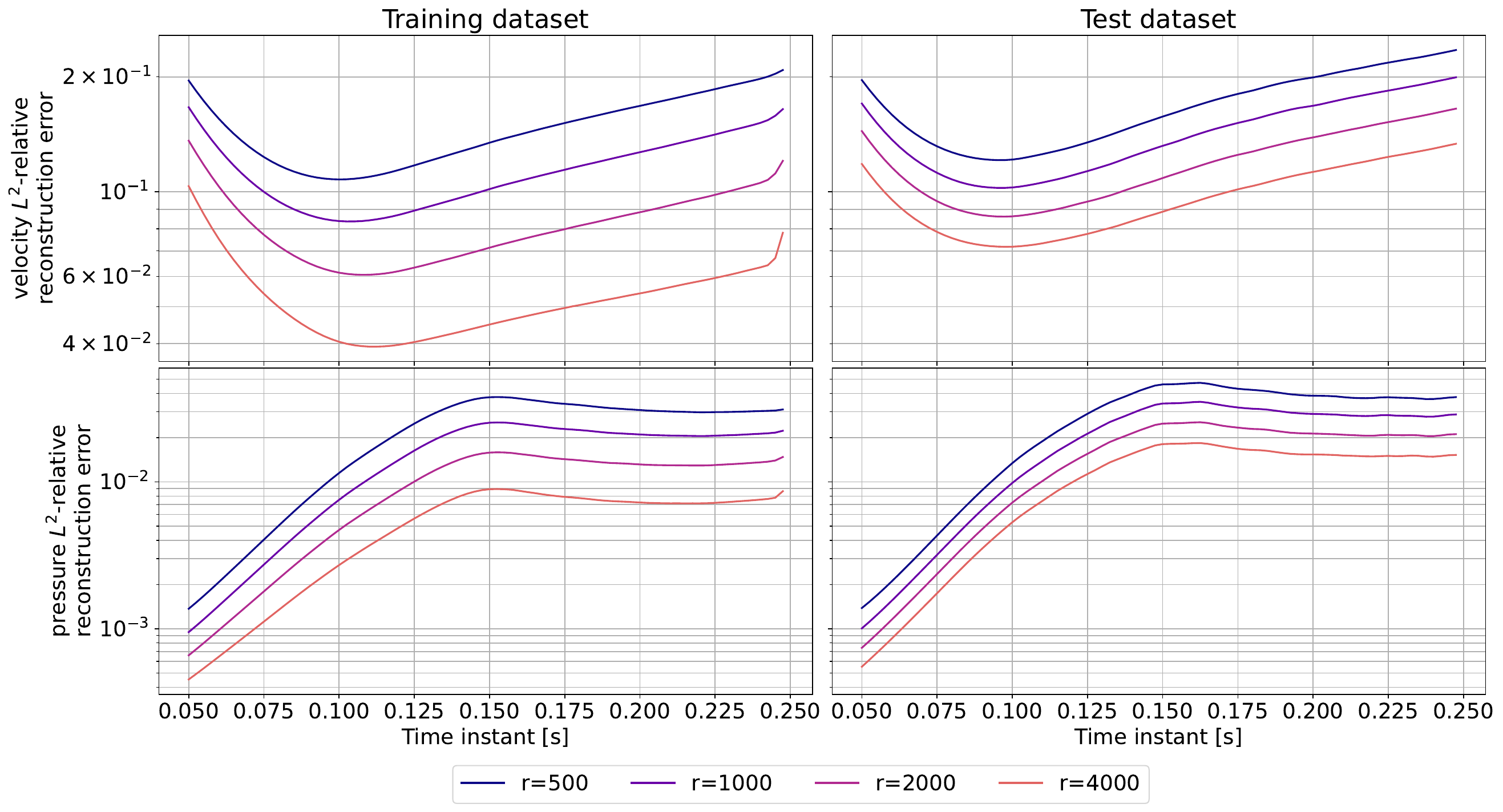}
  \caption{Average among the training or test datasets, of the relative $L^2$-reconstruction errors \eqref{eq:rec} of the velocity 
  and pressure fields in the time interval $[0.05s,0.25]$ for different rSVD ranks.}
  \label{fig:recerr}
\end{figure}

The relative $L^2$-reconstruction errors \eqref{eq:rec} are shown in figure~\ref{fig:recerr} for $r\in\{500, 1000, 2000, 4000\}$, showing that a very high number of modes
is required to obtain approximation errors of the order of 10\% for the velocity field. The error for the pressure field is of the order of 1-5\% for all considered dimensions.
\begin{rmk}[Time window]
  \label{rmk:timewindow}
  Outside the considered time window $t \in [0.05s,0.25s]$, the velocity was poorly approximated. This might be due to the additional complexity in the flow patterns during flow deceleration, such as arise of vortices, which depend very strongly on the geometrical details and are not accurately reprodicible with a linear basis such as $\Phi_{\mathbf u}$. This is the reason why we restrict our data assimilation studies to the time window $t \in [0.05s,0.25s]$.
\end{rmk}

The current results suggest that the size of the reduced space might not be suitable for designing reduced order models, restricting the finite element spaces of the variational formulation to the linear space spanned by the rSVD basis.
Although this approach has been proposed and applied in related contexts using much lower reduced space dimensions~\cite{guibert2014group,PEGOLOTTI2021113762}, the additional complexity
considered in the current setting (general geometries, high variability of outlet boundary dimension, variable Windkessel parameters depending on flow split and measured inlet flow rate, presence of 
stenosis which leads to more complex patterns, and the incorporation of turbulence modelling) leads to the need of a much larger space for satisfactory approximations.

\begin{figure}[!htp]
  \centering
  \includegraphics[width=0.7\textwidth]{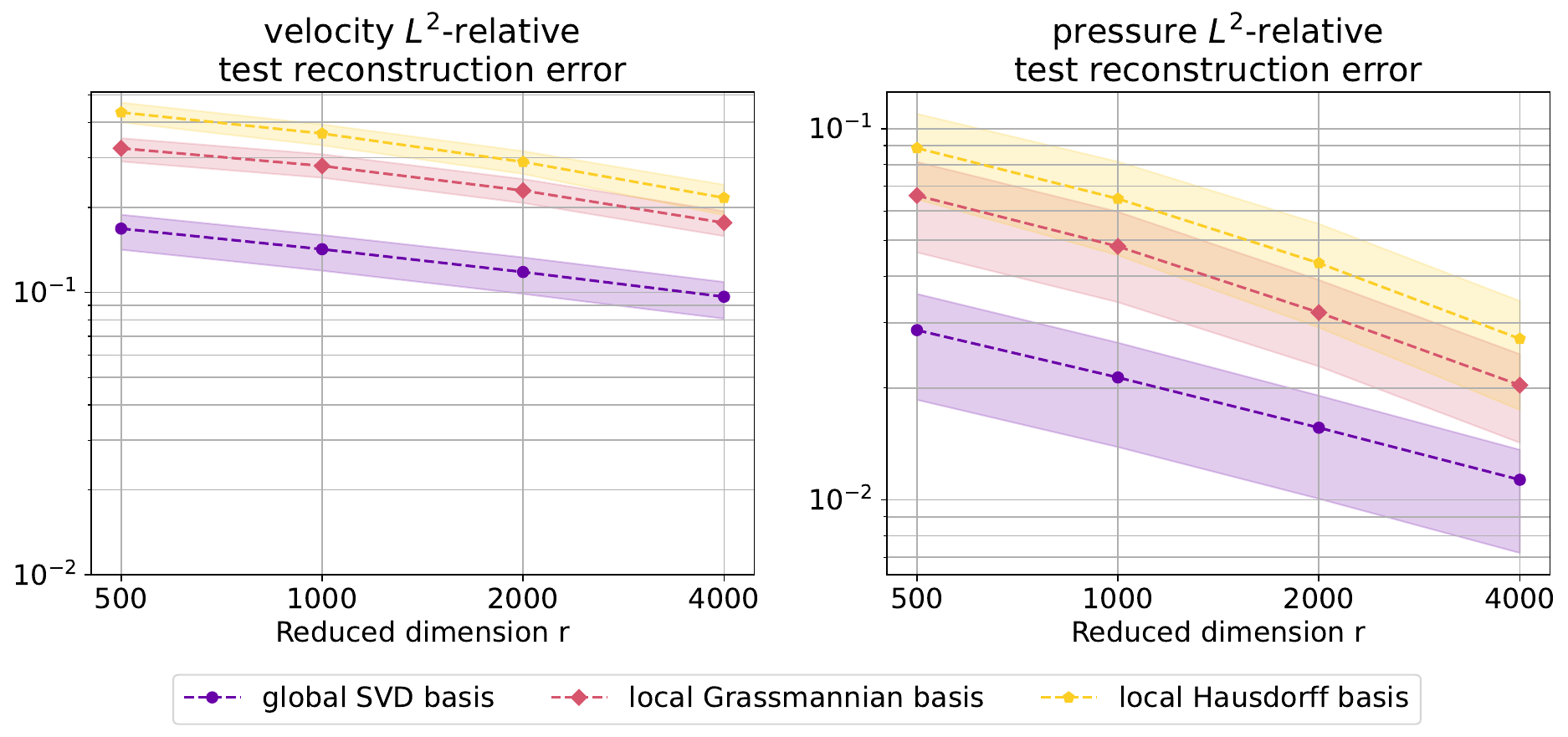}
  \caption{Relative $L^2$-reconstruction errors \eqref{eq:rec} for velocity and pressure averaged over the computational domain and over the considered
  time interval $[0.05s,0.25]$ varying the rSVD rank $r$ and considering a global rSVD and local Hausdorff and Grassman rSVD bases.
	The $25\%$ and $75\%$ percentile are also shown.}
  \label{fig:reclocal}
\end{figure}

Handling shape variability in the context of reduced-order modeling and data assimilation has been also recently discussed in~\cite{galarce2022state} in the context of parametric domains, 
proposing to employ  a \textit{local} rSVD basis, i.e., first clustering the different shapes using a multidimensional scaling (MDS) clustering algorithm, and then assembling the reduced-order model only considering the closest instances.
%
To test an analogous approach, we clustered the training geometries based on the MDS using the dissimilarity matrices $K^{\mathbf u}_{\text{Gr}}$ and $K^{\mathbf u}_{d_{H}}$. Then, given a test geometry
and a fixed value or $r$, we collect the snapshots of the closest $n_{\text{local}}:=\left\lceil r/n_T \right \rceil$ training shapes 
\begin{equation}
  \label{eq:trainsnap}
X^{\mathbf u}_{\text{train}, H}\in\mathbb{R}^{\dofu \times(n_{\text{local}}\,n_T)},\quad X^p_{\text{train}, H}\in\mathbb{R}^{\dofp\times(n_{\text{local}}\,n_T)} ,
\end{equation}
and
\begin{equation}
X^{\mathbf u}_{\text{train}, {\rm Gr}}\in\mathbb{R}^{\dofu \times(n_{\text{local}}\,n_T)},\quad X^p_{\text{train}, {\rm Gr}}\in\mathbb{R}^{\dofp\times(n_{\text{local}}\,n_T)} ,
\end{equation}
depending on the considered dissimilarity metric for the clustering, and compute the corresponding local rSVD bases.
Figure~\ref{fig:reclocal} shows the $L^2$-reconstruction errors \eqref{eq:rec} for the global rSVD basis and for the local ones
In the considered range of dimensions, the global rSVD achieves a better accuracy than the local approaches, using either the Hausdorff or the Grassmann metric. 
As observed above, this might reflect the high geometrical variability of the considered dataset, for which a larger amount of local geometries are required. 

\begin{figure}[!htp]
  \centering
  \includegraphics[width=1\textwidth]{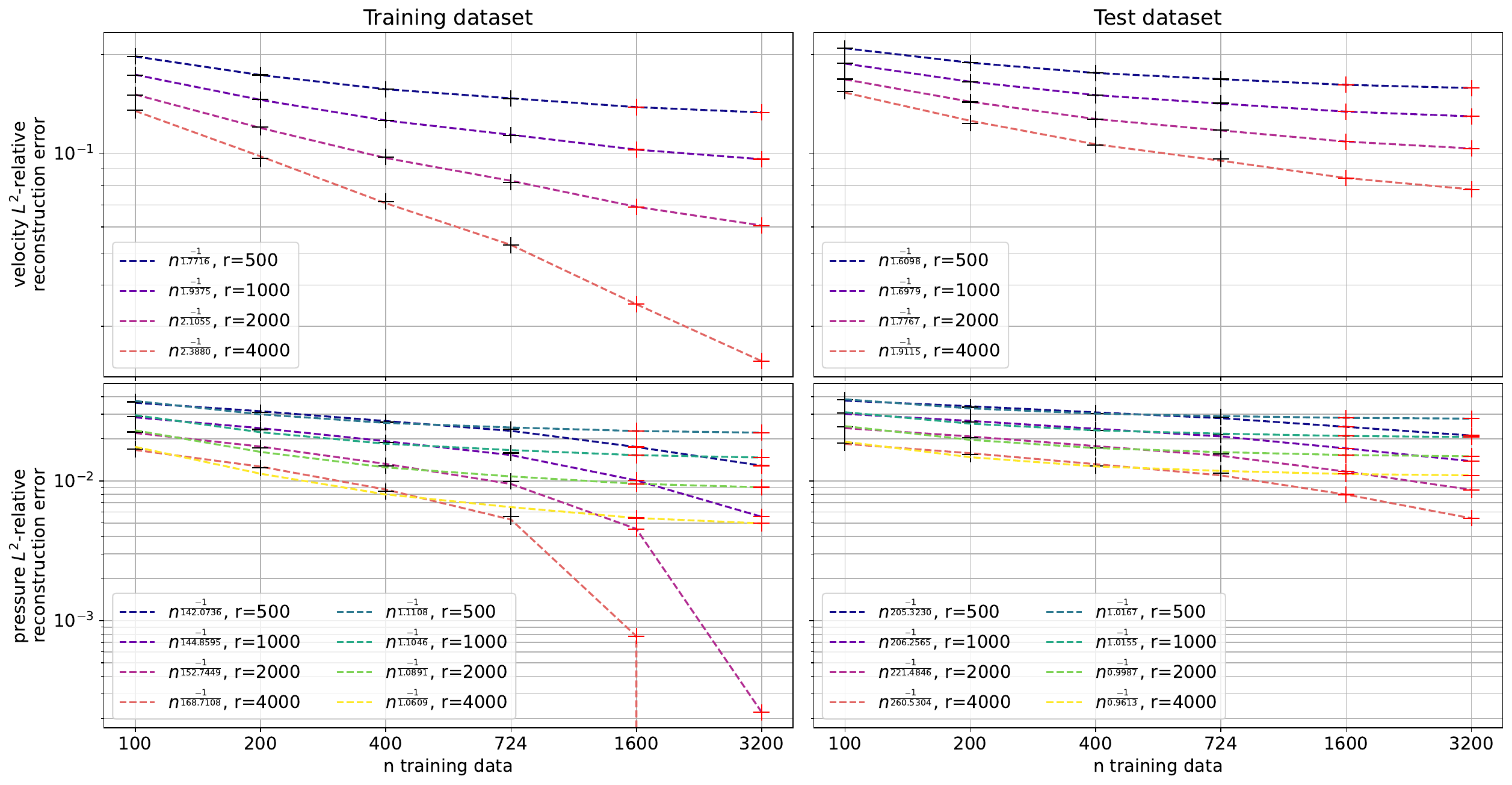}
  \caption{Study of the asymptotic behaviour of the relative $L^2$ reconstruction errors \eqref{eq:rec} increasing the size of the training dataset.
  A nonlinear fit with the function $n_{\text{train}}^{-a}+b$ and multiple initial values for the parameters $a>0,b$ is also shown: red crosses represent the extrapolated values. 
}
  \label{fig:recerr_ntrain}
\end{figure}

Figure~\ref{fig:recerr_ntrain} shows the dependency of the reconstruction error on the number of training data employed. We observe that the global rSVD basis 
is not an efficient approximant of the solution manifold, yielding a decay of the velocity reconstruciton error of the order of $n_{\text{train}}^{-1/2}$. 
Increasing the training dataset with additional geometries is expected to improve the local rSVD error, and local rSVD basis and non-linear dimension reduction methods should be preferred. 

\section{EPD-GNN trained with registered solutions}\label{ssec:pres-gnn}
Building on the shape registration algorithm, we propose a new framework for inference with neural networks on different meshes. 
The dataset is represented by the collection of registered velocity and pressure fields supported on the reference shape.

We employ encode-process-decode graph neural networks (EPD-GNN), introduced in~\cite{pfaff2020learning}, that 
represent the state of the art GNN architectures to perform inference on computational meshes. 
To reduce the computational cost, the reference mesh has been coarsened using TetGen~\cite{Si2015}, reducing the number of vertices from $n_{p, \mathcal{S}}=110676$ to $n_{\text{vertices}}=5181$.
The velocity and pressure fields are transported from the fine to the coarse template meshes and back through RBF interpolation. The metrics (equation~\eqref{eq:l2relerr}) used to validate the results are always evaluated on the fine target meshes. We employ the nearest-neighbours algorithm to enrich each vertex with $e\in\{6, 9, 12\}$ edges to the closest vertices, for a total of $n_{\text{edges}}\in\{36648, 53748, 66135\}$ edges, respectively. Only undirected graphs will be employed.
We consider two inference problems.
\paragraph*{Geometry to velocity (\textit{gnn-gv}) and geometry to pressure (\textit{gnn-gp}) inference}
The input represents a geometrical encoding of the target computational domains with additional velocity b.c.. For each target domain, we evaluate the scalar field that represents the distance from the centerline $Z_{\shape{i}}\in\mathbb{R}^{n_{\text{vertices}}}$, from equation~\eqref{eq:enc_geo}. The pullback of this scalar field to the reference geometry through the registration map together with the pushforwarded coordinates of the vertices of the template geometry $\phi^i_1(\XS)\in\mathbb{R}^{n_{\text{vertices}\times 3}}$ is our $4$-dimensional geometrical encoding. This geometrical encoding is then embedded with a Fourier positional encoding~\cite{sutherland2015error} with $10$ features, through the maps $\{\cos(2^{i}z_j), \sin(2^{i}z_j)\}_{i=0, j=0}^{9, 3}$, where $\mathbf{z}=\{z_i\}_{i=0}^3\in\mathbb{R}^4$ is an arbitrary input vector, for a total of $n_{\text{feat}}=80=10\cdot 2\cdot 4$ geometrical input features. To these it is added the velocity field at $n_{t,\text{GNN}}=8$ times
\begin{equation}\label{eq:int}
  t\in\{0.05s, 0.075s, 0.1s, 0.125s, 0.15s, 0.2s, 0.225s\}=I_t,
\end{equation}
restricted at the boundaries $\Gamma_{\text{in}}$ and $\Gamma_i$ with $i\in\{1, 2, 3, 4\}$ of the target domains and then pulled back to the reference geometry. The value of the velocity boundary field is zero inside the computational domain $\overline{\Omega}\backslash \left(\Gamma_{\text{in}}\cup \left(\cup_{i=1}^{4}\Gamma_i\right)\right)$. The total dimension of the inputs is thus $n_{\text{fnodes}} = 104 =n_{\text{feat}}+n_{t,\text{GNN}}\cdot 3$, where $3$ refers to the number of components of the velocity field. Each edge between the vertices $\mathbf{x}_i$ and $\mathbf{x}_j$ has as features, the vector $\mathbf{x}_i-\mathbf{x}_j$, its $L^2$-norm $\lVert\mathbf{x}\rVert_2$, and the difference between the values of the input vector at the nodes $\mathbf{x}_i$ and $\mathbf{x}_j$, for a total of $n_{\text{fedges}}=4+n_{\text{fnodes}}$ edge features. The output of the GNNs for the problems \textit{gnn-gv} or \textit{gnn-gp} are the velocity field $u\in\mathbb{R}^{n_{\text{vertices}}\times 3n_{t,\text{GNN}}}$ or the pressure field $p\in\mathbb{R}^{n_{\text{vertices}}\times n_{t,\text{GNN}}}$ respectively, evaluated at $n_{t,\text{GNN}}=8$ time instants $t\in I_t$ and supported on the coarse reference mesh.

\paragraph*{Velocity to pressure (\textit{gnn-vp}) inference}
The input is the velocity field $u\in\mathbb{R}^{n_{\text{vertices}}\times 3n_{t,\text{GNN}}}$ at $n_{t,\text{GNN}}=8$ times $t\in I_t$ supported on the coarse reference mesh with $n_{\text{vertices}}=5181$ vertices and $n_{\text{edges}}\in\{36648, 53748, 66135\}$ edges, depending on the number of adjacent nodes $e\in\{6, 9, 12\}$. Each edge between the vertices $\mathbf{x}_i$ and $\mathbf{x}_j$ has, as features, the vector $\mathbf{x}_i-\mathbf{x}_j$, its $L^2$-norm $\lVert\mathbf{x}\rVert_2$, and the difference between the values of the velocity field at $n_{t,\text{GNN}}=8$ time instances at the nodes $\mathbf{x}_i$ and $\mathbf{x}_j$, for a total of $n_{\text{fedges}}=28=4+n_{t,\text{GNN}}\cdot 3$ edge features, $3$ stands for the components of the velocity field. The output is the pressure field $p\in\mathbb{R}^{n_{\text{vertices}}\times n_{t,\text{GNN}}}$ at $n_{t,\text{GNN}}=8$ times $t\in I_t$ supported on the coarse reference mesh. \newline

As our EPD-GNN model we choose the \textit{MeshGraphNet} architecture implemented in NVIDIA-Modulus~\cite{modulus}, based on \texttt{pytorch}~\cite{NEURIPS2019_9015}. The hyperparameters are the width of the network $w$, i.e., the number of consecutive EPD layers, the common hidden dimension $h$ of the node encoder and decoders and the edge encoder, and also the number of edges $e$ of each node. The loss is the relative mean squared error. We apply the ADAM stochastic optimization method~\cite{kingma2014adam} to train the EPD-GNNs with a scheduler used to halve the learning rate when the validation error does not decrease after $200$ epochs. The initial learning rate value is $0.001$. The $52$ test geometries are not employed during the training. We perform the optimization on a single GPU NVIDIA A100-SXM4 with 40GB of graphics RAM size.

We perform two hyperparameter studies.  Firstly, we consider the values $$(w,h)\in\{(10, 64), (15, 128), (20, 256)\}$$ and fix $e=9$, with $n_{\text{train}}=720$ training geometries and $n_{\text{val}}=4$ validation geometries, and $n_{\text{epochs}}=500$. We then select the best model as the one with lowest validation error. 
Secondly, we fix $(w,h)=(30, 256)$ and choose $e\in\{6, 9, 12\}$, with $n_{\text{train}}=680$ and $n_{\text{val}}=44$, and $n_{\text{epochs}}=1000$. We then select the best model as the one at the last epoch.

\begin{figure}[!htp]
  \centering
  \includegraphics[width=.85\textwidth]{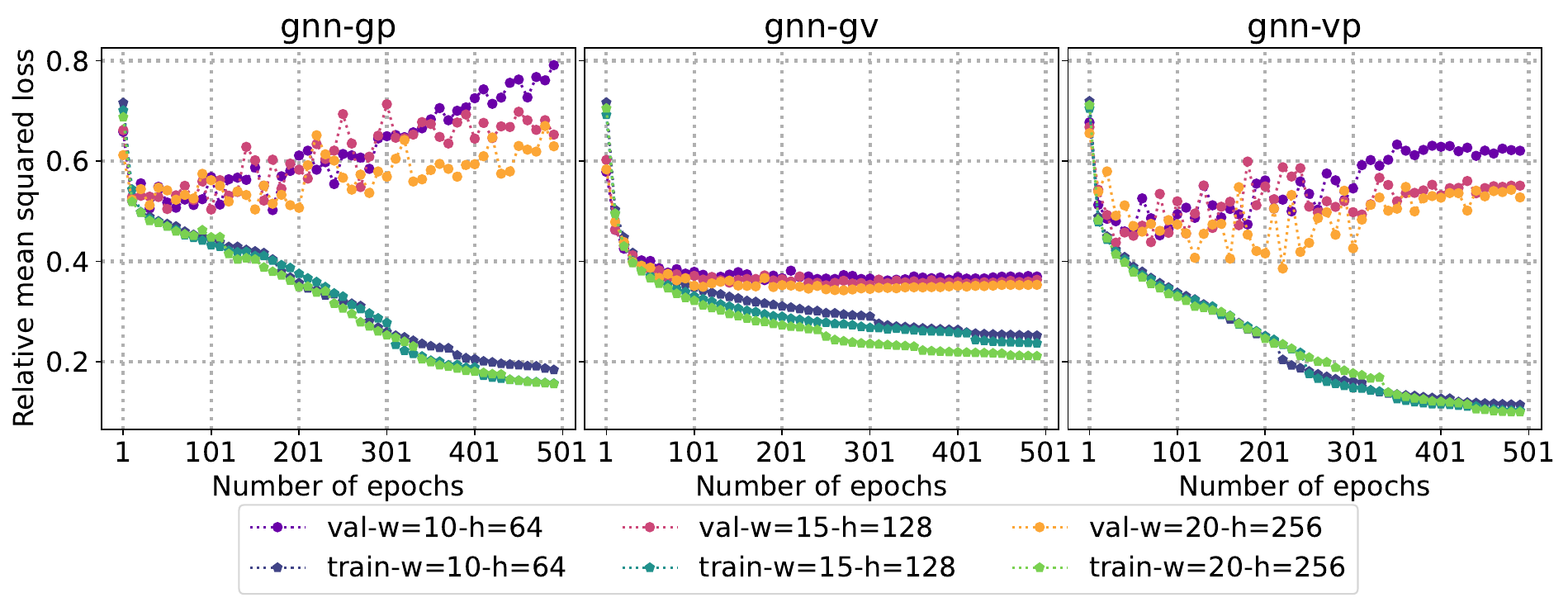}
  \caption{First hyperparameter study: $e=9$, $(w,h)\in\{(10, 64), (15, 128), (20, 256)\}$, $n_{\text{train}}=720$, $n_{\text{val}}=4$, and $n_{\text{epoch}}=500$.}
  \label{fig:overfitting}
\end{figure}

\begin{figure}[!htp]
  \centering
  \includegraphics[width=.85\textwidth]{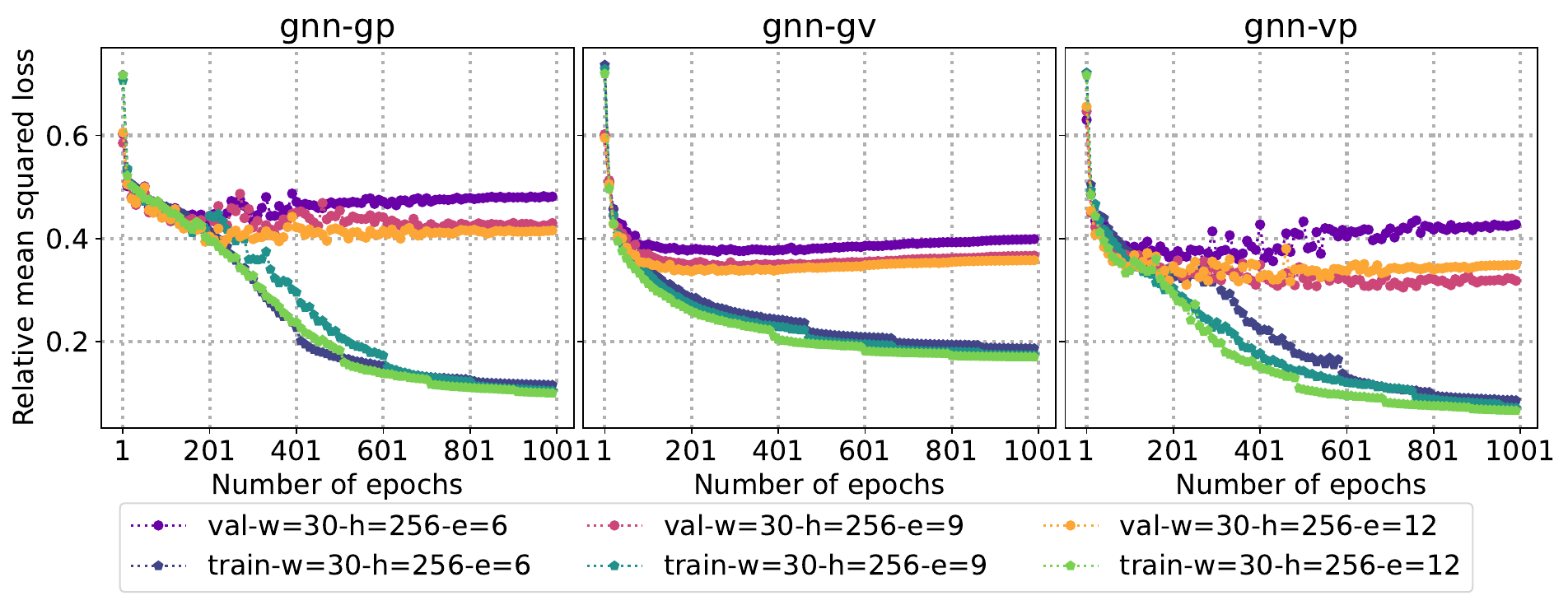}
  \caption{Second hyperparameter study: $e\in\{6, 9, 12\}$, $(w,h)=(30, 256)$, $n_{\text{train}}=680$, $n_{\text{val}}=44$, and $n_{\text{epoch}}=1000$.}
  \label{fig:overfitting_edges}
\end{figure}

The values of the loss during the training for the first and second hyperparameter studies are reported in figures~\ref{fig:overfitting} and~\ref{fig:overfitting_edges}, respectively. A comparison of the training and validation mean squared loss, computed on the coarse mesh with $n_{\text{vertices}}=5181$ vertices, highlights a clear overfitting phenomenon in our limited data regime, with a high generalization error compared to the training error. From the convergence behavior of the training error, we are hopeful that increasing the training dataset would bring better results. The optimal way to increase the training dataset is a future direction of research.

\begin{figure}[!htp]
  \centering
  \includegraphics[width=1\textwidth]{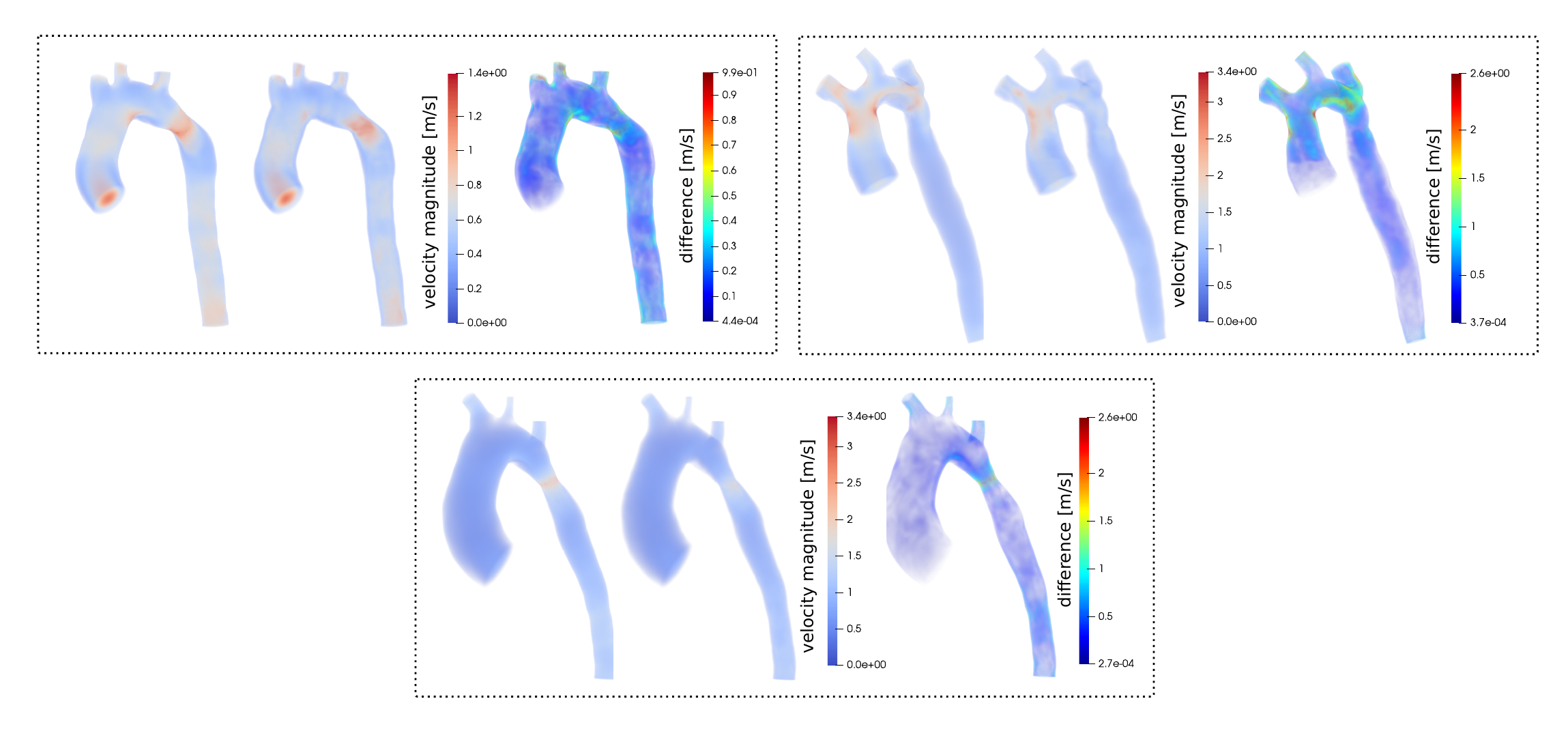}
  \caption{Results on the test dataset of the EPD-GNNs for the problem \textit{gnn-gv} at systolic peak $t=0.125s$: true velocity field magnitude, predicted velocity field magnitude and difference between the two scalar fields. \textbf{Left: }Minimum test $L^2$-relative error $\boldsymbol{\epsilon}_{\widehat{\mathbf u}}=\textbf{0.26}$. \textbf{Right: } Maximum test $L^2$-relative error $\boldsymbol{\epsilon}_{\widehat{\mathbf u}}=\textbf{0.44}$. \textbf{Bottom: } Median test $L^2$-relative error $\boldsymbol{\epsilon}_{\widehat{\mathbf u}}=\textbf{0.31}$.}
  \label{fig:gv}
\end{figure}

\begin{figure}[!ht]
  \centering
  \includegraphics[width=1\textwidth]{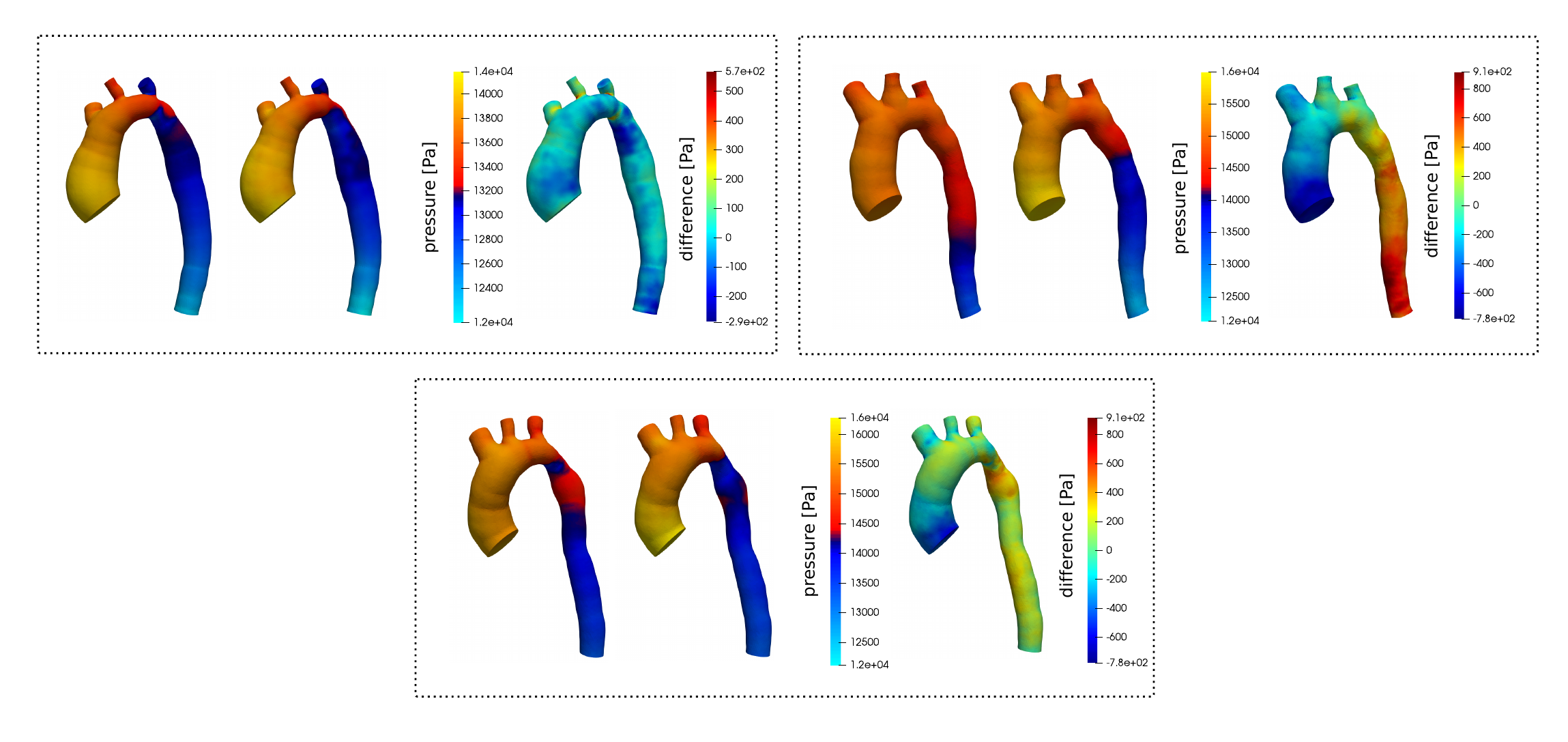}
  \caption{Results on the test dataset of the EPD-GNNs for the problem \textit{gnn-gp} at systolic peak $t=0.125s$: true pressure field, predicted pressure field and difference between the two scalar fields. The values of the pressure fields are rescaled to the same average. \textbf{Left: }Minimum test $L^2$-relative error $\boldsymbol{\epsilon}_{\widehat{p}}=\textbf{0.11}$. \textbf{Right: } Maximum test $L^2$-relative error $\boldsymbol{\epsilon}_{\widehat{p}}=\textbf{0.83}$. \textbf{Bottom: } Median test $L^2$-relative error $\boldsymbol{\epsilon}_{\widehat{p}}=\textbf{0.27}$.}
  \label{fig:gp}
\end{figure}

\begin{figure}[!ht]
  \centering
  \includegraphics[width=1\textwidth]{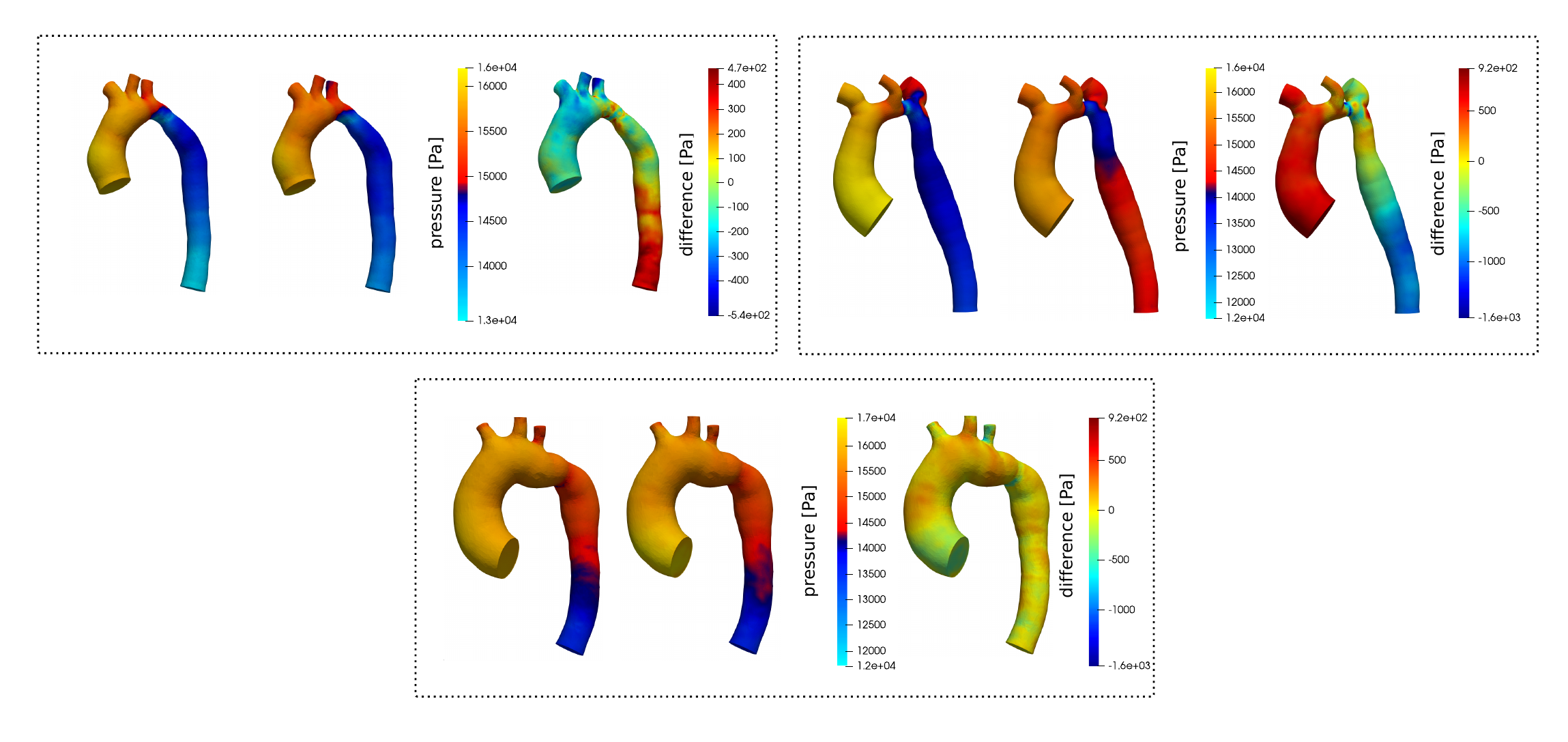}
  \caption{Results on the test dataset of the EPD-GNNs for the problem \textit{gnn-vp} at systolic peak $t=0.125s$: true pressure field, predicted pressure field and difference between the two scalar fields. The values of the pressure fields are rescaled to the same average. \textbf{Left: }Minimum test $L^2$-relative error $\boldsymbol{\epsilon}_{\widehat{\mathbf p}}=\textbf{0.11}$. \textbf{Right: } Maximum test $L^2$-relative error $\boldsymbol{\epsilon}_{\widehat{\mathbf p}}=\textbf{0.54}$. \textbf{Bottom: } Median test $L^2$-relative error $\boldsymbol{\epsilon}_{\widehat{\mathbf p}}=\textbf{0.22}$.}
  \label{fig:vp}
\end{figure}

To evaluate the prediction errors we will consider the $L^2$-relative errors for the velocity $\epsilon_{\widehat{\mathbf{u}}}$ and pressure $\epsilon_{\widehat{p}}$ evaluated on target shapes $\mathcal{T}$:
\begin{equation}
  \label{eq:l2relerr}
  \epsilon_{\widehat{\mathbf{u}}} = \frac{\lVert \widehat{\mathbf{u}}_{\text{true}}-\widehat{\mathbf{u}}\rVert_2}{\lVert \widehat{\mathbf{u}}_{\text{true}}\rVert_2},\qquad\epsilon_{\widehat{p}} = \frac{\lVert \widehat{p}_{\text{true}}-\overline{\widehat{p}}_{\text{true}} -(\widehat{p}-\overline{\widehat{p}})\rVert_2}{\lVert \widehat{p}_{\text{true}}-\overline{{\widehat{p}}}_{\text{true}} \rVert_2},
\end{equation}
where $\widehat{\mathbf{u}}_{\text{true}}$ and $\widehat{p}_{\text{true}}$ are the high-fidelity velocity and pressure fields obtained from the solution of the Navier--Stokes equation \eqref{eq:3dnse} on the target domain, $\widehat{\mathbf{u}}$ and $\widehat{p}$ are the predicted velocity and pressure fields, and
$\overline{\widehat{p}}_{\text{true}}\in\mathbb{R}$ and $\overline{\widehat{p}}\in\mathbb{R}$ denote the averages of the pressure fields. Hat symbols denote quantities defined on the target geometries. The minimum, maximum and median relative $L^2$-errors for the problems $\textit{gnn-gv}$, $\textit{gnn-gp}$, and $\textit{gnn-vp}$, are reported in Table~\ref{tab:gnns}: we evaluated these errors from the best models selected from the hyperparameter studies. The fields associated to the minimum, maximum and median values are shown in figure~\ref{fig:gv} for \textit{gnn-gv}, figure~\ref{fig:gp} for \textit{gnn-vp}, and figure~\ref{fig:vp} for \textit{gnn-vp}. In the following sections we will compare these results with PBDW (section~\ref{sec:da}) and pressure estimators (section~\ref{sec:prec}), using only the best model selected from the first hyperparameter study.
\begin{table}[H]
  \centering
  \begin{tabular}{l|ccc|ccc|ccc}
 & \multicolumn{3}{c}{$\textit{gnn-gv}~(\epsilon_{\widehat{\mathbf{u}}})$} & \multicolumn{3}{|c}{$\textit{gnn-gv}~(\epsilon_{\widehat{\mathbf{u}}})$} & \multicolumn{3}{|c}{$\textit{gnn-gv}~(\epsilon_{\widehat{\mathbf{u}}})$}\\ \hline
& \textbf{min} & \textbf{max} & \textbf{median} & \textbf{min} & \textbf{max} & \textbf{median} & \textbf{min} & \textbf{max} & \textbf{median}\\
       \hline
      First hp study  & $0.26$ & $0.44$ & $0.31$ & $0.11$ & $0.83$ & $0.27$ & $0.11$ & $0.54$ & $0.22$ \\
      \hline
      Second hp study & $0.25$ & $0.45$ & $0.31$ & $0.11$ & $0.93$ & $0.30$ & $0.11$ & $0.67$ & $0.24$ \\
      \hline
  \end{tabular}\hspace{5mm}
  \caption{Relative $L^2$-errors of the velocity and pressure predictions corresponding the best architectures from the first and second hyperparameter (hp) studies.}
  \label{tab:gnns}
\end{table}

EPD-GNNs could be seen as a potential approach to reduce the effort and the time required for accurate experimental acquisition of 4DMRI data. However, in our case, the results suggest that the training of EPD-GNNs requires more data to achieve a better accuracy. 
Similar problems have been recently addressed considering GNNs or a combination of NNs and SVD as data-driven surrogate models in simpler settings, i.e., 
considering only healthy geometries, neglecting the secondary branches (LBCA, LCCA, LSA), employing a simplified physical model~\cite{pajaziti2023shape}(where $20$ and $57$ SVD modes for pressure and velocity are sufficient in their case to achieve a good reconstruction error with a different registration method), or employing 1D graphs instead of full 3D geometries~\cite{iacovelli2023novel,pegolotti2024learning}.

The overall low level of accuracy of the EPD-GNNs predictions is confirmed in figures~\ref{fig:gv},~\ref{fig:gp}, and~\ref{fig:vp}, showing the minimum, the maximum and the median $L^2$-relative error for the EPD-GNNs predictions.
These results might suggest that the variability of stenotic aortic geometries requires necessarily a larger training dataset
than only the $724$ training geometries used in the present study. 
This causes noticeable overfitting, see, e.g.,  the value of the loss on the validation set in figure~\ref{fig:overfitting}, and a still high training error.

\begin{figure}[!htp]
  \centering
  \includegraphics[width=0.7\textwidth]{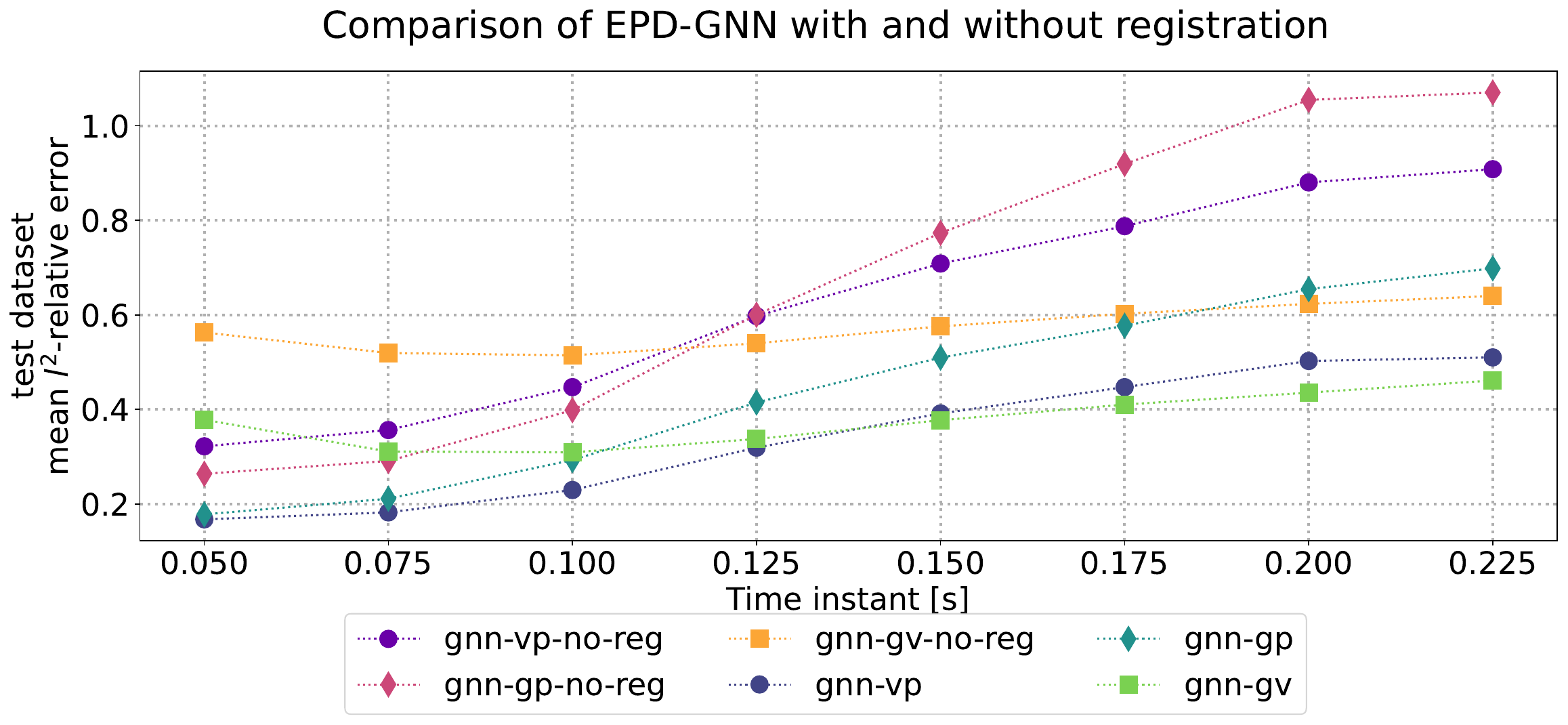}
  \caption{$L^2$-relative errors $\epsilon_{\widehat{\mathbf{u}}}$ and $\epsilon_{\widehat{p}}$ computed on the target computational domains of the predicted velocity (\textit{gnn-gv},\textit{gnn-gv-no-reg}) and pressure fields (\textit{gnn-gp},\textit{gnn-gp-no-reg}, \textit{gnn-vp},\textit{gnn-vp-no-reg}) using EPD-GNN models with and without registration of the datasets.}
  \label{fig:epdgnn_reg_no_reg}
\end{figure}

In figure~\ref{fig:epdgnn_reg_no_reg}, we compare the mean $L^2$-relative error of the new framework over the $52$ test target geometries against the results of an EPD-GNN architecture trained on the original target geometries without registration, and coarsening them with TetGen~\cite{Si2015} as it has been done previously for the template mesh. The definition of the loss, inputs, outputs, and optimization are the same for registered and non-registered datasets. The difference is that the datasets are not supported on the same graph anymore. 
We can observe that the registration represents an efficient encoding of the geometrical and physical features of the problem: for each inference problem \textit{gnn-vp}, \textit{gnn-gp}, and \textit{gnn-gv}, the EPD-GNNs trained on registered data perform better than their alternatives on non-registered data \textit{gnn-gv-no-reg}, \textit{gnn-gp-no-reg}, and \textit{gnn-vp-no-reg}.

\section{Data assimilation of the velocity field}
\label{sec:da}
In subsection~\ref{subsec:sml_rec}, we have shown how to obtain a global rSVD basis $\Phi_{\mathbf u}\in\mathbb{R}^{d_{\mathbf u}\times r_{\mathbf u}}$ for the velocity field on the reference geometry, combining CFD solutions from a database of patient geometries with registration. 
In this section, we address the reconstruction of the velocity field associated to a new patient geometry from a set of velocity observations acquired via 4D flow MRI on a lower resolution.
The data assimilation problem is solved using the Parametrized-Background Data-Weak (PBDW) method, in which, given the observations,
the velocity reconstruction is computed solving a modified least squares problem minimizing the distance of the reconstruction from a physics-informed
linear space and with an additional correction that accounts for the discrepancy with the available measurements.
The method was originally proposed in \cite{MPPY2015} and further analyzed and extended in \cite{cohen2022_nonlinearSpaces,gong2019pbdw}. 

We consider a physics-informed space defined by the global rSVD basis on the template, obtained from different shapes and CFD solutions. Additionally, we extend the approach 
proposed in~\cite{gong2019pbdw} for homogeneous noise to the case of heteroscedastic noise, to handle real applications which require assimilation techniques robust against real data.

\newcommand{\pbdw}[1]{#1_{\rm PBDW}}
\subsection{Parametrized-Background Data-Weak with heteroscedastic noise}
Let us consider a patient geometry $\mathcal T$ and its computational mesh $\Omega_{\mathcal T}$.
We denote with $\phi:[0,1]\times\mathbb{R}^3\rightarrow\mathbb{R}^3$
the map to register the patient $\mathcal T$ on the reference shape. 
Given the velocity rSVD basis $\Phi_{\mathbf u}\in\mathbb{R}^{d_{\mathbf u}\times r_{\mathbf u}}$ on the reference shape, we use the registration map $\phi_1$ and RBF interpolation 
onto the finite element space on $\Omega_{\mathcal T}$  to compute the transported basis $\widehat{\Phi}_{\mathbf u}\in\mathbb{R}^{\widehat{d}_{\mathbf u}\times r_{\mathbf u}}$  on the new patient shape using the pushforward operator \eqref{eq:pushforward}. In particular, $\widehat{d}_{\mathbf u}$ denotes the corresponding number of velocity degrees of freedom of the velocity in the
computational domain $\Omega_{\mathcal T}$.

We assume to have available a set of velocity observations gathered from medical imaging, modelled as linear operators. Given a grid of voxels $\{Q_i\}_{i=1}^{M_{\text{voxels}}}$ s.t. $Q_i=\times_{i=1}^3 [a_i, b_i]$, $b_i>a_i,\ a_i,b_i\in\mathbb{R}_+$, with centers $\mathbf{c}^{\text{vox}}_i\in\mathbb{R}^3$ and vertices $\{\mathbf{x}_i^{\text{vox}}\}_{i=1}^8\subset\mathbb{R}^3$:
\begin{equation}
  \label{eq:voxel}
  l_i(\widehat{\mathbf{v}}) = \frac{1}{9}\left(\sum_{i=1}^{8}\widehat{\mathbf{v}}(\x^{\text{vox}}_i)+\widehat{\mathbf{v}}(\mathbf{c}^{\text{vox}}_i)\right)\approx \int_{Q_i}\widehat{\mathbf{v}}(\x)\ d\x,\qquad l_i:\mathbb{R}^{\widehat{d}_{\mathbf u}}\rightarrow\mathbb{R},
\end{equation}
with $\widehat{\mathbf{v}}:\Omega_{\mathcal{T}}\rightarrow\mathbb{R}^3$. Moreover, we introduce the divergence operator
\begin{equation}
  \label{eq:div}
  l_{\text{div}}(\widehat{\mathbf{v}})=\int_{\Omega_h}\text{div}(\widehat{\mathbf{v}}(\x))\ d\x,\qquad l_{\text{div}}:\mathbb{R}^{\widehat{d}_{\mathbf u}}\rightarrow\mathbb{R},
\end{equation}

\begin{rmk}
  The divergence operator \eqref{eq:div} can be evaluated exactly for velocity fields belonging to piecewise linear finite element space.
\end{rmk}

Let $M := M_{\text{voxels}}+1$ and let us denote with $\mathcal{Z}_{\mathbf u}\in\mathbb{R}^{\widehat{d}_{\mathbf u}\times M}$ the matrix whose columns are the Riesz representers of 
the operators $\{l_i\}_{i=1}^{M_{\text{voxels}}}$ and $ l_{\text{div}}$ with respect to the discrete $\ell^2$ norm.
Moreover, we will denote with 
$y \in \mathbb R^M$ the available set of measurements, and  with $\widehat{\mathbf u}^{\text{true}}\in\mathbb{R}^{\widehat{d}_{\mathbf u}}$ the \textit{true} velocity field, i.e., the unknown field
from which the available measurements are taken, i.e., such that
\begin{equation*}\label{eq:noisy_y}
y = \mathcal{Z}_{\mathbf u} \widehat{\mathbf u}^{\text{true}} + \epsilon_{y},\;\epsilon_y\sim\mathcal{N}(0, S),
\end{equation*}
where $S\in\mathbb{R}^{M\times M}$ is the measurements noise covariance matrix.
 
\begin{rmk}
The operator \eqref{eq:div} will be used to impose incompressibility of the reconstructed velocity field as a fictitious measurement.
Alternatively, the divergence constraint can be imposed exactly with the Piola transform (see, e.g.~\cite{guibert2014group}) acting on the registration map $\phi$.
\end{rmk}

Let us now define the matrices
\begin{align*}
  L=\mathcal{Z}_{\mathbf u}^T\widehat{\Phi}_{u},\quad L\in\mathbb{R}^{ M\times r_{\mathbf u}},\qquad K=\mathcal{Z}_{\mathbf u}^T\mathcal{Z}_{\mathbf u},\quad K\in\mathbb{R}^{ M\times  M}\,.
\end{align*}

In the case of homogeneous noise, the PBDW approach proposed in~\cite{gong2019pbdw} seeks the reconstruction 
in the form of $\pbdw{\widehat{\mathbf u}} = \widehat{\Phi}_{u} \pbdw{z} +  \mathcal{Z}_{\mathbf u} \pbdw{\eta}$ solving
  \begin{equation*}
    \label{eq:pbdw_homo}
    (z_{\text{PBDW}}, \eta_{\text{PBDW}}) = \argmin_{(z, \eta)\in\mathbb{R}^{r_{\mathbf u}}\times \mathbb{R}^M}\xi^{-1} \lVert\eta\rVert^2+\lVert Lz + K\eta-y\rVert^{2}_{S^{-1}},
  \end{equation*}
where $S\in\mathbb{R}^{M\times M}$ is the measurements diagonal covariance matrix and $\xi >0$ is set from a validation dataset and proportional to 
$\lVert y-l(\widehat{\mathbf u}^{\text{true}})\rVert_2/\lVert\mathcal{Z}_{\mathbf u}^{\perp}\widehat{\mathbf u}^{\text{true}}\rVert_2$. The special choices
$\xi=0$ and $S=\text{Id}$ yield the original PBDW formulation.

\begin{problem}[PBDW with heteroscedastic noise] We consider the following extension: given $y \in \mathbb R^M$, find $\pbdw{\widehat{\mathbf u}} = \widehat{\Phi}_{u} \pbdw{z} +  \mathcal{Z}_{\mathbf u} \pbdw{\eta}$ such that 
  \begin{equation}
    \label{eq:pbdw_hetero}
    (z_{\text{PBDW}}, \eta_{\text{PBDW}}) = \argmin_{(z, \eta)\in\mathbb{R}^{r_{\mathbf u}}\times \mathbb{R}^M}\lVert\eta\rVert^2_{R^{-1}}+\lVert Lz + K\eta-y\rVert^{2}_{S^{-1}},
  \end{equation}
where a matrix $R\in\mathbb{R}^{M\times M}$, instead of a single parameter, needs to be set.
\end{problem}

\begin{theorem}[PBDW reconstruction]
  \label{theo:pbdw}
Let us assume that $R$ is chosen such as $R^{-1}=KS^{-1}$. Then the solution to problem \eqref{eq:pbdw_hetero} can be obtained solving the sub-problems:
  \begin{eqnarray}
    z_{\text{PBDW}} & = \argmin_{z\in\mathbb{R}^{r_{\mathbf u}}} \lVert Lz-y\rVert^2_{S^{-1}W^{-1}},\label{eq:pbdw_sub1}\\
    \eta_{\text{PBDW}} & = \argmin_{\eta\in \mathbb{R}^M}\lVert\eta\rVert^2_{R^{-1}}+\lVert K\eta-y_{\text{err}}\rVert^2_{S^{-1}}, \label{eq:pbdw_sub2}
  \end{eqnarray}
  where $W = (K + \text{Id})$ and $y_{\text{err}}=y-Lz_{\text{PBDW}}$.
\end{theorem}
\begin{proof}
  The proof is reported in the appendix~\ref{appendix:pbdw}.
\end{proof}
\begin{rmk}
 Our choice for $R^{-1}$ results in the choice of the prior distribution for $\eta$. According to the Gauss-Markov theorem, the solution of \eqref{eq:pbdw_sub1} is 
 $z_{\text{PBDW}}\sim\mathcal{N}(m_{z_{\text{PBDW}}}, \Sigma_{z_{\text{PBDW}}})$ with
 \begin{equation}
       \label{eq:1pbdw}
       m_{z_{\text{PBDW}}}=\underbrace{(L^T S^{-1}W^{-1} L)^{-1} L^T S^{-1}W^{-1}}_{:=H_{z_{\text{PBDW}}}} y,\;
       \Sigma_{z_{\text{PBDW}}}=(L^T S^{-1}W^{-1} L)^{-1}.
   \end{equation} 
       This results can be interpreted as follows. Let us assume that there exists a reconstruction $z_{\text{true}}$ on the reduced-order space that fits the measurements, i.e.,  
       $y\approx Lz_{\text{true}}+\epsilon_{z}$, up to a noise $\epsilon_{z}\sim\mathcal{N}(0, WS)$. Then, the estimate is unbiased, i.e.,
        $\mathbb{E}[z_{\text{PBDW}}]=z_{\text{true}}$, and it minimizes $\mathbb{E}[\lVert z-z_{\text{true}}\rVert^2_2]$ as well as the covariance $\mathbb{E}[(z-z_{\text{true}})\otimes(z-z_{\text{true}})]$.

 Let $y_{\text{err}}=y-Lz_{\text{PBDW}}=(I-LH_{z_{\text{PBDW}}})y$. Then, the solution to the problem \eqref{eq:pbdw_sub2} 
 can be interpreted as an inverse problem in the Bayesian framework with resulting posterior distribution:
 $$
 \eta_{\text{PBDW}}|y_{\text{err}}\sim\mathcal{N}(m_{\eta_{\text{PBDW}}}, \Sigma_{\eta_{\text{PBDW}}}),
 $$ 
 where
 \begin{equation}\label{eq:2pbdw}
 \begin{aligned}
 m_{\eta_{\text{PBDW}}} & = \underbrace{(KS^{-1}K+KS^{-1})^{-1}KS^{-1}}_{:= H_{\eta_{\text{PBDW}}}}y_{\text{err}},\\
 \Sigma_{\eta_{\text{PBDW}}} & =(KS^{-1}K+KS^{-1})^{-1} = \left[R^{-1}\left(K + \text{Id}\right)\right]^{-1} = W^{-1} R.
 \end{aligned}
 \end{equation}
 Equation \eqref{eq:2pbdw} is equivalent to assuming that there exists a correction on the measurement space 
 $\eta_{\text{true}}$, such that $y_{\text{err}}\approx K\eta_{\text{true}}+\epsilon_{\eta}$, with $\epsilon_{\eta}\sim\mathcal{N}(0, S)$.
 
 Using \eqref{eq:1pbdw} and \eqref{eq:2pbdw}, one obtains that the solution to \eqref{eq:pbdw_hetero} is Gaussian distributed,
 $\pbdw{\widehat{\mathbf u}} = \widehat{\Phi}_{\mathbf u}z_{\text{PBDW}}+\mathcal{Z}_{\mathbf u}\eta_{\text{PBDW}}$,
 \[
 \pbdw{\widehat{\mathbf u}} \sim\mathcal{N}(m_{\pbdw{\widehat{\mathbf u}}}, \Sigma_{\pbdw{\widehat{\mathbf u}}})
 \]
 with 
   \begin{equation}\label{eq:cov}
   \begin{aligned}
       &m_{\widehat{\mathbf{u}}_{\text{PBDW}}} = \widehat{\Phi}_{\mathbf u}  m_{z_{\text{PBDW}}} + \mathcal{Z}_{\mathbf u} m_{\eta_{\text{PBDW}}} = 
 [\widehat{\Phi}_{\mathbf u}  H_{z_{\text{PBDW}}}+\mathcal{Z}_{\mathbf u} H_{\eta_{\text{PBDW}}}-\mathcal{Z}_{\mathbf u} H_{\eta_{\text{PBDW}}}LH_{z_{\text{PBDW}}}]y = H_{\widehat{\mathbf{u}}_{\text{PBDW}}}y  \\
     &\Sigma_{\widehat{\mathbf{u}}_{\text{PBDW}}} = H_{\widehat{\mathbf{u}}_{\text{PBDW}}}S H_{\widehat{\mathbf{u}}_{\text{PBDW}}}^T.
   \end{aligned}
   \end{equation}
\end{rmk}

The next result builds on the estimate of \cite{gong2019pbdw} and, taking into account additional sources of error coming from the registration step, 
provides an error estimate for the PBDW reconstruction.
\begin{theorem}[Error estimate for PBDW reconstruction with heterogeneous noise]
  \label{theo:pbdwmsq}
Let $P_X:\mathbb{R}^{\widehat{d}_{\mathbf u}}\rightarrow\mathbb{R}^{\widehat{d}_{\mathbf u}}$ denote the linear projections onto a linear subspace $X\subset\mathbb{R}^{\widehat{d}_{\mathbf u}}$. Let $H_{\widehat{\mathbf{u}}_{\text{PBDW}}}:\mathbb{R}^{\widehat{d}_{\mathbf u}}\rightarrow\mathbb{R}^{\widehat{d}_{\mathbf u}}$ be the matrix defined in \eqref{eq:cov}, and
$H_l := H_{\widehat{\mathbf{u}}^{\text{PBDW}}} \, \mathcal Z_{\mathbf u}^T$. 
The following estimate holds: 
  \begin{linenomath}\begin{align*}
    \mathbb{E}[\lVert &\widehat{\mathbf{u}}^{\text{true}}-\widehat{\mathbf{u}}_{\text{PBDW}}\rVert_2]\leq &\\
    &\leq\text{tr}(H_{\widehat{\mathbf{u}}_{\text{PBDW}}}SH_{\widehat{\mathbf{u}}_{\text{PBDW}}}^T)^{\tfrac{1}{2}} &\text{(noise error)}\\
    &+\lVert(\text{Id}-H_l)\circ P_{\text{Im}(H_l)}\rVert_2\lVert \widehat{\mathbf{u}}^{\text{true}}\rVert_2&\text{(PBDW bias)}\\
    &+\lVert\text{Id}-H_l\rVert_2 \ \cdot \inf_{{\mathbf u}^{\text{best}}\in\text{col}(X_{\text{train}}^{\mathbf{u}})}\big(&\text{(PBDW stability constant)}\\
    &\lVert(\phi_{\text{RBF}})^{\#}({\mathbf u}^{\text{best}})-\widehat{\mathbf{u}}^{\text{true}}-P_{{\text{Im}(H_l)}}((\phi_{\text{RBF}})^{\#}({\mathbf u}^{\text{best}})-\widehat{\mathbf{u}}^{\text{true}})\rVert_2&\text{(manifold approximation error)}\\
    &+ C\cdot\lVert {\mathbf u}^{\text{best}}-P_{\text{Im}(\Phi_{\mathbf u})} {\mathbf u}^{\text{best}}\rVert_2&\text{(template rSVD approximation error)}\\
    &+\lVert P_{\text{Im}(H_l)}\left((\phi_{\text{RBF}})^{\#}({\mathbf u}^{\text{best}})\right)-(\phi_{\text{RBF}})^{\#}(P_{\text{Im}(\Phi_{\mathbf u})}{\mathbf u}^{\text{best}})\rVert_2\big)&\text{(registration degradation error)}
  \end{align*}\end{linenomath}
  where the matrix $X_{\text{train}}^{u}\in\mathbb{R}^{d_{\mathbf u}\times (n_{\text{train}}n_T)}$ contains, columnwise, the set of training snapshots registered on the reference shape.
  
\end{theorem}
\begin{proof}
  The proof is reported in appendix~\ref{appendix:pbdw} along with an interpretation of the various sources of error.
\end{proof}

\newcommand{\snrho}{\text{SNR-ho}}
\newcommand{\snrhe}{\text{SNR-he}}

\subsection{Heteroscedastic noise model}
In the case of 4DMRI images, the observations often present velocity gradients whose accuracy degrades close to the vessel boundaries (see, e.g.~\cite{zingaro2024advancing, IRARRAZAVAL2019250}).
To account for this aspect, we consider a heteroscedastic noise model depending on three parameters:
the signal-to-noise homoscedastic ration ($\snrho$), the signal-to-noise heteroscedastic ratio ($\snrhe$), and the divergence observations operator's variance ($\sigma_{\text{div}}^2$).
We then subdivide the domain $\Omega_{\mathcal T}$ into a boundary layer $\Omega_{T, \text{he}}$, where measurements are affected by non-homogeneous variance noise, and 
an inner domain $\Omega_{\mathcal T, \text{ho}}$, whose measurements are characterized by standard homogeneous noise.

The covariance matrix $S\in\mathbb{R}^{(3M_{\text{voxels}}+1)\times (3M_{\text{voxels}}+1)}$ in equation~\eqref{eq:pbdw_hetero} is defined as a block matrix
\begin{equation*}
  S = \begin{pmatrix}
    S_{\text{obs}} & 0\\
    0 & \sigma_{\text{div}}^2
    \end{pmatrix},\end{equation*}
where the first block $S_{\text{obs}}\in\mathbb{R}^{3M_{\text{voxels}}\times 3M_{\text{voxels}}}$ is associated to observation operators $\{\{l_i\}_{i=1}^{M_{\text{voxels}}}\}$, and the last diagonal 
entry is associated to the divergence operator ($l_{\text{div}}$).
Let $M_{\text{ho}}$ and $M_{\text{he}}$ denote the amount of voxels whose centers belong to $\Omega_{\mathcal T, \text{ho}}$ and $\Omega_{\mathcal T, \text{he}}$, respectively.
Moreover, let us denote with $\mathbf{c}^{\text{vox}}_i$ the center of voxel $i$.
In our approach, the heteroscedastic noise is modeled as a spatially correlated multiplicative scalar to each velocity observation in the boundary layer, which
only affects the magnitude of the observed velocity vectors.

The covariance matrix $S_{\text{obs}}$ is then split into a homoscedastic ($S_{\text{ho}}$) and a heteroscedastic ($S_{\text{he}}$) block, 
associated to the observation operators in each subdomain
\begin{align*}
  S_{\text{obs}} = \begin{pmatrix}
    S_{\text{ho}} & 0\\
    0 & S_{\text{he}}
    \end{pmatrix},\quad
  S_{\text{ho}} := \left(\frac{\bar{\widehat{\mathbf{u}}}}{\text{SNR-ho}}\right)^2\text{Id}_{3M_{\text{ho}}},\quad
  S_{\text{he}} :=\left(\frac{\bar{\widehat{\mathbf{u}}}}{\text{SNR-he}}\right)^2PCP^T,\;
\end{align*}
where $\bar{\widehat{\mathbf{u}}}=\frac{1}{M_{\text{voxels}}}\sum_{i=1}^{M_{\text{voxels}}}l_i(\widehat{\mathbf{u}})$, $P\in\mathbb{R}^{3M_{\text{he}}\times M_{\text{he}}}$ is the operator 
projecting the velocity vector into its norm and $C\in\mathbb{R}^{M_{\text{he}}\times M_{\text{he}}}$ is the Gramian matrix of the radial basis function kernel
\begin{equation*}
  k(\mathbf{c}^{\text{vox}}_i, \mathbf{c}^{\text{vox}}_j) = \exp\left(-\frac{\lVert\mathbf{c}^{\text{vox}}_i-\mathbf{c}^{\text{vox}}_j\rVert^2_2}{2l_{\mathcal T}}\right)+\epsilon^2\delta_{ij},\qquad
  \forall (i,j) \mid \mathbf{c}^{\text{vox}}_i, \mathbf{c}^{\text{vox}}_j \in \Omega_{T, \text{he}},
\end{equation*}
with length scale $l_{\mathcal T}>0$ and additive homogeneous noise variance $\epsilon^2>0$. 

These parameters depend on the measurement procedure employed. In what follows, we set 
$l_{\mathcal T}=\text{diam}(\Omega_{\mathcal T})/12$ and $\epsilon^2=0.1$.

\subsection{Numerical results and comparison with GNNs}
\begin{figure}[!htp]
  \centering
  \includegraphics[width=0.8\textwidth]{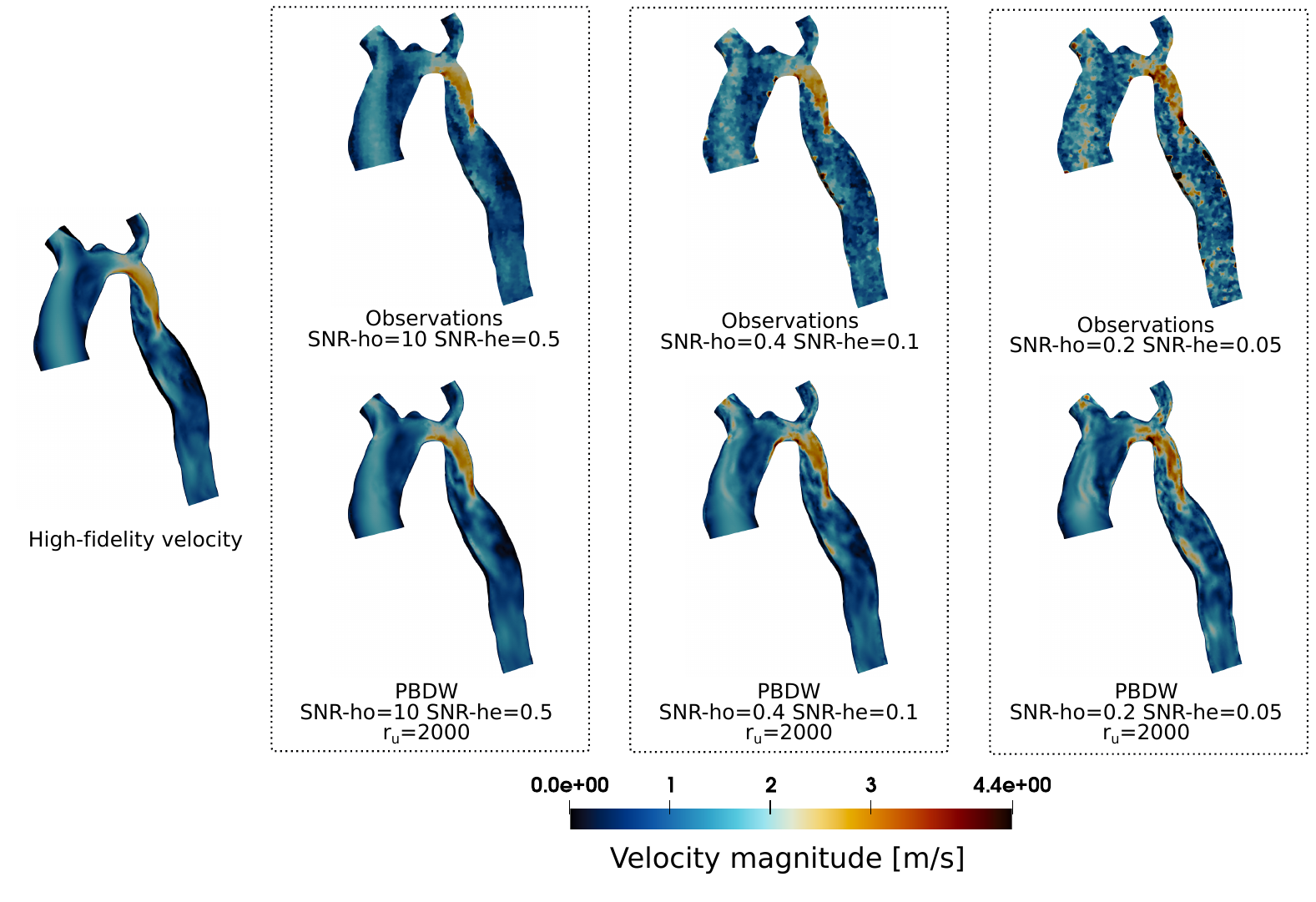}
  \caption{Data assimilation of the velocity field at time $0.1$s from noisy velocity measurements for the test geometry $n=12$ (see figure~\ref{fig:cluster_v}). 
\textbf{Left.} High-fidelity velocity field from the CFD simulations on the corresponding domain. 
\textbf{Right.} Observations with the three considered noise intensities (top) and PBDW predictions with fixed number of template velocity modes $r_{\mathbf u} = 2000$ (bottom).}
  \label{fig:noise_pbdw}
\end{figure}

We consider three levels of noise: $(\text{SNR-ho}, \text{SNR-he})\in \{(10, 0.5), (0.4, 0.1), (0.2, 0.05)\}$ and voxels of resolution $\SI{2e-3}{\meter^3}$.
The velocity observations are computed approximating the integral in equation~\eqref{eq:voxel} with the average of the values of the finite element function on the voxel centers and on the vertices.
Figure \ref{fig:noise_pbdw} shows the resulting observations for the test geometry $n=12$
(the closest to the training velocity solution manifold with respect to the Grassmann metric on the velocity fields, see figure~\ref{fig:cluster_v}),
at a selected time instant for the different noise intensities, together with the corresponding PBDW reconstruction.

\begin{figure}[!ht]
  \centering
  \includegraphics[width=0.9\textwidth]{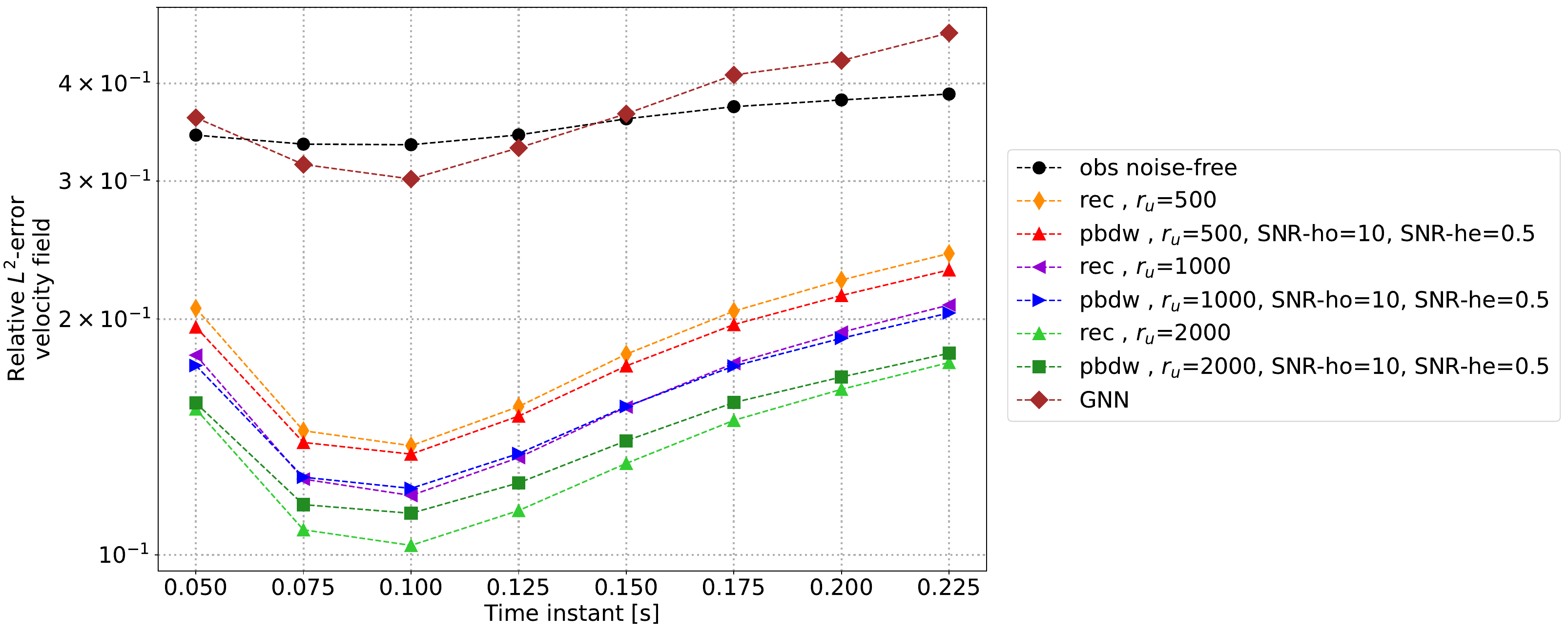}\\
  \includegraphics[width=0.9\textwidth]{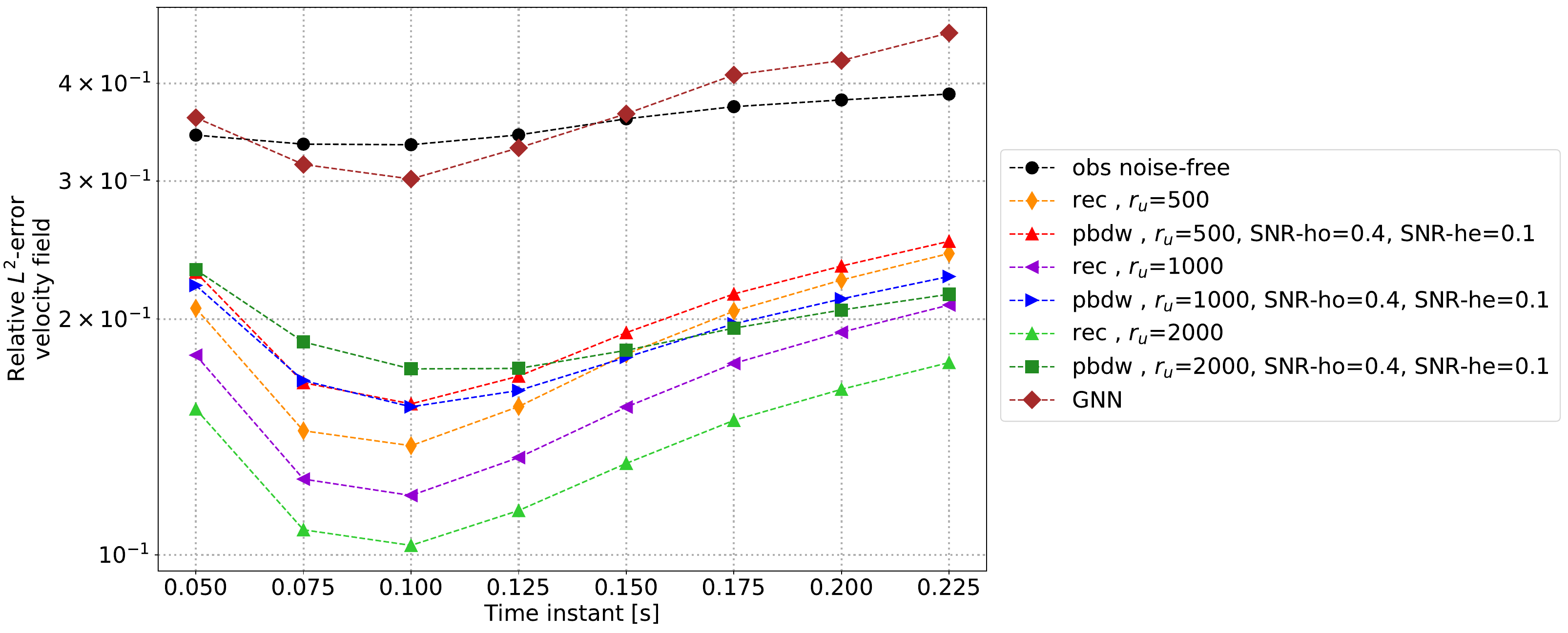}\\
  \includegraphics[width=0.9\textwidth]{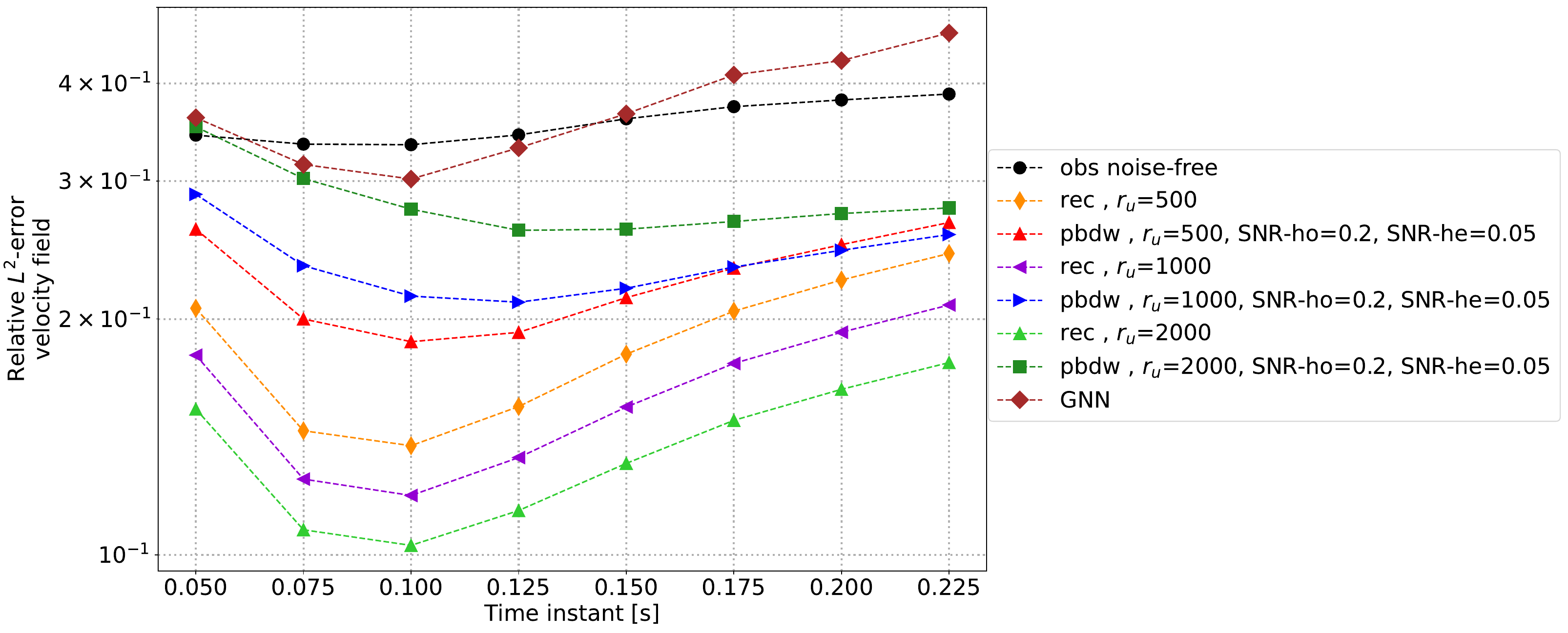}\\
  \caption{Average $L^2$-relative error $\epsilon_{\mathbf u}$ of the velocity fields evaluated on the $52$ target-test geometries: the accuracy of data assimilation with PBDW from velocity observations is compared with respect to direct inference with EPD-GNNs. The rSVD reconstruction error (\textit{rec}) and the error of the noise-free observations (\textit{obs noise-free}) are also shown. \textbf{Top: } $(\text{SNR-ho}, \text{SNR-he})=(10, 0.5)$. \textbf{Middle: } $(\text{SNR-ho}, \text{SNR-he})=(0.4, 0.1)$. \textbf{Bottom: } $(\text{SNR-ho}, \text{SNR-he})=(0.2, 0.05)$.}
  \label{fig:pbdw_vs_gnn_v}
\end{figure}

We compare the reconstructed velocity field with those obtained with the EPD-GNN surrogate models from the inference problem \textit{gnn-gv} (geometry $\mapsto$ velocity) introduced in Section~\ref{ssec:pres-gnn}, which delivers a velocity prediction solely from the geometry data.
The results are shown in figure~\ref{fig:pbdw_vs_gnn_v}, depicting the $L^2$ average relative error $\epsilon_{\mathbf u}$ on the test dataset of $52$ geometries, using PBDW and EPD-GNNs.
The errors are evaluated on the target geometry, after transporting the predicted velocity fields with the registration map in the case of EPD-GNNs.

Notice that the 4DMRI data are not used in the inference problem with EPD-GNNs. The comparison with PBDW is shown to underline, in the case of limited data, the difference in accuracy between a purely data-driven inference problem, such as the EPD-GNN,  and a physics-based data assimilation method that incorporates a state space
constructed using geometrical and physical information, as well as an additional set of observations.

Figure~\ref{fig:pbdw_vs_gnn_v} shows that increasing the level of noise affects the stability properties of the rSVD basis used in PBDW. For the lowest noise
$(\text{SNR-ho}, \text{SNR-he})=(10, 0.5)$  the best results are obtained with $r_{\mathbf u} = 2000$, while $r_{\mathbf u} = 500$ is the best performing case for $(\text{SNR-ho}, \text{SNR-he})=(0.2, 0.05)$.
In the ideal, noise-free case, the optimality properties of the PBDW guarantee that the reconstruction error is lower than the sole rSVD approximation error 
for specific $r<M$. For the case with the lowest noise $(\text{SNR-ho}, \text{SNR-he})=(10, 0.5)$, it can be observed that the PBDW reconstruction error
is lower than the rSVD approximation error for $r_{\mathbf u}=500$, but higher for $r_{\mathbf u}=2000$.

\begin{figure}[!htp]
  \centering
  \includegraphics[width=0.85\textwidth]{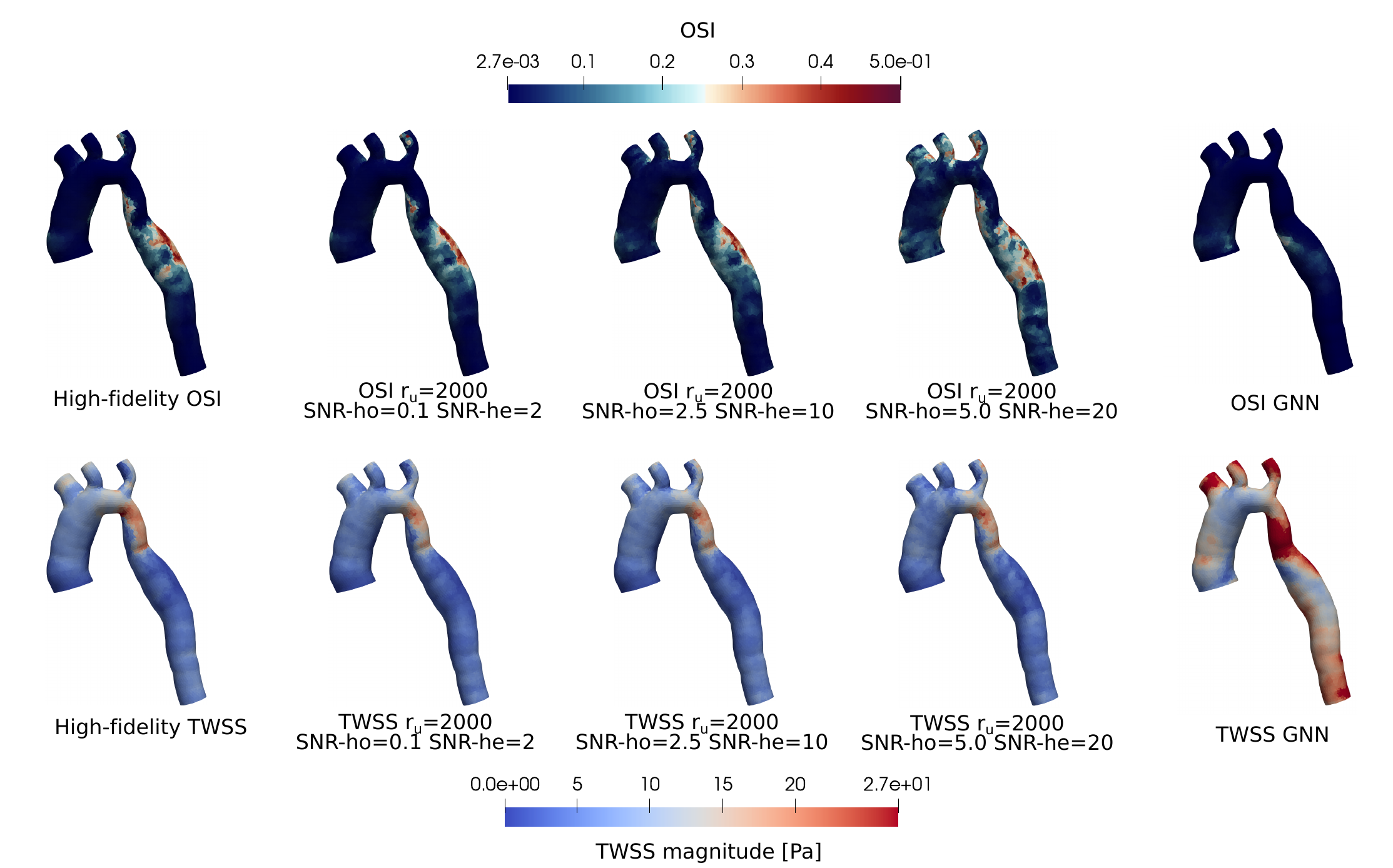}
  \caption{Time-averaged wall shear stress (TWSS, \textbf{top}) and oscillatory index (OSI, \textbf{bottom}) for test geometry $n=12$, figure~\ref{fig:cluster_v}. The results are shown for different noise levels $(\text{SNR-ho}, \text{SNR-he})\in \{(10, 0.5), (0.4, 0.1), (0.2, 0.05)\}$ and fixed rSVD rank $r_{\mathbf u}=2000$.}
  \label{fig:twss_12}
\end{figure}

\begin{figure}[!htp]
  \centering
  \includegraphics[width=0.8\textwidth]{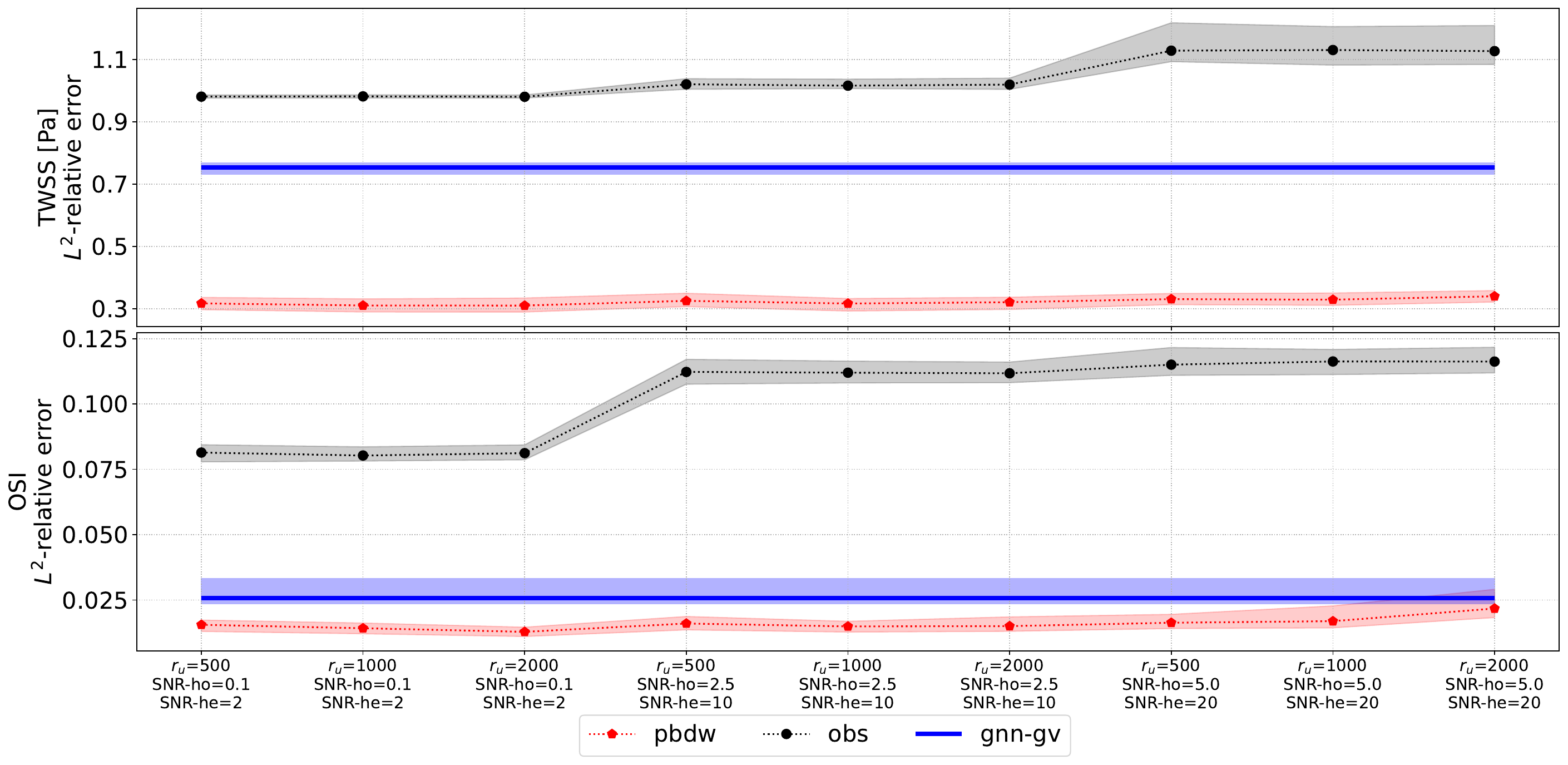}
  \caption{Average $L^2$-relative error of the time averaged wall shear stress (TWSS, top) and of the oscillatory index (OSI, bottom) on the $52$ test geometries,
calculated from the voxel observations (\textit{obs}), the PBDW reconstruction (\textit{pbdw}) and the EPD-GNN reconstruction (\textit{gnn-gv}).
  The results are shown for different noise levels $(\text{SNR-ho}, \text{SNR-he})\in \{(10, 0.5), (0.4, 0.1), (0.2, 0.05)\}$ and rSVD ranks $r_{\mathbf u}\in\{500, 1000, 2000\}$. The $25\%$ and $75\%$ percentile are shown with shaded regions.}
  \label{fig:twss_osi}
\end{figure}

From the high-resolution reconstructed velocity field, obtained through PBDW or EPD-GNNs, clinically relevant biomarkers such as the time-averaged wall shear stress (TWSS)
\begin{align*}
  \tau_{\text{wss}}(\mathbf{x}, t)&=\mu\frac{\partial}{\partial\mathbf{n}(\mathbf{x}, t)}\left(\mathbf{u}(\mathbf{x}, t)-(\mathbf{u}(\mathbf{x}, t)\cdot\mathbf{n}(\mathbf{x}, t))\mathbf{n}(\mathbf{x}, t)\right),\\
  \tau_{\text{twss}}(\mathbf{x})&=\int_{t=0.05s}^{t=0.225s}\tau_{\text{wss}}(\mathbf{x}, t)\,dt\approx\frac{1}{8}\sum_{i=0}^{7}\tau_{\text{wss}}(\mathbf{x}, 0.05+i\cdot\Delta t),
\end{align*}
and the oscillatory shear index OSI$_I$ relative to the time interval $I=[0.05s,0.225s]$
\[
  OSI_{I}(\mathbf{x})=\frac{1}{2}\left(1-\frac{\left.|\int_{t=0.05s}^{t=0.225s}\tau_{\text{wss}}(\mathbf{x}, t)\,dt\right|}{\int_{t=0.05s}^{t=0.225s}|\tau_{\text{wss}}(\mathbf{x}, t)|\,dt}\right)\approx\frac{1}{2}\left(1-\frac{\left.|\frac{1}{8}\sum_{i=0}^{7}\tau_{\text{wss}}(\mathbf{x},  0.05+i\cdot\Delta t)\right|}{\frac{1}{8}\sum_{i=0}^{7}|\tau_{\text{wss}}(\mathbf{x},  0.05+i\cdot\Delta t)|}\right)
\]
can be computed.
The $L^2$-relative error for the test geometry $n=12$ are shown in figure~\ref{fig:twss_12} for different noise intensities,
while the results of the quantitative study on all the $52$ test geometries are presented in figure~\ref{fig:twss_osi}. 
PBDW achieves a satisfactory accuracy in all considered cases, while the prediction with EPD-GNNs fails, confirming that the model is not able to capture the
high geometric variability. As previously noticed, the EPD-GNN only relies on the geometry data. 
The purpose of the comparison with PBWD is to underline the problematics of the EPD-GNN approach in this clinical context, 
in order to address them in future studies with a higher computational budget and a higher amount of data.

\section{Pressure estimation}
\label{sec:prec}
This section focuses on the problem of estimating the pressure field from velocity data.
We benchmark the pressure reconstruction method based on the EPD-GNN (section \ref{ssec:pres-gnn}) against two variational-based approaches, the pressure Poisson estimator (PPE) and the Stokes estimator (STE).
A complete comparison of methodologies is presented in \cite{bertoglio2018relative}. We choose to consider PPE due to its simple implementation and already popularized use, and STE, because it has been benchmarked as the best method in \cite{bertoglio2018relative}. As it has been already discussed, joint reconstructions with PBDW for velocity and pressure, as in \cite{galarce2023displacement}, will not be applied due to the loss of accuracy in the overall reconstruction.

The pressure estimation problem is considered in a given target shape $\mathcal T$. 
We will denote with $\Omega_{\mathcal T}$ the corresponding computational domain, with $\mathcal{T}_h$ its triangulation, 
and with $\widehat{d}_{\mathbf u}$ and $\widehat{d}_p$ the degrees of freedom of the underlying finite element spaces for velocity and pressure, respectively. The diameter of a generic element $K \in \mathcal T_h$ is $h_K$. 
The input velocity field at a time $t_n$ will be denoted by $\widehat{\mathbf u}^n \sim\mathcal{N}(\mathbf m^n_{\widehat{\mathbf u}}, \Sigma^n_{\widehat{\mathbf u}})$ and it is assumed to be Gaussian distributed, 
with $\mathbf m^n_{\widehat{\mathbf u}}\in\mathbb{R}^{\widehat{d}_{\mathbf u}}$ and $\Sigma^n_{\widehat{\mathbf u}}\in\mathbb{R}^{\widehat{d}_{\mathbf u}\times \widehat{d}_{\mathbf u}}$. For example, $\widehat{\mathbf u}^n$ could be obtained from 4D-flow MRI data with heteroscedastic PBDW.

\subsection{Variational-based pressure estimators}
\label{subsec:ppestedef}

In the pressure-Poisson estimator (PPE) \cite{ebbers2001}, the velocity field is directly inserted in the right-hand-side variational form of the incompressible Navier--Stokes equations, and a suitable pressure field is recovered solving the resulting problem.
\begin{problem}[PPE]
  \label{def:ppe}
  Given three consecutive velocity time steps $\widehat{\ub}^{n}$, $\widehat{\ub}^{n+1/2}$, and $\widehat{\ub}^{n+1}$, find the
pressure at the intermediate time $\widehat{p}^{n+1/2}_{\text{PPE}}\in\mathbb{P}^1(\mathcal{T}_h)$ such that
  \begin{align}
    \int_{\Omega_{\mathcal T}} \nabla \widehat{p}^{n+1/2}_{\text{PPE}}\cdot\nabla q  =
     -\frac{\rho}{\tau}\int_{\Omega_{\mathcal T}} (\widehat{\ub}^{n+1}-\widehat{\ub}^{n})\cdot\nabla q - \rho\int_{\Omega_{\mathcal T}} (\widehat{\ub}^{n+1/2}\cdot\nabla\widehat{\ub}^{n+1/2})\cdot\nabla q,\;
     \forall q\in\mathbb{P}^1(\mathcal{T}_h),
     \label{eq:ppe-weak}
  \end{align}
  with boundary conditions $\widehat{p}^{n+1/2}_{\text{PPE}}=q=0$ on $\partial \Omega_{\mathcal T}$. 
  Problem \eqref{eq:ppe-weak} can be equivalently written in matrix form as
  \begin{equation} \label{eq:ppe}
  A_{\text{PPE}}\,\widehat{p}^{n+1/2}_{\text{PPE}} = M_{\text{PPE}}\widehat{\ub}_{n+1} - M_{\text{PPE}}\widehat{\ub}_{n} + Q_{\text{PPE}}(\widehat{\ub}^{n+1/2}, \widehat{\ub}^{n+1/2}),
  \end{equation}
with natural definition of the stiffness matrix $A_{\text{PPE}}$, the mass matrix $M_{\text{PPE}}$,  and the advection term $Q_{\text{PPE}}$. 
 \end{problem} 

A bias correction for the estimator is obtained as the solution to the following problem:
Find $b_{\text{PPE}}\in\mathbb{P}^1(\mathcal{T}_h)$ such that 
  \begin{equation}\label{eq:ppe-bias}
    \int_{\Omega_{\mathcal T}} \nabla b_{\text{PPE}}\cdot\nabla q =-\rho\sum_{T\in\mathcal{T}_h}\sum_{i=1}^{\widehat{d}_{\mathbf u}}\sum_{j=1}^{\widehat{d}_{\mathbf u}}\int_T \Psi_{i,j}\nabla q,
    \; \forall q\in\mathbb{P}^1(\mathcal{T}_h)\,,
  \end{equation}
where we introduced the notaton
\begin{equation*}\label{eq:psi-bias-corr}
\Psi_{i,j} : = \left((\phib_i\cdot\nabla\phib_j)\odot\Sigma^{n+1/2}_{ij}\right), \; i,j=1,\hdots,\widehat{d}_{\mathbf u}
\end{equation*}
and $\Sigma^{n+1/2}_{ij}$ stands for the $3\times3$ $(i,j)$-subblock of the covariance matrix $\Sigma^{n+1/2}_{\widehat{\mathbf u}}$ corresponding to the support points 
of the finite element shape functions $\phib_i$,and $\phib_j$, and $\odot$ stands for the element-wise Hadamard product of two matrices.


In the Stokes estimator (STE) \cite{svihlova_2016}, the velocity field is inserted in the right-hand-side of the Navier--Stokes equations, but, unlike the PPE, the 
variational problem for the pressure field is formulated as a Stokes projection, including an additional corrector. As shown, e.g.,  in \cite{bertoglio2018relative}, this approach allows to 
obtain more robust results, especially against noisy velocity data.
\begin{problem}[STE]
  \label{def:ste}
    Given three consecutive velocity time steps $\widehat{\ub}^{n}$, $\widehat{\ub}^{n+1/2}$, and $\widehat{\ub}^{n+1}$, find $(\wb,\widehat{p}^{n+1/2}_{\text{STE}})\in [\mathbb{P}^1(\mathcal{T}_h)]^d \times \mathbb{P}^1(\mathcal{T}_h)$ such that
  \begin{equation}\label{eq:ste-weak}
  \begin{aligned}
    \int_{\Omega_{\mathcal T}} \nabla\wb:\nabla\mathbf{z} &- \int_{\Omega_{\mathcal T}} \widehat{p}^{n+1/2}_{\text{STE}}(\nabla\cdot\mathbf{z})+\int(\nabla\cdot\wb)q\\
     &+ \sum_{K\in\mathcal{T}_h}C_s h^2_K\int_K \nabla \widehat{p}^{n+1/2}_{\text{STE}}\cdot\nabla q=\\
     &-\frac{\rho}{\tau}\int_{\Omega_{\mathcal T}} (\widehat{\ub}^{n+1}-\widehat{\ub}^{n})\cdot\mathbf{z} - \rho\int_{\Omega_{\mathcal T}} (\widehat{\ub}^{n+1/2}\cdot\nabla\widehat{\widehat{\ub}^{n+1/2}})\cdot\wb-\mu\int_{\Omega_{\mathcal T}} \nabla\widehat{\ub}^{n+1/2}:\nabla\mathbf{z}\\
     &+\sum_{K\in\mathcal{T}_h}C_s h^2_K\rho\int_K (\mu\Delta\widehat{\ub}^{n+1/2}-\widehat{\ub}^{n+1/2}\cdot\nabla\widehat{\ub}^{n+1/2})\cdot\nabla q,
\;     \forall (\mathbf{z},q) \in [\mathbb{P}^1(\mathcal{T}_h)]^d \times \mathbb{P}^1(\mathcal{T}_h)
  \end{aligned}
  \end{equation}
    with boundary conditions $\wb=\mathbf{z}=0$ on $\partial \Omega_{\mathcal T}$. 
Problem \eqref{eq:ste} can be equivalently formulated in matrix form as
  \begin{equation}\label{eq:ste}
    \begin{aligned}
     A_{\text{STE}}\,\wb + B \widehat{p}^{n+1/2}_{\text{STE}} &= M_{\text{STE}}\widehat{\ub}_{n+1} - M_{\text{STE}}\widehat{\ub}_{n} + Q_{\text{STE}}(\widehat{\ub}^{n+1/2}, \widehat{\ub}^{n+1/2}) + M_{\text{STE}}\widehat{\ub}_{n+1/2},\\
     B^T\wb &= 0,
  \end{aligned}
  \end{equation}
with natural definition of the stiffness matrix $A_{\text{STE}}$, of the mass matrix $M_{\text{STE}}$, of 
the matrix associated to the grad-div term $B$, and of the advection term $Q_{\text{STE}}$. 
 \end{problem} 
  
A bias correction for the STE can be obtained as solution to the following problem: Find $(\wb,b_{\text{STE}}) \in[\mathbb{P}^1(\mathcal{T}_h)]^d \times \mathbb{P}^1(\mathcal{T}_h)$ such that 
  \begin{equation}\label{eq:ste-bias}
    \begin{aligned}
    \int_{\Omega_{\mathcal T}} \nabla\wb:\nabla\mathbf{z} &- \int_{\Omega_{\mathcal T}} b_{\text{STE}}(\nabla\cdot\mathbf{z})+\int(\nabla\cdot\wb)q =\rho\sum_{T\in\mathcal{T}_h}\sum_{i=1}^{\widehat{d}_{\mathbf u}}\sum_{j=1}^{\widehat{d}_{\mathbf u}}\int_{\Omega_{\mathcal T}} \Psi_{i,j} \cdot\mathbf{z}.
   \end{aligned}
    \end{equation}

The bias corrections introduced here, have been extended to the general case of heteroscedastic noise starting from~\cite{bertoglio2018relative}.

\subsection{Numerical results}
\label{subsec:resultspressureestimators}
This section is dedicated to the comparison of the performance of the variational-based estimators (PPE and STE)
against the GNNs. The results are evaluated considering a global error on the pressure fluctuation on the whole target domain $\epsilon_{\widehat{p}}$ (equation~\eqref{eq:l2relerr}).

A further biomarker of clinical interest is the pressure drop across the coarctation. We consider two cross-sections $\Gamma^{\text{sec}}_\text{in}$ close to the inlet $\Gamma_{\rm in}$ and a cross-section $\Gamma^{\text{sec}}_4$ close to the outlet $\Gamma_4$ of the target geometries (the position of the cross-sections depends on the centerline encoding of the geometries, as in~\cite{katz2023impact}) and define the pressure drop as
\begin{equation}
  \widehat{p}_{4\text{-}\text{in}} = \frac{1}{\left|\Gamma^{\text{sec}}_4\right|}\int_{\Gamma^{\text{sec}}_4}\widehat{p}\,d\sigma-\frac{1}{|\Gamma^{\text{sec}}_\text{in}|}\int_{\Gamma^{\text{sec}}_\text{in}}\widehat{p}\,d\sigma,
\end{equation}
where $| \cdot |$ stands for the area of the corresponding surface. 


We will consider the pressure fields computed with the PPE and STE from the high-fidelity velocity $\widehat{\mathbf{u}}^{\text{true}}$ (\textit{ppe}, \textit{ste}), 
the pressure fields computed with the PPE and STE from the observed velocity field with values $\{l_i\}_{i=1}^{M_{\text{voxels}}}$ and support points $\{\mathbf{c}_i^{\text{vox}}\}_{i=1}^{M_{\text{voxels}}}$ (\textit{ppeobs}, \textit{steobs}), 
the pressure fields computed with the PPE and STE from the PBDW velocity $\widehat{\mathbf{u}}^{\text{true}}_{\text{PBDW}}$ (\textit{ppefom}, \textit{stefom}), and the pressure fields 
computed with the reduced order models of the PPE and STE from the PBDW velocity $\widehat{\mathbf{u}}^{\text{true}}_{\text{PBDW}}$ (\textit{pperom}, \textit{sterom}), as described in appendix~\ref{apendix:rom}. 
Additionally, the pressure inferred with GNNs from the geometrical encoding is denoted by \textit{gnn-gp}, the pressure inferred with GNNs from the high-fidelity velocity field $\widehat{\mathbf{u}}^{\text{true}}$ by \textit{gnn-vp}, and the pressure inferred with GNNs from the PBDW velocity field $\widehat{\mathbf{u}}^{\text{true}}_{\text{PBDW}}$ by \textit{gnn-vp-pbdw}. The inference problem \textit{gnn-vp} and \textit{gnn-vp-pbdw} employ the same architecture defined in section~\ref{ssec:pres-gnn}.

\begin{figure}[!ht]
  \centering
  \includegraphics[width=0.9\textwidth]{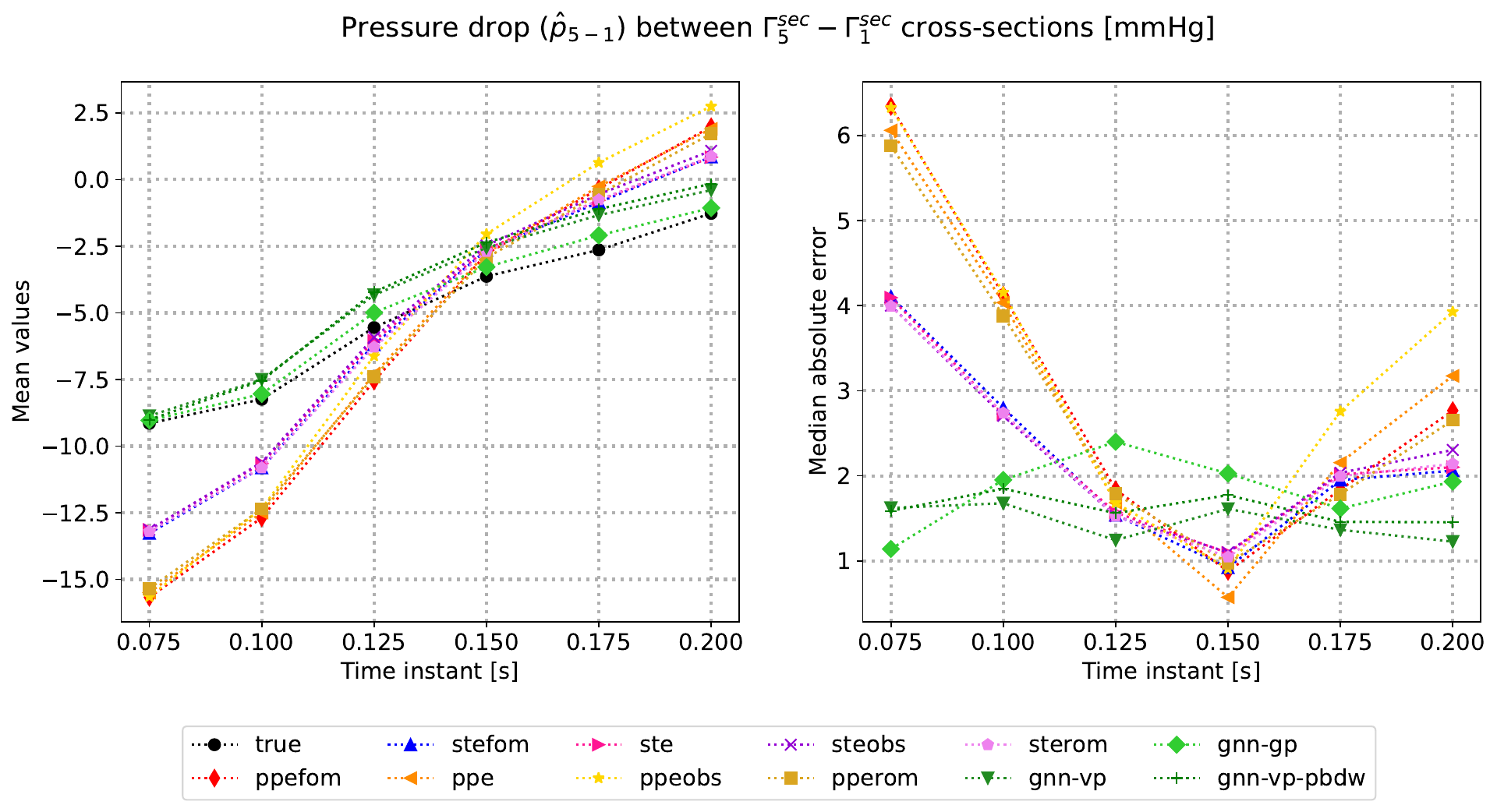}\\
  \includegraphics[width=0.45\textwidth]{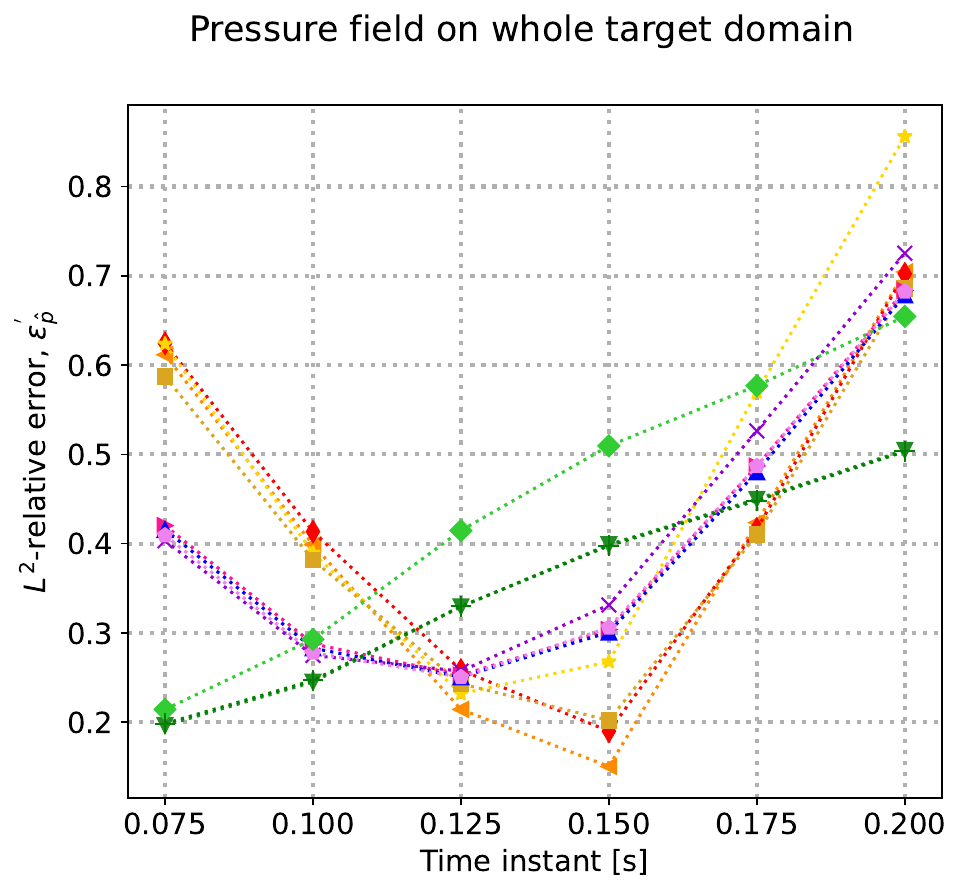}
  \caption{\textbf{Top left}: average value of the pressure drop $\overline{\widehat{p}_{4\text{-}\text{in}}-\widehat{p}_{\text{true},4\text{-}\text{in}}}$ between cross-sections $\Gamma_4^{\text{sec}}-\Gamma_\text{in}^{\text{sec}}$ over all the $52$ test geometries at time instants $t\in\{0.075s, 0.1s, 0.125s, 0.15s, 0.175s, 0.2s\}$. \textbf{Top right}: median of the absolute error $|\widehat{p}_{4\text{-}\text{in}}-\widehat{p}_{\text{true},4\text{-}\text{in}}|$ of the predicted pressure drop with respect to the high-fidelity pressure drop over all the $52$ test geometries. \textbf{Bottom: } average of the $L^2$-relative error $\epsilon_{\widehat{p}}$ defined in equation~\ref{eq:l2relerr} over all the $52$ test geometries. The description of the labels is reported in the text.}
  \label{fig:pressure_res}
\end{figure}

\begin{figure}[!ht]
  \centering
  \includegraphics[width=0.9\textwidth]{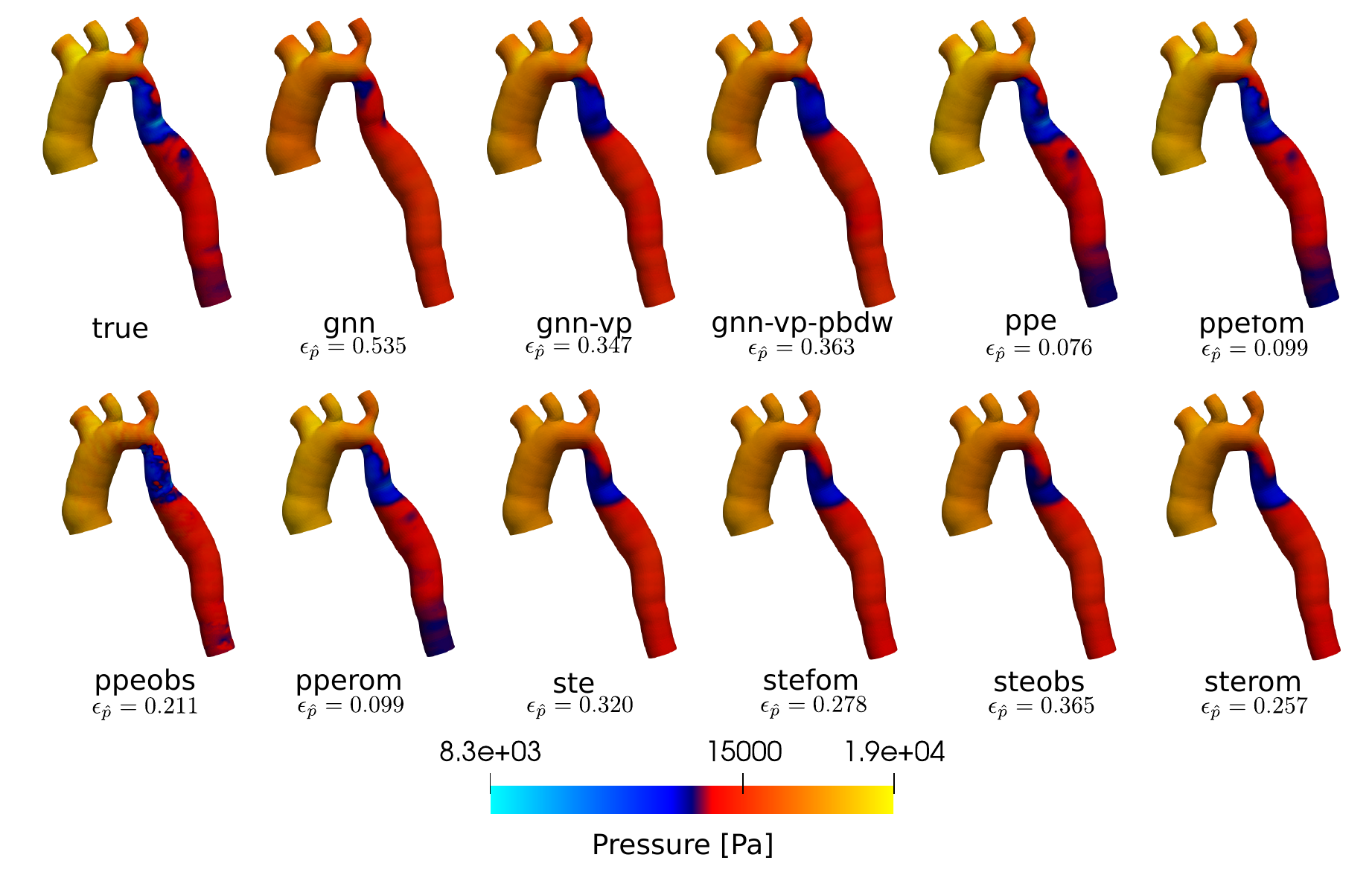}\\
  \caption{Predicted pressure at time $t=0.125s$ for test geometry $n=12$ from figure~\ref{fig:cluster_v} and associated $L^2$-relative errors $\epsilon_{\widehat{p}}$.
  }
  \label{fig:pres_12}
\end{figure}

The average pressure drops, median absolute pressure drops errors and the average $L^2$-relative errors (equation \eqref{eq:l2relerr}) across all test geometries are shown in figure~\ref{fig:pressure_res}.
For a qualitative comparison, the different predicted pressure fields are shown for the test geometry $n=12$ in figure~\ref{fig:pres_12}. 
In this case the results of the EPD-GNNs are comparable to those of PPE and STE. However, GNN are expected to deliver better results increasing the amount of training data,
due to the high geometric variability of the dataset: a symptom is the low training error in figure~\ref{fig:overfitting}. The test cases corresponding to the minimum, maximum and median $L^2$-relative errors are reported in figure~\ref{fig:gp} for the \textit{gnn-gp} problem. 

The best accuracy on the pressure field approximation and pressure drop corresponds to the time instants $t=0.125s$ and $t=0.15s$ for the PPE and STE estimators, the same time instants associated also to the best rSVD reconstruction error for the pressure fields in figure~\ref{fig:recerr}. However, the worse accuracy is not due to the rSVD reconstruction error but to the definition of the PPE and STE, as the same accuracy is associated to the estimated pressure fields \textit{ppe} and \textit{ste} obtained from the high-fidelity velocity. Moreover, the accuracy of PPE and STE does not seem to depend on the choice of observed velocity field: be it obtained from 4D-flow MRI observations (\textit{ppeobs}, \textit{steobs}), high-fidelity velocity fields (\textit{ppe}, \textit{ste}), PBDW velocity (\textit{ppefom}, \textit{stefom}), or reduced-order PPE and STE (\textit{pperom}, \textit{sterom}), the predictions achieve almost the same accuracy. Possibly, to improve the accuracy, the time resolution of $\Delta t=0.025s$ should be reduced to the high-fidelity simulations' time step $\Delta t=0.0025s$.

\subsection{Forward uncertainty quantification}
\label{subsec:uqpressureestimators}
Since PBDW models the uncertainty on the predicted velocity field from the coarse measurements $\widehat{\mathbf{u}}_{\text{PBDW}}\sim\mathcal{N}(m_{\widehat{\mathbf{u}}_{\text{PBDW}}}, \Sigma_{\widehat{\mathbf{u}}_{\text{PBDW}}})$, we want to study the uncertainty propagation to the pressure field through the pressure estimators and the inference with GNNs. In the following studies, we will keep the test geometry fixed and equal to test case $n=12$ from figure~\ref{fig:cluster_v}.

\begin{figure}[!htp]
  \centering
  \includegraphics[width=0.9\textwidth]{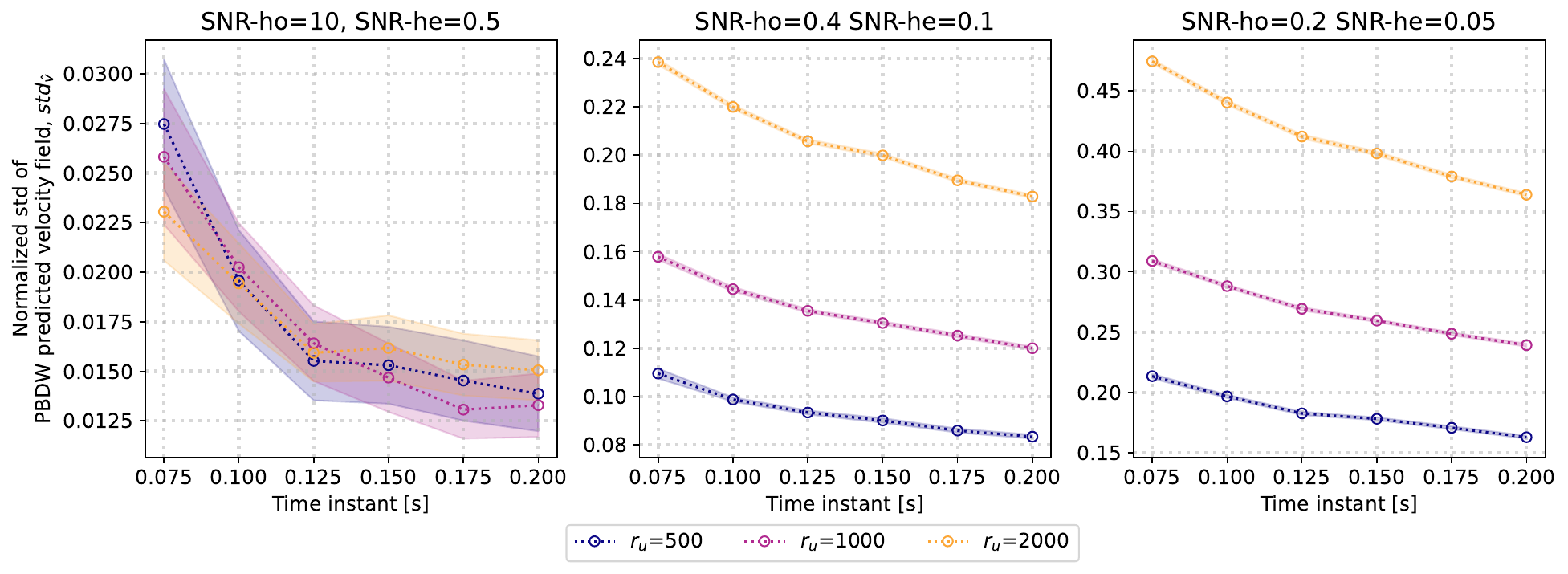}\\
  \caption{Normalized standard deviation ($\text{std}_{\widehat{\mathbf{u}}}$, equation~\ref{eq:std}) of PBDW velocity $\widehat{}_{\text{PBDW}}$ for different values of signal-to-noise ratios $(\text{SNR-ho}, \text{SNR-he})\in \{(10, 0.5), (0.4, 0.1), (0.2, 0.05)\}$ and number of velocity modes $r_{\mathbf u}\in\{500, 1000, 2000\}$. We consider only test case $n=12$, $n_{\text{samples}}=100$ velocity fields were sampled to compute $\text{std}_{\widehat{\mathbf{u}}}$.}
  \label{fig:snr}
\end{figure}

To measure the velocity field variability, we evaluate the normalized standard deviation at different values of $(\text{SNR-ho}, \text{SNR-he})\in \{(10, 0.5), (0.4, 0.1), (0.2, 0.05)\}$ in figure~\ref{fig:snr}:
\begin{equation}
  \label{eq:std}
  \text{std}_{\widehat{\mathbf{u}}} = \frac{\sum^{n_{\text{samples}}}_{i=1} \lVert \widehat{\mathbf{u}}_i-\tfrac{1}{n_{\text{samples}}}\sum_{i=1}^{n_{\text{samples}}}\widehat{\mathbf{u}}_i\rVert_2}{\sum^{n_{\text{samples}}}_{i=1} \lVert \tfrac{1}{n_{\text{samples}}}\sum_{i=1}^{n_{\text{samples}}}\widehat{\mathbf{u}}_i\rVert_2},
\end{equation}
where $n_{\text{samples}}=100$ are the number of samples from the Gaussian distribution $\widehat{\mathbf{u}}_{\text{PBDW}}\sim\mathcal{N}(m_{\widehat{\mathbf{u}}_{\text{PBDW}}}, \Sigma_{\widehat{\mathbf{u}}_{\text{PBDW}}})$ of the velocity predicted by PBDW with $r_{\mathbf u}\in\{500, 1000, 2000\}$ velocity modes.

Since the computational cost of a forward uncertainty problem is high due to the high number of forward evaluations, $n_{\text{samples}}=100$ in our case, we employ a reduced order model of the PPE and STE (\textit{pperom, sterom}), as described in appendix~\ref{apendix:rom}. We keep the number of pressure modes $r_p=1000$ fixed and vary the number of velocity modes $r_{\mathbf u}\in\{500, 1000, 2000\}$ and signal-to-noise ratios $(\text{SNR-ho}, \text{SNR-he})\in \{(10, 0.5), (0.4, 0.1), (0.2, 0.05)\}$.
\begin{rmk}
Notice that, unlike for classical reduced-order models, the assembly of the matrices cannot be performed offline, since the registration map is needed to transport the rSVD modes from the reference to the target geometry. 
However, once the registration map is evaluated and the rSVD bases have been transported to the target geometry, the assembly of the reduced systems can be performed in parallel and solved with less computational costs thanks to the lower dimensionality of the reduced systems. 
\end{rmk}
%

\begin{rmk}
For enhancing the stability of the reduced STE, the velocity rSVD modes are enriched with the supremizer technique~\cite{ballarin2015supremizer}. 
In this work, we evaluated the supremizers directly on the test geometries. The enrichment could be also performed offline on the reference geometry and transported on the new target, but this might
affect the overall stability of the formulation. 
\end{rmk}

\begin{figure}[!htp]
  \centering
  \includegraphics[width=0.65\textwidth]{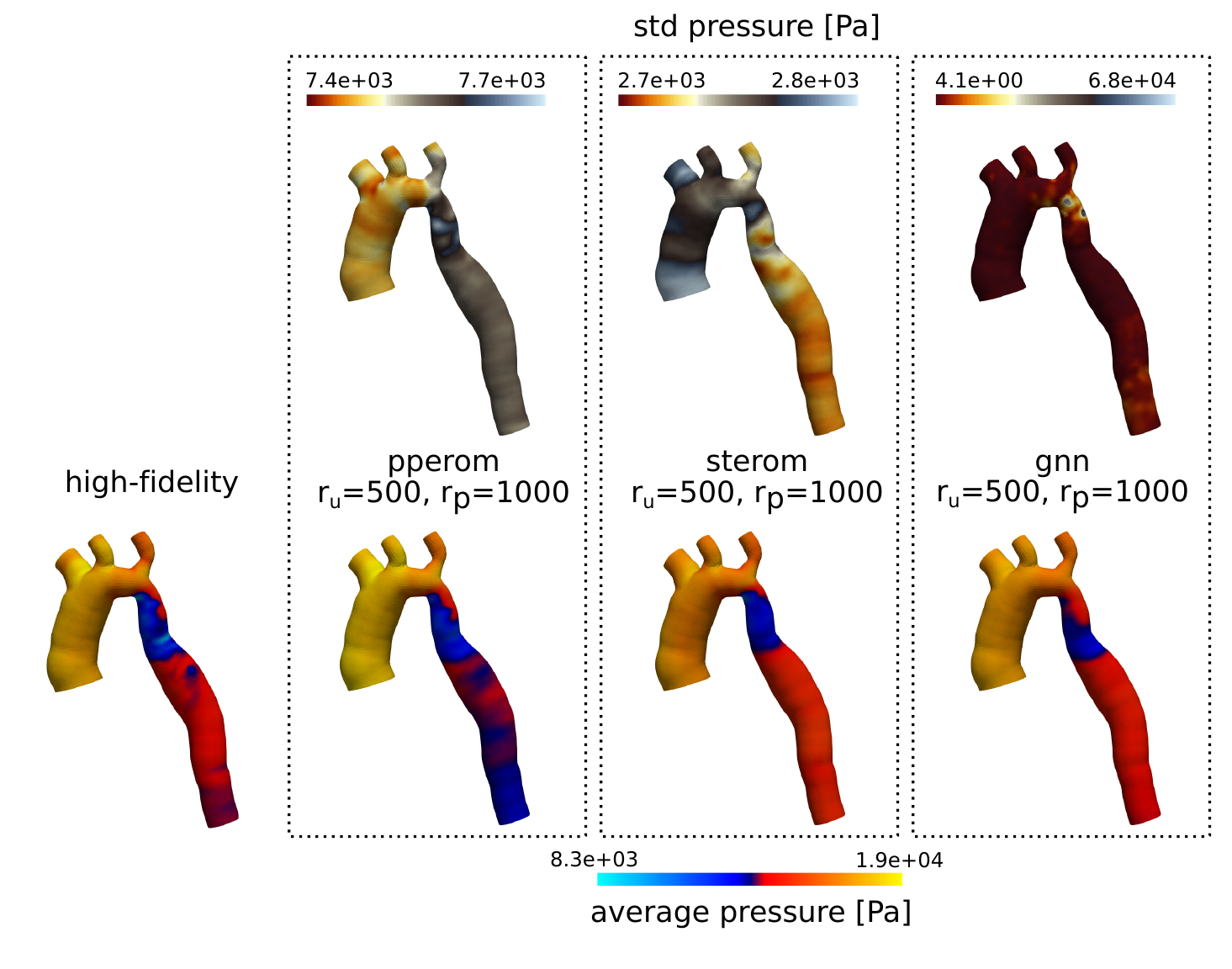}\\
  \caption{Average and standard deviation of pressure field over $n_{\text{samples}}=100$ samples, obtained forwarding the uncertainty with \textit{pperom, sterom} or \textit{gnn-vp-pbdw} from $n_{\text{samples}}=100$ PBDW velocity samples $\widehat{\mathbf{u}}_{\text{PBDW}}$. The test/target geometry is fixed $n=12$. The results correspond to the upper left block of figure~\ref{fig:uq_p}.}
  \label{fig:uq_p}
\end{figure}

In figure~\ref{fig:uq_p}, we compare the results of the PPE and STE with the predictions from the GNNs that compute the pressure field based on the PBDW velocity field as input, \textit{gnn-vp-pbdw}. The errors are evaluated on the test/target geometry with respect to the metrics in equation~\ref{eq:l2relerr}. The PPE estimator is omitted for $(\text{SNR-ho}, \text{SNR-he})\in \{(0.4, 0.1), (0.2, 0.05)\}$ as the relative error $\epsilon_{\widehat{\mathbf{u}}}$ goes above the value $1$ for all time instants.

\begin{figure}[!htp]
  \centering
  \includegraphics[width=0.85\textwidth]{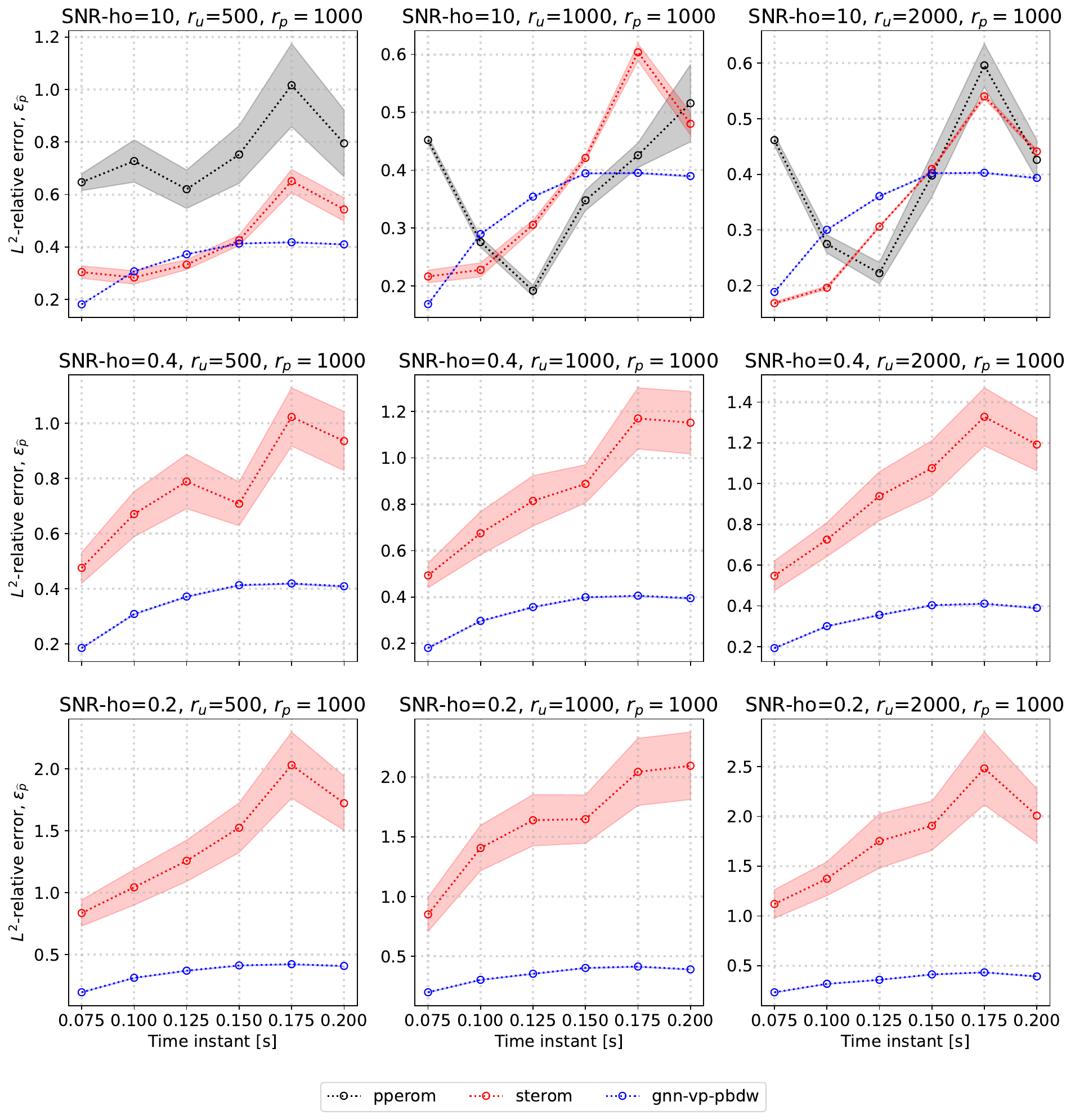}\\
  \caption{$L^2$-relative error $\epsilon_{\widehat{p}}$ and $0.95$ confidence intervals of the pressure predicted on the target geometry $n=12$ with PPE and STE and a GNN model, for different values of signal-to-noise ratios $(\text{SNR-ho}, \text{SNR-he})$ and number of velocity modes $r_{\mathbf u}\in\{500, 1000, 2000\}$. The pressure modes are fixed at $r_p=1000$.}
  \label{fig:uq_cfd}
\end{figure}

In comparison to PPE and STE results, the GNNs' predictions are robust (or rather overconfident), to the uncertainty in the velocity field. It can be shown also looking at the standard deviation of the pressure field predicted with  \textit{pperom, sterom} or \textit{gnn-vp-pbdw} from $n_{\text{samples}}=100$ PBDW velocity samples $\widehat{\mathbf{u}}_{\text{PBDW}}$, in figure~\ref{fig:uq_cfd}: the magnitude of the normalized standard deviation in equation~\eqref{eq:std} is high only on a small subdomain for the GNN models, underlying that perturbations of the inputs do not considerably affect the outputs.

\section{Discussion and limitations}
\label{sec:discussions}
\paragraph{Increasing the accuracy}
The correlations between velocity and pressure, and the comparison of the PPE and STE with GNNs (\textit{gnn-vp-pbdw}) suggest that GNNs can be used to infer the pressure from the velocity field in presence of a sufficiently high amount of data. Employing 4DMRI to acquire the velocity field observations should reduce the source of errors due to intrinsic uncertainties on the geometry and 
boundary condition acquisition and also on the choice of physical model and related approximations (i.e. rigid walls, turbulence model). The inference of the pressure and the velocity from the geometrical encoding (\textit{gnn-gp}, \textit{gnn-gv}) represents a more difficult task that may require a substantial increment of the available training dataset. Local linear and nonlinear dimension reduction techniques for solution manifold learning should be employed.

\paragraph{Limitations of the model} Apart from physical modelling assumptions, such as the neglecting fluid-structure interaction, the modelling of the aortic valve, and external tissue support on the vessel walls, we considered only parabolic inflow and a Windkessel model tuned according to an estimated flow split. These last two approximations, could be removed enriching the training dataset at the price of an additional computational cost as the intrinsic dimensionality of the solution manifold would also increase. Arbitrarily increasing the accuracy of the physical model does not necessarily results in more accurate predictions in a data assimilation context due to the epistemic uncertainties associated to the boundary conditions, to the physics, and to the geometry definition.

\paragraph{GNNs memory-bound} We employed coarse meshes for GNNs instead of the full-mesh due to memory constraints: GNNs built on graphs that correspond to meshes of large-scale applications have a high memory footprint limiting the size of the NN models. This is an active field of research: possible solutions include \textit{multigrid} strategies, distributed training, mini-batches subsamplers and gradient checkpointing. We remark that despite our limited computational budget, the low accuracy on the test dataset depends on the high geometric variability rather than on the employment of coarse meshes, as can be deduced from the high overfitting error and high training error on the coarse mesh in figures~\ref{fig:overfitting} and ~\ref{fig:overfitting_edges}.

\section{Conclusions}
\label{sec:conclusions}
In this work, we propose a robust shape registration method for complex geometries, such as aortic shapes, and its application to reduced-order modeling, data assimilation problems, and the training of EPD-GNNs for the inference of velocity and pressure fields from available observations.

The registration is based on a multigrid ResNet-LDDMM approach, trained using a shape database created with Statistical Shape Modeling (SSM) from an initial cohort of healthy subjects and aortic coarctation patients. The optimization is based on a modified Chamfer distance, tailored to computational meshes. By refining the surface mesh over the epochs, we show that realistic computational meshes can be efficiently handled.

The registration allows the definition of a bijective and non-parametric mapping between shapes, independent of mesh connectivity.
This enables the development of projection-based ROMs on different geometries after the pushforward of the SVD modes from a reference shape.
Our results showed that this might be challenging in practice due to the high number of velocity and pressure modes needed to achieve satisfactory accuracy.
However, in other applications with a larger amount of data or a solution manifold characterized by a lower intrinsic dimensionality, registration could be effectively used to develop physics-based surrogates, as in~\cite{guibert2014group}. 
We studied the correlation between geometric encodings and velocity/pressure fields on a common reference geometry,
emphasizing the impact of solution manifold learning metrics (Hausdorff vs. Grassmann). 

The employment of registration leads to a substantial improvement in the training of EPD-GNNs for the inference of pressure and velocity fields from geometry encoding, as well as for the inference of the pressure field from velocity data.
Next, we studied the application of the global rSVD basis in the context of the Parametrized-Background Data-Weak (PBDW) method.
Focusing on measurement data mimicking 4D MRI measurements, we proposed a natural extension of the PBDW formulation to heteroscedastic noise models, accounting for higher uncertainties in velocity gradients near vessel walls.
The results showed that physics-driven methods, such as PBDW, yield better results for velocity reconstruction than purely data-driven methods, such as EPD-GNN, which infers the velocity field solely from the geometry encoding.
Regarding the reconstruction of the pressure field and pressure drops from velocity data, the results of EPD-GNN—whether from coarse observations or in combination with the PBDW reconstruction of a high-resolution velocity field—outperformed widely used approaches such as the Pressure-Poisson estimator (PPE) and the Stokes estimator (STE).
Finally, we explored the use of the PPE and STE, combined with projection-based reduced-order modeling, for forward uncertainty quantification.

For relevance in clinical applications and image-based diagnostic, surrogates models and data assimilation methods shall be able to predict biomarkers of interest accurately and efficiently, 
validating the results with medical data. Towards this direction, in the specific context of estimating pressure drops in patients affected by aortic coarctation, 
the high variability of the three-dimensional fluid dynamics requires a higher amount of data. Reducing the surrogate models to one-dimensional centerline supports instead of predicting 
full three-dimensional fields might yield a more efficient and sufficiently accurate approach. 
However, although three-dimensional models are not more computationally expensive and with a higher intrinsic dimensionality, these are also more interpretable and have the possibility to evaluate
additional quantity of interests such as secondary flow degree, normalized flow displacement, as well as wall shear stresses.

Our pipeline could be extended to other anatomical shapes (e.g., heart or liver), as well as beyond biomedical applications (e.g., for shape optimization in general context). 
From the point of view of surrogate modelling, different approaches to tackle the shape variability could be compared considering the reconstruction errors on the template geometry, thus measuring the complexity of the different solution manifolds. 
In this way, approaches that are effective on test cases that require few SVD modes to be well-approximated, but struggle on more difficult ones, could be identified.

Future directions of research include greedy data augmentation strategies to efficiently increase the training dataset and design better solution manifold approximants, the implementation of nonlinear shape generative models from machine learning instead of SSM (which is a linear method based on PCA), and the optimization of EPD-GNNs architectures for large-scale applications, with a focus on memory efficiency.

\section*{Acknowledgements}
This research has been partially funded by the Deutsche Forschungsgemeinschaft (DFG, German Research Foundation) under Germany's Excellence Strategy - MATH+: the Berlin Mathematical Research Center [EXC-2046/1 - project ID: 390685689], project AA5--10 \textit{Robust data-driven reduced-order models for cardiovascular
imaging of turbulent flows}.

\appendix
\setcounter{equation}{0}
\renewcommand\theequation{\arabic{equation}}
\section{Numerical solution of the 3D-0D blood flow model}
\label{appendix:weak}
At each time step $t_{n+1}\in[2\Delta t,0.7s]$ we solve the following linear system for $(\mathbf{u}_{n+1}, p_{n+1}, \{\pi_{i,n+1}\}_{i=1}^{4})$, given the solution for these quantities at time
steps $t_{n}$ and $t_{n-1}$. 

\begin{problem}[Weak formulation]
  \label{def:weak_ins_w}
  Find $(\mathbf{u}_{n+1}, p_{n+1}, \{\pi_{i,n+1}\}_{i=1}^{4})\in[\mathbb{P}^1(\mathcal{T}_h)]^3\times\mathbb{P}^1(\mathcal{T}_h)\times\mathbb{R}^4$ such that $\forall(\mathbf{v},q)\in [\mathbb{P}^1(\mathcal{T}_h)]^3\times\mathbb{P}^1(\mathcal{T}_h)$:
  \begin{subequations}\begin{align*}
    \rho\left(\frac{\tfrac{3}{2}\mathbf{u}_{n+1}-\mathbf{u}_{n}+\tfrac{1}{2}\mathbf{u}_{n-1}}{\Delta t},\mathbf{v}\right)_{\Omega}&+\left(2\mu\nabla^S\mathbf{u}_{n+1}, \nabla^S\mathbf{v}\right)_{\Omega} + \rho\left((2\mathbf{u}_n-\mathbf{u}_{n-1})\cdot\nabla\mathbf{u}_{n+1},\mathbf{v}\right)_{\Omega}+\\
    &-(p_{n+1},\nabla\cdot\mathbf{v})_{\Omega}+(\nabla\cdot\mathbf{u}_{n+1}, q)_{\Omega}+\\
    &+\mathcal{S}(\mathbf{u}_{n+1}, \mathbf{u}_{n}, \mathbf{u}_{n-1},p_{n+1}, p_{n}, p_{n-1},\mathbf{v}, q) =\\
    &= \sum_{i=1}^4\left(R_{p,i}\int_{\Gamma_i}\mathbf{u}_{n}\cdot\mathbf{n}+\pi_{i,n+1},\mathbf{v}\right)_{\Gamma_i}+\tfrac{\rho}{2}\left(\min\{(2\mathbf{u}_{n}-\mathbf{n}_{n-1})\cdot\mathbf{n},0\},\mathbf{v}\right)_{\Gamma_i},\\
    \text{for}\quad i\in\{1, 2,3,4\},\qquad &\int_{\Gamma_i}\mathbf{u}_{n}\cdot\mathbf{n}=C_{d,i}\frac{\pi_{i, n+1}-\pi_{i, n}}{\Delta t}+\frac{\pi_{i,n}}{R_{d,i}},
  \end{align*}\end{subequations}
  where $\rho,\mu, \{C_{d, i}\}_{i=1}^4, \{R_{d, i}\}_{i=1}^4, \{R_{p, i}\}_{i=1}^4$ are positive constants. The stabilization $\mathcal{S}$ includes the SUPG-PSPG and VMS-LES terms:
  \begin{subequations}\begin{align*}
    &\mathcal{S}(\mathbf{u}_{n+1}, \mathbf{u}_{n}, \mathbf{u}_{n-1},p_{n+1}, p_{n}, p_{n-1};\mathbf{v}, q) = \sum_{K\in\mathcal{T}_h}(\mathcal{S}_{\text{SUPG-PSPG}}^K+\mathcal{S}_{\text{VMS}}^K+\mathcal{S}_{\text{LES}}^K),\\
    &\mathcal{S}_{\text{SUPG-PSPG}}^K=(\tau_M\cdot \mathbf{r}_M(\mathbf{u}_{n+1}, p_{n+1}),\rho \mathbf{u}_{n+1}^{*}\cdot\nabla\mathbf{v}+\nabla q)_K+(\tau_C\cdot r_C,\nabla\cdot\mathbf{v})_K,\\
    &\mathcal{S}^K_{\text{VMS}} = (\tau_M\cdot \mathbf{r}_M(\mathbf{u}_{n+1}, p_{n+1}),\rho \mathbf{u}_{n+1}^{*}\cdot\nabla^T\mathbf{v})_K,\\
    &\mathcal{S}^K_{\text{LES}} = -(\tau_M \mathbf{r}_M(\mathbf{u}^{\text{EXT}}_{n+1}, p^{\text{EXT}}_{n+1})\otimes\tau_M \mathbf{r}^{\text{LHS}}_M(\mathbf{u}_{n+1}, p_{n+1}),\rho\nabla\mathbf{v})_K+(\tau_M \mathbf{r}_M(\mathbf{u}_{n+1}, p_{n+1})\otimes \tau_M \mathbf{r}^{\text{RHS}}_M, \rho\nabla\mathbf{v})_K,\\
      \end{align*}\end{subequations}
    with
\begin{subequations}\begin{align*}    
    &\mathbf{r}_M(\mathbf{u}, p) = \mathbf{r}^{\text{LHS}}_M(\mathbf{u}, p)-\mathbf{r}^{\text{RHS}}_M,\\
    &\mathbf{r}^{\text{LHS}}_M(\mathbf{u}, p) =\rho\left(\frac{\tfrac{3}{2}\mathbf{u}}{\Delta t}+\mathbf{u}^{*}_{n+1}\cdot\nabla\mathbf{u}\right)+\nabla p-\mu\Delta\mathbf{u},\\
    &\mathbf{r}^{\text{RHS}}_M=\rho\frac{\mathbf{u}_{n}-\tfrac{1}{2}\mathbf{u}_{n-1}}{\Delta t},\\
    &r_C = \nabla\cdot\mathbf{u}_{n+1},\\
    &\tau_M = \left(\frac{4\rho^2}{\Delta t^2}+\rho^2\mathbf{u}^{*}_{n+1}\cdot G\mathbf{u}^{*}_{n+1}+30\mu^2 G:G\right)^{\tfrac{1}{2}},\\
    &\tau_C = (\tau_M\mathbf{g}\cdot\mathbf{g})^{-1},\\
  \end{align*}\end{subequations}
  with $\mathbf{u}^{*}_{n+1} = 2\mathbf{u}_{n}-\mathbf{u}_{n-1}$, $p^{*}_{n+1}=2p_{n}-p_{n-1}$, $\mathbf{u}^{\text{EXT}}_{n+1}=2\mathbf{u}_{n}-\mathbf{u}_{n-1}$,$p^{\text{EXT}}_{n+1}=2p_{n}-p_{n-1}$, and where $G=J^{-T}J^{-1}$ and $\mathbf{g}=J^{-T}\mathbf{1}$ denote the metric tensor and vector obtained from the Jacobian $J$ of the map from reference element to current finite element.
\end{problem}
The first time step $t=\Delta t$ is solved with an explicit Euler scheme. We use an uniform time step with $\Delta t=\SI{2.5e-03}{\second}$.

Additional information can be found in the documentation of \texttt{lifex-cfd}~\cite{AFRICA2024109039} and references~\cite{Fumagalli2020, Fedele2017, bazilevs2007variational, forti2015semi}. To solve the linear system, we employ GMRES with the SIMPLE preconditioner, see the documentation \texttt{lifex-cfd}~\cite{AFRICA2024109039} and references~\cite{Deparis2014}: the AMG preconditioner is used to invert the Shur complement and the $(\mathbf{u}_{n+1},\mathbf{v})$ matrix block in the discretized block matrix obtain from definition~\ref{def:weak_ins_w}, neglecting the Windkessel model.

\begin{figure}[!htp]
  \centering
  \includegraphics[width=1\textwidth,trim=0cm 0cm 0cm 1cm, clip=true]{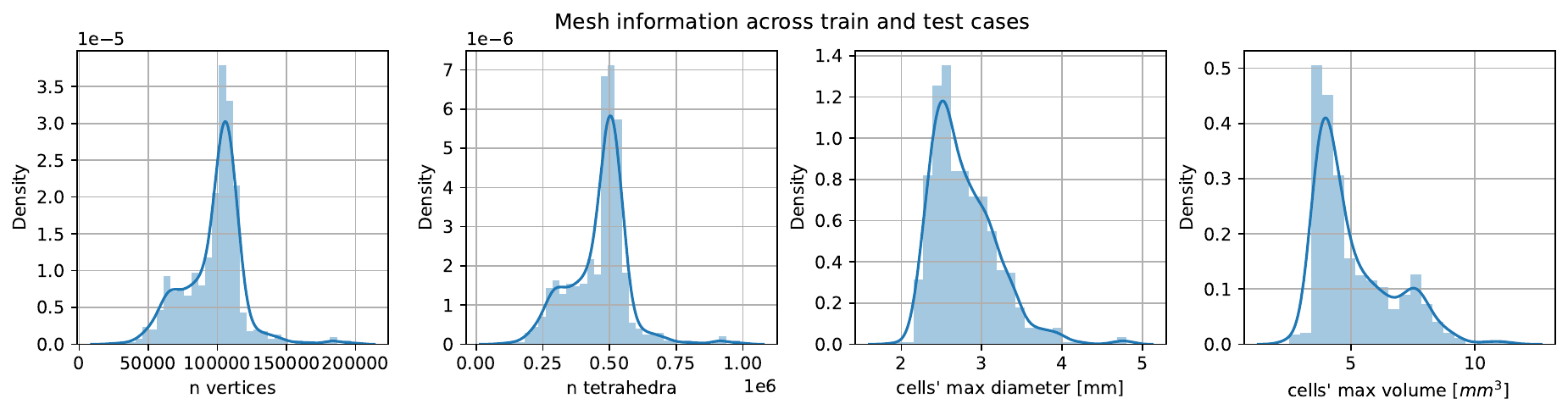}
  \caption{Mesh information across the $776$ train and test geometries.}
  \label{fig:mesh_info}
\end{figure}

\begin{figure}[!htp]
  \centering
  \includegraphics[width=0.85\textwidth]{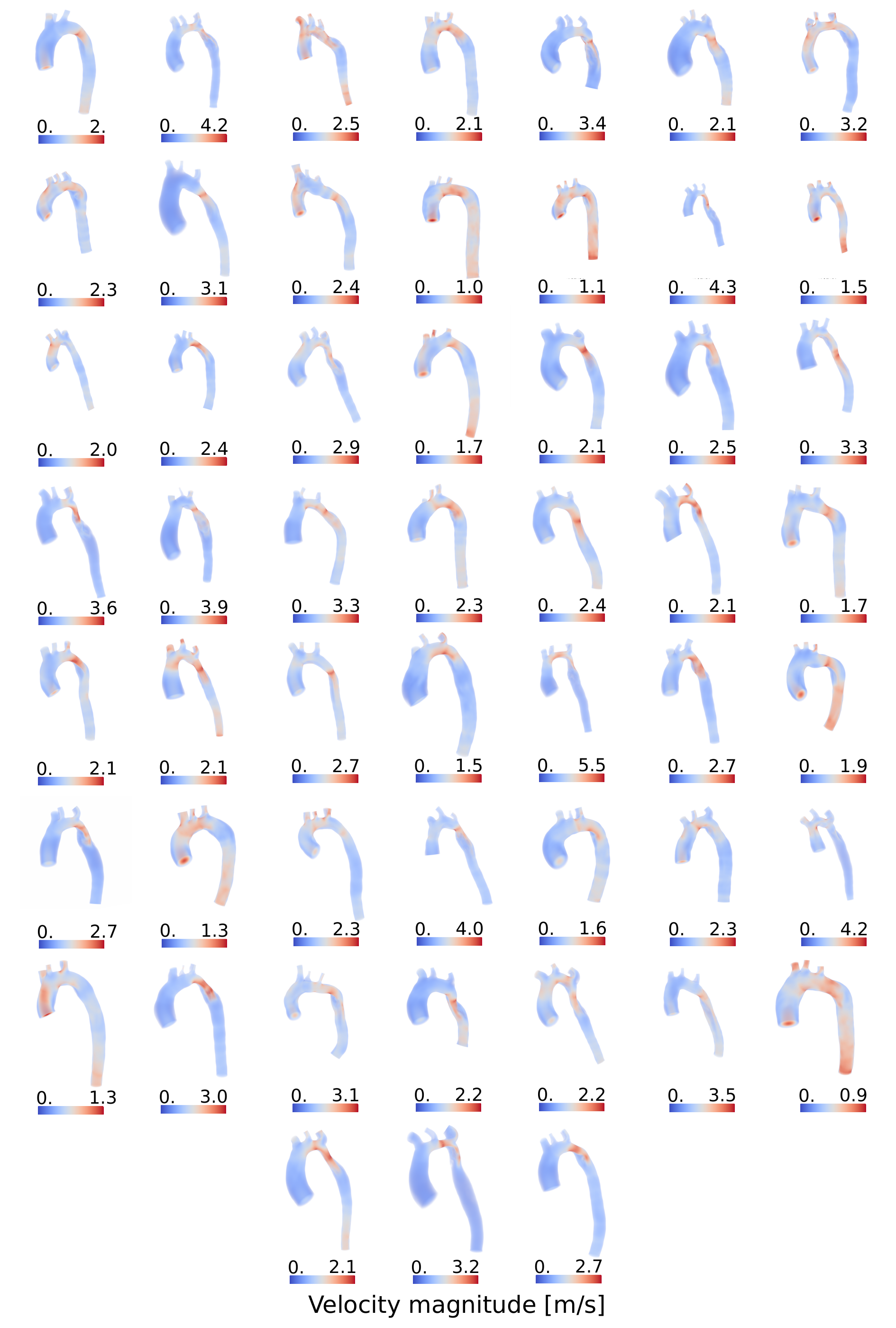}
  \caption{High-fidelity velocity fields at systolik peak $t=0.125s$ on the $52$ test geometries.}
  \label{fig:snapsu}
\end{figure}

The meshes are generated with VMTK~\cite{antiga2008image} from surface representations obtained with the ZIB Amira software~\cite{stalling2005amira}. Through VMTK one boundary layer with thickness factor $0.3$ is added. The number of vertices of each mesh, the number of tetrahedra of each mesh, the maximum diameter across all tetrahedra of each mesh and the maximum volume across all tetrahedra of each mesh is reported in figure~\ref{fig:mesh_info}: the distribution of these quantities is reported over all train and test cases for a total of $776$ geometries. Some velocity snapshots are shown in figure~\ref{fig:snapsu}.

The CFD simulations were run in parallel on $64$ cores in a distributed-memory fashion. Multiple compute servers at Weierstraß Institute for Applied Analysis and Stochastic (WIAS, Berlin) were employed to run in serial the CFD simulations for each training and test geometry: HPE Synergy 660 Gen10 with four Intel Xeon Gold CPUs each with 18 cores and HPE Synergy 480 Gen10 with two Intel Xeon Gold CPUs each with 18 cores. The models of the Intel Xeon Gold CPUs vary.

\section{Convergence of the discrete registration problem}
\label{appendix:convergence}
Proof of theorem~\ref{theo:existreg}.
\begin{proof}
For all $0\leq r\leq s$, it holds
  \begin{align*}
    \lVert\phi_\epsilon(t)-\phi(t)\rVert_{H^r}&\leq \int_0^t \lVert f_\epsilon\circ\phi_\epsilon-f\circ\phi_\epsilon\rVert_{H^r}+\lVert f\circ\phi_\epsilon-f\circ \phi\rVert_{H^r}\ d\tau \\
    &\lesssim \int_{0}^t\lVert f_\epsilon-f\rVert_{H^r}+\lVert f\rVert_{H^r}\lVert\phi_\epsilon-\phi\rVert_{H^r}\ d\tau,
  \end{align*}
  thanks to the fact that $\phi_\epsilon$ is a diffeomorphism.
  By Gronwall's inequality it also holds
  \begin{equation}
    \label{eq:phif}
    \lVert \phi_\epsilon - \phi\rVert_{L^2(I, H^{r})}\rVert \leq C_1\lVert f_\epsilon - f\rVert_{L^2(I, H^{r})},  \;\forall r,\  0\leq r\leq s.
  \end{equation}
  The first thesis follows for the universal approximation theorem of ReLU neural networks. The second estimate follows from theorem 4.1~\cite{doi:10.1142/S0219530522500014}.
\end{proof}
We study the convergence of the discrete registration problem. We define the functional $F:L^2(I, H^s)\rightarrow\mathbb{R}$,
\begin{align}
  \label{eq:funccont}
   F(f) &= \int_{0}^{1}\lVert f(t, \cdot)\rVert^2_{H^s}\ dt + \int_G |S(\phi(1, \x))-T(\x)|^2\ d\x,\\\
    \phi(t, \x)&=\x + \int_{0}^{t}f(s, \x)\,ds,\qquad (t,\x)\in I\times G
\end{align}
with $S, T$ as in theorem~\ref{def:regpb}, and the approximants $F_N:L^2(I, H^s)\rightarrow\mathbb{R}$, without loss of generality with the same discretization parameter $N\in\mathbb{N}$ in time and space,
\begin{align}
  \label{eq:funcdis}
  F_N(f) &=\begin{cases}
    \frac{1}{N}\sum_{i=1}^{N}\lVert f(\tfrac{1}{N}, \cdot)\rVert^2_{H^s} + \frac{1}{N}\sum_{j=1}^{N} |S(\phi_N(1, \x_j))-T(\x_j)|^2,\qquad &f\in \mathcal{C}^{0, 1}(I, \mathcal{C}^{0, 1})\\\
    +\infty,\qquad &f\notin \mathcal{C}^{0, 1}(I, \mathcal{C}^{0, 1})
  \end{cases}\\
  \phi_N(t, \x)&=\x + \int_{0}^{t}f(s, \x)\,ds,\qquad (t,\x)\in I\times G,
\end{align}
where $(\mathbf{x}_j)_{j=1}^N$ are sampled in $G$ such that $\sum_{j=1}^{N} |S(\phi_N(1, \x_j))-T(\x_j)|^2\chi_{A_j}$ are simple functions converging to $|S(\phi_N(1, \x))-T(\x)|^2$ with supports $(A_j)_{j=1}^{N}\subset G$ of Lebesgue measure $\tfrac{1}{N}=\lambda_{\mathcal{L}}(A_j),\ \forall j$.

The next result states the $\Gamma$-convergence of $F_N$ to $F$.
\begin{theorem} Let $F_N$ and $F$ be defined as in equations~\eqref{eq:funccont} and~\eqref{eq:funcdis}, then
 \begin{equation*}
  F = \Gamma{\text -}\!\lim_{N\to\infty} F_N
 \end{equation*}
\end{theorem}
\begin{proof}
  To obtain the $\Gamma$-convergence it is sufficient to prove the limsup and liminf inequalities:
  \begin{align*}
    \forall f\in L^2(I, H^s),\ \exists (f_N)_N: f_N\to f,\quad \text{s.t.}\quad &F(f)\geq \limsup_{N\to\infty}F_N(f_N),\\
    \forall f\in L^2(I, H^s),\ \text{for every}\ (f_N)_N,\ f_N\to f,\qquad &F(f)\leq\liminf_{N\to\infty} F_N(f_N).
  \end{align*}
  The function $S\circ\phi(1, \cdot)$ is measurable thanks to the fact that $\phi$ is a diffeomorphism, so the second term of the continuous functional in equation~\eqref{eq:funccont} is well-defined. Let us prove the limsup inequality with a recovery sequence first. Due to the density of $\mathcal{C}^{0, 1}(I, \mathcal{C}^{0, 1})$ in $L^2(I, H^s)$ for every $f\in L^2(I, H^s)$ there exists $(f_N)_N\subset\mathcal{C}^{0, 1}(I, \mathcal{C}^{0, 1})$ s.t. $f_N\to f$. Thanks to the convergence of the norms, we have that the first term in equations~\eqref{eq:funccont} and~\eqref{eq:funcdis} converges. To prove that the second term converges, we need that $\phi_N\rightarrow \phi$ thanks to equation~\eqref{eq:phif}. Up to passing to subsequences, since $S$ is continuous a.e., we have the pointwise convergence a.e. of $|S(\phi_N(1, \x))-T(\x)|^2$ to $|S(\phi(1, \x))-T(\x)|^2$. Since $S, T$ and $G$ are bounded, by the dominated convergence theorem we obtain also the convergence of the second term of the functionals in~\eqref{eq:funccont} and~\eqref{eq:funcdis}.
  
  We now prove the liminf. Thanks to the convergence of the norms, we have the convergence to the first terms of the functionals in equation~\eqref{eq:funccont} and~\eqref{eq:funcdis}. The convergence $f_N\to f$ in $L^2(I, H^s)$ implies, through equation~\eqref{eq:phif}, $\phi_N\to \phi$ in $L^2(I, H^s)$ which implies, by Sobolev embedding for $s>d/2 + 1$, $\phi_N\to \phi$ in $\mathcal{C}^{1, \alpha}(I, \mathcal{C}^{1, \alpha})$ with $\alpha=s - \tfrac{d}{2} - 1>0$. Using that $S$ is continuous a.e., by Fatou's lemma,
\begin{align*}
  F(f)&=\lVert f\rVert^2_{L^2(I, H^s)}+\int_G |S(\phi(1, \x))-T(\x)|^2\ d\x \\
  &= \lVert f\rVert^2_{L^2(I, H^s)}+\int_G \liminf_{N\to\infty}|S(\phi_N(1, \x))-T(\x)|^2\ d\x\\
   &\leq\lim_{N\to\infty}\left(\int_G |S(\phi_N(1, \x))-T(\x)|^2\ d\x-\frac{1}{N}\sum_{j=1}^{N} |S(\phi_N(1, \x_j))-T(\x_j)|^2\right) + \liminf_{N\to\infty} F_N(f_N),
\end{align*}
by dominated convergence the lim of the right-hand side converges.
\end{proof}
  
  \begin{corol} [Corollary 7.20, \cite{dal2012introduction}]
    For every $N\in\mathbb{N}$, let $f_N$ be a minimizer (or $\epsilon_N$-minimizer) of $F_N$ in $L^2(I, H^s)$. If $f$ is a cluster point of $(f_N)_N$, then $f$ is a minimizer of $f$ and $F(f)=\limsup_{N\to\infty}F_N(f_N)$. If $(f_N)_N$ converges to $f$ in $L^2(I, H^s)$, then $f$ is a minimizer of $F$ and $F(f)=\lim_{N\to\infty}F_N(f_N)$.
  \end{corol}

\section{PBDW with heteroscedastic noise}
\label{appendix:pbdw}
\begin{proof}{Theorem~\ref{theo:pbdw}}
  For a fixed $z\in\mathbb{R}^{r_{\mathbf u}}$ in~\eqref{eq:pbdw_hetero}, the minimum with respect to $\eta$ is given by
  \begin{equation}\label{eq:pbdw_eta_sub}
\left(R^{-1}+KS^{-1}K\right)\pbdw{\eta} + KS^{-1}(Lz-y) = 0,
    \end{equation}
    i.e., the minimizer of \eqref{eq:pbdw_sub2}.
From \eqref{eq:pbdw_eta_sub}, using $R=SK^{-1}$, it also follows 
  \begin{equation}\label{eq:pbdw_eta_sub3}
     \pbdw{\eta} = \left[ SK^{-1}\left(R^{-1}+KS^{-1}K\right) \right]^{-1}(y-Lz) = \left( \text{Id} + K \right)(y-Lz) = W (y-Lz)\,.
    \end{equation}
Inserting \eqref{eq:pbdw_eta_sub3} into~\eqref{eq:pbdw_hetero}, it follows that $\pbdw{z}$ minimizes the function
  \begin{equation*}
  \begin{aligned}
    & (y-Lz)^TW^{-1}(R^{-1}+(SK^{-1}R^{-1})^TS^{-1}(SK^{-1}R^{-1}))W^{-1}(y-Lz) \\
    & = (y-Lz)^TW^{-1}\underbrace{\left(KS^{-1}+S^{-1}\right)}_{W S^{-1}}W^{-1}(y-Lz) = \lVert Lz-y\rVert^2_{S^{-1}W^{-1}},
  \end{aligned}
    \end{equation*}
    which completes the proof.
\end{proof}

\begin{rmk}
  The linear estimator $H_{\widehat{\mathbf{u}}}$ is biased:
  \begin{equation*}
    \mathbb{E}[H_{\widehat{\mathbf{u}}}(y+\epsilon_{y})-\Tilde{H}_{\widehat{\mathbf{u}}_{\text{PBDW}}}y] = [\Phi_r (H_{z_{\text{PBDW}}}-\Tilde{H}_{\Tilde{z}}) + \mathcal{Z}_{\mathbf u} (H_{\eta_{\text{PBDW}}}-\Tilde{H}_{\Tilde{\eta}}) + \mathcal{Z}_{\mathbf u} (H_{\eta_{\text{PBDW}}}LH_{z_{\text{PBDW}}}-\Tilde{H}_{\Tilde{\eta}}L\Tilde{H}_{\Tilde{z}})]y,
  \end{equation*}
  where 
  \begin{align*}
    &\Tilde{H}_{\widehat{\mathbf{u}}_{\text{PBDW}}} = \Phi_r \Tilde{H}_{\Tilde{z}}+\mathcal{Z}_{\mathbf u} \Tilde{H}_{\Tilde{\eta}}-\mathcal{Z}_{\mathbf u} \Tilde{H}_{\Tilde{\eta}}L\Tilde{H}_{\Tilde{z}}\\
    &\Tilde{H}_{\Tilde{z}}=H_{\Tilde{z}}y = (L^T K^{-1} L)^{-1} L^T K^{-1}\\
    &\Tilde{H}_{\Tilde{\eta}}= K^{-1}
  \end{align*}
\end{rmk}
\begin{proof}{Theorem~\ref{theo:pbdwmsq}}
  From the definition of covariance (equation~\eqref{eq:cov}) of the PBDW prediction $\widehat{\mathbf{u}}_{\text{PBDW}}$:
  \begin{align*}
    \mathbb{E}[\lVert \widehat{\mathbf{u}}^{\text{true}}-\widehat{\mathbf{u}}_{\text{PBDW}}\rVert^2_2]&\leq \mathbb{E}[\lVert \widehat{\mathbf{u}}^{\text{true}}-\mathbb{E}[\widehat{\mathbf{u}}_{\text{PBDW}}]\rVert^2_2]+\mathbb{E}[\lVert \mathbb{E}[\widehat{\mathbf{u}}_{\text{PBDW}}]-\widehat{\mathbf{u}}_{\text{PBDW}}\rVert^2_2]\\
    &= \lVert \widehat{\mathbf{u}}^{\text{true}}-\mathbb{E}[\widehat{\mathbf{u}}_{\text{PBDW}}]\rVert^2_2+\text{trace}(H_{\widehat{\mathbf{u}}_{\text{PBDW}}}SH_{\widehat{\mathbf{u}}_{\text{PBDW}}}^T).
  \end{align*}
  We consider the operator $H_l = H_{\widehat{\mathbf{u}}^{\text{PBDW}}}\circ l:\mathbb{R}^{\widehat{d}_{\mathbf u}}\rightarrow\mathbb{R}^{\widehat{d}_{\mathbf u}}$,
  \begin{align*}
    \lVert \widehat{\mathbf{u}}^{\text{true}}-\mathbb{E}[\widehat{\mathbf{u}}_{\text{PBDW}}]\rVert&=\lVert \widehat{\mathbf{u}}^{\text{true}}-H_{\widehat{\mathbf{u}}^{\text{PBDW}}}(l(\widehat{\mathbf{u}}^{\text{true}}))\rVert \\
    &= \lVert \left((\text{Id}-H_l)\circ P_{\text{Im}(H_l)}\right)(\widehat{\mathbf{u}}^{\text{true}})\rVert_2+\lVert (\text{Id}-H_l)(P_{\text{Im}(H_l)^{\perp}}\widehat{\mathbf{u}}^{\text{true}})\rVert_2\\
    &\leq \lVert(\text{Id}-H_l)\circ P_{\text{Im}(H_l)}\rVert_2\lVert \widehat{\mathbf{u}}^{\text{true}}\rVert_2+\lVert\text{Id}-H_l\rVert_2 \lVert P_{\text{Im}(H_l)^{\perp}}\widehat{\mathbf{u}}^{\text{true}}\rVert_2
  \end{align*}
  Summing and subtracting the terms $(\phi_{\text{RBF}})^{\#}({\mathbf u}^{\text{best}})$, $P_{{\text{Im}(H_l)}}(\phi_{\text{RBF}})^{\#}({\mathbf u}^{\text{best}})$, and $(\phi_{\text{RBF}})^{\#}(P_{\text{Im}(\Phi_{\mathbf u})} {\mathbf u}^{\text{best}})$:
  \begin{align*}
    \lVert P_{\text{Im}(H_l)^{\perp}}\widehat{\mathbf{u}}^{\text{true}}\rVert_2 &\leq\lVert(\phi_{\text{RBF}})^{\#}({\mathbf u}^{\text{best}})-\widehat{\mathbf{u}}^{\text{true}}-P_{{\text{Im}(H_l)}}((\phi_{\text{RBF}})^{\#}({\mathbf u}^{\text{best}})-\widehat{\mathbf{u}}^{\text{true}})\rVert_2\\
    &+\lVert (\phi_{\text{RBF}})^{\#}({\mathbf u}^{\text{best}})-(\phi_{\text{RBF}})^{\#}(P_{\text{Im}(\widehat{\Phi}_{u})} {\mathbf u}^{\text{best}})\rVert_2\\
    &+\lVert P_{\text{Im}(H_l)}(\phi_{\text{RBF}})^{\#}({\mathbf u}^{\text{best}})-(\phi_{\text{RBF}})^{\#}(P_{\text{Im}(\Phi_{\mathbf u})}{\mathbf u}^{\text{best}})\rVert_2,
  \end{align*}
  we remark that $\text{Im}(H_l) = \text{Im}(\widehat{\Phi}_{u})+\left(\text{Im}(\widehat{\Phi}_{u})^{\perp}\cap\text{Im}(\mathcal{Z}_{\mathbf u})\right)$. The RBF interpolation operator with thin-plate basis $I_{\text{RBF}}:\mathbb{R}^3\rightarrow\mathbb{R}^3$ is Lipshitz with constant $C_{\text{RBF}}$. Setting $C=C_{\text{RBF}}\sup_{\mathbf{x}\in \Omega_{S}} \lVert\text{det}(\nabla \phi_1(\mathbf{x}))\rVert_2$:
  \begin{align*}
    \lVert (I_{\text{RBF}}\circ\phi_1)^{\#}({\mathbf u}^{\text{best}})-(I_{\text{RBF}}\circ\phi_1)^{\#}(P_{\text{Im}(\widehat{\Phi}_{u})} {\mathbf u}^{\text{best}})\rVert_2\leq C\lVert {\mathbf u}^{\text{best}}-P_{\text{Im}(\widehat{\Phi}_{u})} {\mathbf u}^{\text{best}}\rVert_2,
  \end{align*}
  applying the change of variables formula.
\end{proof}

\begin{rmk}[Sources of error]
  We will describe the different sources of error in the right hand side of theorem~\ref{theo:pbdwmsq}. The fist term $\text{trace}(H_{\widehat{\mathbf{u}}_{\text{PBDW}}}SH_{\widehat{\mathbf{u}}_{\text{PBDW}}}^T)$ accounts for the uncertainty $\epsilon_y$ in the observations $y$: for exact observations it is zero.
  The second term $\lVert(\text{Id}-H_l)\circ P_{\text{Im}(H_l)}\rVert_2$ accounts for the additional term $\lVert \eta\rVert_{R^{-1}}$ in equation~\ref{eq:pbdw_hetero}: if this term would be removed than we would have $\lVert(\text{Id}-H_l)\circ P_{\text{Im}(H_l)}\rVert_2=0$ as $H_l$ would be the exact linear projection into $\text{Im}(H_l) = \text{Im}(\widehat{\Phi}_{u})+\left(\text{Im}(\widehat{\Phi}_{u})^{\perp}\cap\text{Im}(\mathcal{Z}_{\mathbf u})\right)$, see also~\cite{gong2019pbdw}.
  The PBDW stability constant $\lVert\text{Id}-H_l\rVert_2$ reduces to $\beta^{-1}(\Phi_{\widehat{\mathbf{u}}},  \mathcal{Z}_{\mathbf u})$ with $\beta(\Phi_{\widehat{\mathbf{u}}},  \mathcal{Z}_{\mathbf u})= \text{inf}_{z\in\text{Im}(\widehat{\Phi}_{u})}\sup_{\eta\in\text{Im}(\mathcal{Z}_{\mathbf u})} \frac{(z,\eta)}{\lVert z\rVert_2\lVert \eta\rVert_2}$, when we remove the regularization term $\lVert \eta\rVert_{R^{-1}}$ in equation~\ref{eq:pbdw_hetero}, obtaining the original noiseless PBDW formulation, see also~\cite{gong2019pbdw}.
  The fourth term:
  \begin{equation*}
    \lVert(\phi_{\text{RBF}})^{\#}({\mathbf u}^{\text{best}})-\widehat{\mathbf{u}}^{\text{true}}-P_{{\text{Im}(H_l)}}((\phi_{\text{RBF}})^{\#}({\mathbf u}^{\text{best}})-\widehat{\mathbf{u}}^{\text{true}})\rVert_2
  \end{equation*}
  represents the approximation error of the velocity fields in the new geometry $\widehat{\mathbf{u}}^{\text{true}}$ with the transported template manifold $(\phi_{\text{RBF}})^{\#}(X^{\text{train}}_{ut})$: it is zero if either $(\phi_{\text{RBF}})^{\#}({\mathbf u}^{\text{best}})=\widehat{\mathbf{u}}^{\text{true}}$, that is the velocity field $\widehat{\mathbf{u}}^{\text{true}}$ is contained in the transported solution manifold, or $\widehat{\mathbf{u}}^{\text{true}}=(\phi_{\text{RBF}})^{\#}({\mathbf u}^{\text{best}})+x$, with $x\in\text{Im}(\widehat{\Phi}_{u})$, if the error committed in approximating $\widehat{\mathbf{u}}^{\text{true}}$ with the transported solution manifold does not affect PBDW.
  The fifth term $C\cdot\lVert {\mathbf u}^{\text{best}}-P_{\text{Im}(\widehat{\Phi}_{u})} {\mathbf u}^{\text{best}}\rVert_2$ includes a multiplicative constant $C>0$, that accounts for the volume deformation from the template to the target geometry, and the reconstruction rSVD error in the template geometry studied in section~\ref{subsec:sml_rec}. The last term
  \begin{equation*}
    P_{\text{Im}(\widehat{\Phi}_{u})^T\cap\text{Im}(\mathcal{Z}_{\mathbf u})}(\phi_{\text{RBF}})^{\#}({\mathbf u}^{\text{best}})+P_{\text{Im}(\widehat{\Phi}_{u})}(\phi_{\text{RBF}})^{\#}({\mathbf u}^{\text{best}})-(\phi_{\text{RBF}})^{\#}(P_{\text{Im}(\Phi_{\mathbf u})}{\mathbf u}^{\text{best}}),
  \end{equation*}
  is a compatibility condition related to the commutativity of the diagram:
  \[ \begin{tikzcd}
    \mathbb{R}^{d_{\mathbf u}} \arrow{r}{P_{\text{Im}(\Phi_{\mathbf u})}} \arrow[swap]{d}{\phi_{\text{RBF}}} & \text{Im}(\Phi_{\mathbf u})\subset\mathbb{R}^{d_{\mathbf u}} \arrow{d}{\phi_{\text{RBF}}} \\%
    \mathbb{R}^{\widehat{d}_{\mathbf u}} \arrow{r}{P_{\text{Im}(\widehat{\Phi}_{u})}}& \text{Im}(\widehat{\Phi}_{u})\subset\mathbb{R}^{\widehat{d}_{\mathbf u}}
  \end{tikzcd}
  \]
\end{rmk}

\begin{rmk}[Registration map error]
  In theorem~\ref{theo:pbdwmsq}, we have not included the RBF interpolation error and the error in the approximation of the exact registration map $\phi_{1}$ with $\phi_{1,N}$, with the notations of theorem~\ref{theo:existreg}: let $u:\Omega_S\subset\mathbb{R}^3\rightarrow\mathbb{R}^3$ be a velocity field on the template computational domain $\Omega_S$:
  \begin{align*}
    \lVert &(I_{\text{RBF}}\circ \phi_{1,N})^{\#}(u)-(\phi_{1})^{\#}(u)\rVert_{L^2(\phi_1(\Omega_S))}\leq\\
    &\lVert (I_{\text{RBF}}\circ \phi_{1,N})^{\#}(u)-(\phi_{1,N})^{\#}(u)\rVert_{L^2(\phi_1(\Omega_S))}+
    \lVert(\phi_{1,N})^{\#}(u)-(\phi_{1})^{\#}(u)\rVert_{L^2(\phi_1(\Omega_S))}\\
    &\leq C h^2_{\phi_1(X_S), \Omega_T}\rho_{\phi_1(X_S), \Omega_T}^{1/2}\lVert(\phi_{1})^{\#}(u)\rVert_{L^2(\phi_1(\Omega_S))} + \sup_{\mathbf{x}\in\Omega_S}\lVert \nabla u(\mathbf{x})\rVert_{L^2(\phi_1(\Omega_S))}\cdot \lVert \phi_{1, N}^{-1}-\phi_1^{-1}\rVert_{L^2(\phi_1(\Omega_S))}\\
    &\leq C h^2_{\phi_1(X_S), \Omega_T}\rho_{\phi_1(X_S), \Omega_T}^{1/2}\lVert(\phi_{1})^{\#}(u)\rVert_{L^2(\phi_1(\Omega_S))} + \sup_{\mathbf{x}\in\Omega_S}\lVert \nabla u(\mathbf{x})\rVert_{L^2(\phi_1(\Omega_S))}\cdot C_3\cdot N^{-\frac{m}{2(d+1)}},\qquad \forall m\geq 1,
  \end{align*}
  where we have applied theorem~\ref{theo:existreg} with $n=0$, and to bound the RBF interpolation error we have applied theorem 4.2~\cite{narcowich2006sobolev} with $\beta=2$, since $(\phi_{1,N})^{\#}(u)\in\mathcal{C}^{2, 1}(I, \mathcal{C}^{2, 1})\subset W^{2,2}(\phi_1(\Omega_S))$, and $\tau=5/2$, since the Fourier transform of thin-plate spline basis~\cite{schaback1999improved} decays as $\lVert\omega\rVert_2^{-5}$ in frequency space $\omega\in\mathbb{R}^3$. The constants introduced are defined as:
  \begin{equation*}
    h_{\phi_1(X_S), \Omega_T} = \sup_{\mathbf{x}\in\Omega_T}\inf_{\mathbf{x}_j\in \phi_1(X_S)}\lVert\mathbf{x}-\mathbf{x}_j\rVert_2,\qquad \rho_{\phi_1(X_S), \Omega_T}=\frac{h_{\phi_1(X_S), \Omega_T}}{\frac{1}{2}\min_{j\neq k}\lVert\mathbf{x}_i-\mathbf{x}_j\rVert_2},
  \end{equation*}
  following the notation in~\cite{narcowich2006sobolev}.
\end{rmk}

\section{Reduced order PPE and STE}
\label{apendix:rom}

The reduced PPE and STE can be defined from their matrix forms in Equations~\eqref{eq:ppe} and~\eqref{eq:ste} respectively:
\begin{align*}
  \widehat{\Phi}_{p}^T A_{\text{PPE}}\,\widehat{\Phi}_{p}z_{\widehat{p}^{n+1/2}_{\text{PPE}}} = \widehat{\Phi}_{p}^T M_{\text{PPE}}^{n+1}\widehat{\ub}_{n+1} - \widehat{\Phi}_{p}^T M_{\text{PPE}}^{n}\widehat{\ub}_{n} + \widehat{\Phi}_{p}^T Q_{\text{PPE}}(\widehat{\ub}^{n+1/2}, \widehat{\ub}^{n+1/2}),
\end{align*}
\begin{align*}
  \widehat{\Phi}_{u}^T A_{\text{STE}}\,\widehat{\Phi}_{u}\mathbf{z}_{\wb} + \widehat{\Phi}_{u}^T B \widehat{\Phi}_{p}z_{\widehat{p}^{n+1/2}_{\text{STE}}} &= \widehat{\Phi}_{u}^T M_{\text{STE}}^{n+1}\widehat{\ub}_{n+1} - \widehat{\Phi}_{u}^T M_{\text{STE}}\widehat{\ub}_{n} + \widehat{\Phi}_{u}^T Q_{\text{STE}}(\widehat{\ub}^{n+1/2}, \widehat{\ub}^{n+1/2}) + \widehat{\Phi}_{u}^T M_{\text{STE}}\widehat{\ub}_{n+1/2},\\
     \widehat{\Phi}_{p}^T B^T\widehat{\Phi}_{u}\mathbf{z}_{\wb} &= 0,
\end{align*}
after left multiplication with the corresponding rSVD transported basis $\widehat{\Phi}_{u}\in\mathbb{R}^{r_{\mathbf u}\times \widehat{d}_{\mathbf u}}$ and $\widehat{\Phi}_{p}\in\mathbb{R}^{r_p\times \widehat{d}_p}$ and subsitution of $\widehat{p}^{n+1/2}_{\text{PPE}}\in\mathbb{R}^{\widehat{d}_p},\widehat{p}^{n+1/2}_{\text{STE}}\in\mathbb{R}^{\widehat{d}_p}, \wb\in\mathbb{R}^{\widehat{d}_{\mathbf u}}$ with the respective reduced variables $z_{\widehat{p}^{n+1/2}_{\text{PPE}}}\in\mathbb{R}^{r_p},z_{\widehat{p}^{n+1/2}_{\text{STE}}}\in\mathbb{R}^{r_p}, z_{\wb}\in\mathbb{R}^{r_{\mathbf u}}$.

The velocity rSVD modes $\widehat{\Phi}_{u}$ are enriched with the supremizer technique~\cite{ballarin2015supremizer} with additional $r_p$ modes computed from the pressure rSVD modes $\widehat{\Phi}_{p}$, for a total of $r_{u,\text{sup}}=r_{\mathbf u}+r_p$ velocity modes $\Phi_{\widehat{\mathbf{u}},\text{sup}}\in\mathbb{R}^{r_{u,\text{sup}}\times \widehat{d}_{\mathbf u}}$:
\begin{align*}
  \Phi_{\widehat{\mathbf{u}},\text{sup}}^T A_{\text{STE}}\,\Phi_{\widehat{\mathbf{u}},\text{sup}}\mathbf{z}_{\wb} &+ \Phi_{\widehat{\mathbf{u}},\text{sup}}^T B \widehat{\Phi}_{p}z_{\widehat{p}^{n+1/2}_{\text{STE}}} = \\
  &=\Phi_{\widehat{\mathbf{u}},\text{sup}}^T M_{\text{STE}}\widehat{\ub}_{n+1} - \Phi_{\widehat{\mathbf{u}},\text{sup}}^T M_{\text{STE}}\widehat{\ub}_{n} + \Phi_{\widehat{\mathbf{u}},\text{sup}}^T Q_{\text{STE}}(\widehat{\ub}^{n+1/2}, \widehat{\ub}^{n+1/2}) + \Phi_{\widehat{\mathbf{u}},\text{sup}}^T M_{\text{STE}}\widehat{\ub}_{n+1/2},\\
    &\qquad\,\,\widehat{\Phi}_{p}^T B^T\Phi_{\widehat{\mathbf{u}},\text{sup}}\mathbf{z}_{\wb} = 0.
\end{align*}
In principle, the supremizer enrichment could be performed offline on the template geometry and then transported with the registration map on the new patient geometry. It may occur that the method is yet not stable and additional corrections to the supremizers must be implemented. This is a future research direction. For the moment, we evaluted the supremizers directly on the test geometries.

\section{Computational costs}
\label{appendix:costs}

In Table~\ref{tab:costs} are reported the computational costs of all the numerical procedures introduced. The average cost of high-fidelity simulations with respect to the total number of training and test geometries is reported. To be efficient, a data assimilation method/surrogate model must achieve a considerable speed-up with respect to high-fidelity simulations, possibly employing less computational resources. This is possible because, in the case of surrogate models, the accuracy of the predictions is reduced; in the data assimilation case, the accuracy is inherently dependent on the noise level. Moreover, the offline costs such as rSVD evaluation, \textit{gnn-gp}, \textit{gnn-gv}, \textit{gnn-vp-pbdw} training, and the supremizer enrichment (SUP), are amortized over possibly infinite online queries.

In principle, rSVD computations are performed on 1 CPU but could be implemented in parallel. Data-parallelism on multiple GPUs could be easily applied to the \textit{gnn-gp}, \textit{gnn-gv}, and \textit{gnn-vp-pbdw} training. Model-parallelism could be useful to extend the EPD-GNN architecture to the whole mesh: we remark that we employ a coarser mesh as support of the EPD-GNN architecture because of memory constraints on the single $40GB$ RAM GPU we employ. The timings for a single time instance prediction with \textit{gnn-gp}, \textit{gnn-gv}, and \textit{gnn-vp-pbdw} are reported. They could also be efficiently parallelized for multiple time instances, as well as RBF interpolation.

Our implementation of RBF interpolation (RBFv, RBFp, RBF) between template and target aorta meshes (or coarse and fine meshes in the EPD-GNNs cases) relies on \texttt{SciPy}~\cite{2020SciPy-NMeth}: we employ a fixed number of $30$ neighbours for each evaluation point. However, a distributed-memory parallel implementation of RBF interpolation could be implemented, partitioning the target aorta mesh. The velocity and pressure modes transported from the template to the target geometries are not orthonormal anymore, so an efficient QR orthonormalization step (QRv, QRp) must be applied online. This step could also be parallelized.

The reduced-order models of the PPE and STE pressure estimators need to load and assemble online the right-hand sides (RHS) coming from the velocity field PBDW predictions. To do this efficiently, we have implemented a parallel version of PPE-ROM and STE-ROM. We have reported the costs for the assembly (ASS) of the RHS along with the costs for a single time instance solve (SOL) of the reduced systems: PPE-ROM reduced systems are of size $r_p\times r_p$, while STE-ROM reduced systems are of size $(2r_p+r_{\mathbf u})\times (2r_p+r_{\mathbf u})$ due to the employment of supremizers. As mentioned in appendix~\ref{apendix:rom}, we compute the supremimzers online, but in principle they could be computed offline.

The PBDW predictions for a single time instance are performed on a single CPU. In our context, the dimension of the observations $(M_{\text{voxels}}+1)\sim30000$ is greater than the number of velocity modes employed $(M_{\text{voxels}}+1)\gg r_{\mathbf u}\in\{500, 1000, 2000\}$: we solve for a generalized least-squares problem of dimension $(M_{\text{voxels}}+1)\times r_{\mathbf u}$ in equation~\eqref{eq:1pbdw} and a linear correction of dimension $(M_{\text{voxels}}+1)\times(M_{\text{voxels}}+1)$, in equation~\eqref{eq:2pbdw}. The PBDW predictions could be sped up reducing the resolution of the observations or with distributed-memory parallelism, partitioning the mesh supported on the voxels centroids. Moreover, the employment of a heteroscedastic noise model with respect to a homoscedastic one increases significantly the computational burden.

\begin{table}[hbtp!]
  \centering
  \caption{Computational costs of the numerical procedures introduced along with computational resources employed. Elapsed real times are shown. Offline stage costs are reported in blue, online stage costs in red. ResNet-LDDMM registration is performed both offline and online.}
  \footnotesize
  \begin{tabular}{|c|}
    \hline
    ResNet-LDDMM registration (1 GPU)\\
    \hline
    \hline
    $\sim$ 10 [min] \\
    \hline
\end{tabular}
\vspace{0.25cm}
\begin{tabular}{|>{\columncolor{blue!20}}c|}
    \hline
    Single high-fidelity simulations (72 CPUs)\\
    \hline
    \hline
    2.64 [h] \\
    \hline
\end{tabular}
\begin{tabular}{| c| >{\columncolor{blue!20}}c| >{\columncolor{blue!20}}c |}
  \hline
  & rSVD $u$ (1 CPU)& rSVD $p$ (1 CPU)\\
  \hline
  \hline
  $r=500$ & 34 [min] & 10 [min]\\
  \hline
  $r=1000$ & 36 [min] & 12 [min]\\
  \hline
  $r=2000$ & 38 [min] & 13 [min]\\
  \hline
  $r=4000$ & 51 [min] & 18 [min]\\
  \hline
\end{tabular}

\vspace{0.25cm}
\begin{tabular}{| c | >{\columncolor{blue!20}}c | >{\columncolor{red!20}}c |>{\columncolor{blue!20}}c | >{\columncolor{red!20}}c |}
  \hline
    & \textit{gnn-gp} train (1 GPU)  & \textit{gnn-gp} predict (1 GPU) & \textit{gnn-gv} train (1 GPU)  & \textit{gnn-gv} predict (1 GPU)\\
  \hline
  \hline
  $e=9,w=10,h=64$ &35 [h]& 0.65 [s] + RBF 5.5 [s] &37 [h]& 0.67 [s] + RBF 5.5 [s]\\
  \hline
  $e=9,w=15,h=128$ &39 [h]& 0.89 [s] + RBF 5.5 [s] &34 [h]& 0.92 [s] + RBF 5.5 [s] \\
  \hline
  $e=9,w=20,h=256$ &43 [h]& 1.21 [s] + RBF 5.5 [s] &43 [h]& 1.29 [s] + RBF 5.5 [s]\\
  \hline
\end{tabular}

\vspace{0.25cm}
\begin{tabular}{| c |>{\columncolor{blue!20}}c | >{\columncolor{red!20}}c |}
  \hline
    &\textit{gnn-vp} train (1 GPU) & \textit{gnn-vp} predict (1 GPU)\\
  \hline
  \hline
  $e=9,w=10,h=64$ &21 [h]& 0.8 [s] + RBF 5.5 [s] \\
  \hline
  $e=9,w=15,h=128$ &40 [h]& 1.08 [s] + RBF 5.5 [s] \\
  \hline
  $e=9,w=20,h=256$ &40 [h]& 1.42 [s] + RBF 5.5 [s]\\
  \hline
\end{tabular}

\vspace{0.25cm}
\begin{tabular}{| >{\columncolor{red!20}}c | >{\columncolor{red!20}}c |}
  \hline
   PPE-FOM  (32 CPUs) & STE-FOM (32 CPUs) \\
   \hline
   \hline
    SOL 1.61 [s] & SOL 6.33 [s]\\
   \hline
\end{tabular}

\vspace{0.25cm}
\begin{tabular}{| c | >{\columncolor{red!20}}c | >{\columncolor{red!20}}c |}
  \hline
   & PPE-ROM  ASS (32 CPUs) & PPE-ROM  SOL+RBF+QR (1 CPU) \\
  \hline
  \hline
  $r_p=500$ & ASS 100 [ms] & SOL 0.006 [s] + RBFp 28 [s] + QRp 4.4 [s]\\
  \hline
  $r_p=1000$ & ASS 112 [ms] & SOL 0.03 [s] + RBFp 38 [s] + QRp 17 [s]\\
  \hline
  $r_p=2000$ & ASS 105 [ms] & SOL 0.34 [s] + RBFp 50 [s] + QRp 46 [s]\\
  \hline
\end{tabular}

\vspace{0.25cm}
\begin{tabular}{| c | >{\columncolor{red!20}}c | >{\columncolor{red!20}}c |}
  \hline
  $r_{\mathbf u}=1000$ Fixed & STE-ROM  ASS+SUP (32 CPUs) & STE-ROM SOL+RBF+QR (1 CPU) \\
  \hline
  \hline
  $r_p=500$ & ASS 111 [ms]+SUP 461 [s] & SOL 0.31 [s]+RBFv 32 [s]+ QRv 14 [s]+ RBFp 28 [s] + QRp 4.4 [s]\\
  \hline
  $r_p=1000$ & ASS 130 [ms]+SUP 1125 [s] & SOL 1.02 [s]+RBFv 63 [s]+ QRv 44 [s]+ RBFp 38 [s] + QRp 17 [s]\\
  \hline
  $r_p=2000$ & ASS 123 [ms]+ SUP 2186 [s] & SOL 3.1 [s]+RBFv 112 [s]+ QRv 143 [s]+ RBFp 50 [s] + QRp 46 [s]\\
  \hline
\end{tabular}

\vspace{0.25cm}
\begin{tabular}{| c | >{\columncolor{red!20}}c |}
  \hline
   & PBDW predict single time instant (1 CPU)\\
  \hline
  \hline
  $r_{\mathbf u}=500$ & PBDW 32.4 [s] + RBFv 32 [s] + QRv 14 [s]\\
  \hline
  $r_{\mathbf u}=1000$ & PBDW 74.2 [s] + RBFv 63 [s] + QRv 44 [s]\\
  \hline
  $r_{\mathbf u}=2000$ & PBDW 74.9 [s] + RBFv 112 [s] + QRv 143 [s]\\
  \hline
\end{tabular}
  \label{tab:costs}
\end{table}

\printbibliography


\end{document}